\documentclass[journal, 12pt]{IEEEtran}
\usepackage{booktabs,multirow}
\usepackage{amsmath,amssymb,amsfonts}
\usepackage{graphicx}          
\usepackage{cite}
\usepackage{epstopdf}

\ifCLASSINFOpdf
\else

\fi

\usepackage{multicol}

\usepackage{subfigure}
\usepackage{amsmath}
\usepackage{amssymb}
\usepackage{amsfonts}
\usepackage{graphicx}
\usepackage{cite}
\usepackage{url}
\usepackage{color}
\usepackage{ccaption}
\usepackage{booktabs} 
\usepackage{graphics}
\usepackage{booktabs,multirow}
\usepackage{amsmath,amssymb,amsfonts}
\usepackage{graphicx}          
\usepackage{cite}
\usepackage{epstopdf}
\usepackage{color}
\usepackage{graphicx}
\usepackage{diagbox}        
\usepackage{amsmath}
\usepackage{amssymb}
\usepackage{amsmath}
\usepackage{amsfonts}
\usepackage{cite}
\usepackage{picins}

\newtheorem{assumption}{\bf{Assumption}}
\newtheorem{remark}{\bf{Remark}}
\newtheorem{theorem}{\bf{Theorem}}
\newtheorem{lemma}{\bf{Lemma}}
\newtheorem{corollary}{\bf{Corollary}}
\newtheorem{definition}{\bf{Definition}}
\newtheorem{proposition}{\bf{Proposition}}

\begin{document}
%
\title{Scheduling problems of aircraft on a same runway and dual runways}
\author{Peng Lin, Haopeng Yang, Gui Gui, Mengxiang Zeng and Weihua Gui
	\thanks{Peng Lin, Haopeng Yang, Gui Gui, Mengxiang Zeng and Weihua Gui are with the School of Automation, Central South University, Changsha 410083, China (email: lin\_peng0103@sohu.com, lin\_peng@csu.edu.cn).}
}
\maketitle

\begin{abstract}
In this paper, the scheduling problems of landing and takeoff aircraft on a same runway and on dual runways are addressed. In contrast to the approaches based on mixed-integer optimization models in existing works, our approach focuses on the minimum separation times between aircraft by introducing some necessary assumptions and new concepts including relevance, breakpoint aircraft, path and class-monotonically-decreasing sequence. Four scheduling problems are discussed including landing scheduling problem, takeoff scheduling problem, and mixed landing and takeoff scheduling problems on a same runway and on dual runways with the consideration of conversions between different aircraft sequences in typical scenarios. Two real-time optimal algorithms are proposed for the four scheduling problems by fully exploiting the combinations of different classes of aircraft, and necessary definitions, lemmas and theorems are presented for the optimal convergence of the algorithms. Numerical examples are presented to show the effectiveness of the proposed algorithms. In particular, when $100$ aircraft are considered, by using the algorithm in this paper, the optimal solution can be obtained in less than $5$ seconds, while by using the CPLEX software to solve the mix-integer optimization model, the optimal solution cannot be obtained within $1$ hour.
\end{abstract}


\section{Introduction}
In recent years, the issue of insufficient airport capacity has become increasingly prominent with the rising demand for air transportation, leading to more severe traffic congestion and aircraft delays. Though more runways might be an effective way to relieve this issue, most of the airports have no capacity to construct new runways due to constraints including cost, terrain, and surrounding environmental factors. For this reason, researchers have turned their attention to rescheduling aircraft sequences to improve runway utilization, reduce aircraft delays, and enhance on-time performance.

The core objective of rescheduling aircraft sequences is to reassign each aircraft a runway and scheduled takeoff or landing time to minimize a given objective function, while meeting both time window constraints and wake turbulence separation requirements, when the on-time performance of the aircraft cannot be guaranteed or even when there are conflicts in the flight plans of different airlines.

According to the studied objective functions, existing works on this topic can be mainly categorized into four optimization problems: the minimization problem of the total delay time \cite{b1,b2,b3,b4,b5,b6,b7,b8,b9,b10,b11,b12,b13}, the minimization problem of the total deviation of scheduled takeoff and landing times \cite{b14,b15,b16,b17,b18,b19,b20,b21}, the minimization problem of {the total makespan} for arrival and departure aircraft \cite{b22,b23,b24,b25,b26}, and the minimization problem of overall operational costs \cite{b27,b28,b29}. In existing works, the problems of rescheduling aircraft sequences were usually modeled as a mixed-integer optimization problem, which is essentially a NP-hard problem. Though the mixed-integer optimization models have wide applications and can find the global optimization solutions, the high computational complexity and the high computing time cost make the corresponding algorithms hard to completely match with desired performances.

Based on the mixed-integer optimization models, there are mainly four kinds of algorithms: mixed-integer programming algorithms, dynamic programming algorithms, heuristic algorithms, and metaheuristic algorithms for the aircraft sequencing optimization problem. Mixed-integer programming algorithm is essentially a searching algorithm in the whole space whose computational complexity is usually the highest one among the four kinds of algorithms. In contrast to mixed-integer programming algorithm, where the global optimal solution can be obtained, dynamic programming algorithms, heuristic and meta-heuristic algorithms can offer lower computation complexity but cannot find the global optimal solution and the optimal errors cannot be obtained as well in general due to the limitations of different algorithms. Specifically, for dynamic programming algorithms, some additional assumptions are imposed which can significantly reduce the computational complexity but meanwhile deviate the optimal solutions of the system by unknown errors for  dynamic programming algorithms, while for heuristic and meta-heuristic algorithms, the main idea is to construct initial solutions under given constraints, and then search within the neighborhood of these solutions by iteratively using different operations including transformations and variables exchange without considering the algorithm optimal convergence from a global view point.

In this paper, the scheduling problems of landing and takeoff aircraft on a same runway and dual runways are addressed. In most of existing works, aircraft were considered based on the ICAO separations, which includes $3$ classes of aircraft and there are totally $3\times 3=9$ kinds of pairwise combinations of landing and takeoff aircraft. In this paper, our work can be used to the RECAT-EU systems, which might include $6$ or more classes of aircraft and there might be totally no smaller than $6\times 6=36$ kinds of pairwise combinations of landing and takeoff aircraft.
 The scheduling problems studied in this paper are much more comprehensive and complicated than the existing works based on ICAO separations. In contrast to the approaches based on mixed-integer optimization models in existing works, our approach focuses on the minimum separation times between aircraft by introducing some necessary assumptions and new concepts including relevance, breakpoint aircraft, path and class-monotonically-decreasing sequence. Four scheduling problems are discussed including landing scheduling problem, takeoff scheduling problem, and mixed landing and takeoff scheduling problems on a same runway and on dual runways. Two real-time optimal algorithms are proposed to find the optimal aircraft sequences to minimize the given objective function.
Numerical examples are presented to show the effectiveness of the algorithms. In particular, when $100$ aircraft are considered, by using the algorithms in this paper, the optimal solution can be obtained in less than $5$ seconds, while by using the CPLEX software to solve the mix-integer optimization model, the optimal solution cannot be obtained within $1$ hour.

The main contributions of this paper are mainly in four aspects.
\begin{itemize}
  \item The first contribution is that this paper studies the interaction mechanism of the minimum separation times between aircraft in sequences and establishes a new theoretical framework for scheduling problems of aircraft, which is completely different from the framework of mixed-integer optimization problem.
  \item The second contribution is that optimal solutions can be obtained for scheduling problems of aircraft by excluding most of the non-optimal solutions and narrowing down the search for the optimal solutions within a small range. This is also different from existing works, where optimal solutions can rarely be obtained and the final optimal errors are usually unknown.
  \item The third contribution is that the proposed algorithms exhibit the features of polynomial algorithms and can be applied in actual systems in real time. In existing real-time works, the related algorithms can only be applied under some additional artificial conditions and the resulting optimal errors are unknown.
  \item The fourth contribution is that our work can be used for the RECAT-EU systems, which might include $6$ or more classes of aircraft. In existing works based on the ICAO separations, aircraft are usually classified into $3$ classes. The scheduling problems studied in this paper are much more comprehensive and complicated than the existing works based on ICAO separations.
\end{itemize}

Notations. The operation $A-B$ represents the set that consists of the elements of $A$ which are not elements of $B$; the operation $\phi_0-\phi_1$ represents the sequence which is obtained through modifying the sequence $\phi_0$ by removing the aircraft belonging to the sequence $\phi_1$ and keeping the orders of the rest aircraft unchanged; and the symbol $/$ represents the meaning of ``or".

\section{Problem formulation}
Without considering other airport constraints, in order to ensure the safety of the aircraft, the minimum separation time between each aircraft and its leading aircraft is only related to their own classes.

Suppose that aircraft can be partitioned into $\eta$ classes in descending order of wake impact, represented as $\mathcal{I}=\{1,2, \cdots, \eta \}$, where $\eta$ is a positive integer, and in general the class $1$ usually refers to $\mathrm{A380}$ aircraft.

Let $F(\phi, Sr(\phi))=S_n(\phi)-t_0$ denote the time for a given aircraft sequence $\phi=\langle Tcf_1, Tcf_2, \cdots, Tcf_n\rangle$ to complete all operation (landing and takeoff) tasks, i.e., the total makespan, where $t_0$ denotes the initial operation time, and $Sr(\phi)=(S_1(\phi), S_2(\phi), \cdots, S_n(\phi))$ represents the rescheduled takeoff and landing times of aircraft $Tcf_1, Tcf_2, \cdots, Tcf_n$ with $t_0\leq S_1(\phi)\leq S_2(\phi)\leq \cdots \leq S_n(\phi)$.

The purpose of this paper is to find appropriate operation sequence of aircraft and the corresponding takeoff and landing times to minimize the objective function $F(\phi, Sr(\phi))$ so as to solve the following optimization problem,
\begin{eqnarray}\label{optim1}\begin{array}{lll}\min F(\phi, Sr(\phi))\\
\mbox{Subject to}~~~S_k \in Tf_k=[f_k^{\min}, f_k^{\max}],\\
\hspace {2cm}k =1, \cdots, n,\\
\hspace {2cm}S_{j}-S_i \geq Y_{{i}{j}},~\forall i<j, \\ \hspace{2cm}i ,j=1, \cdots, n,\end{array}\end{eqnarray}
where $Tf_k=[f_k^{\min}, f_k^{\max}]$ represents the set of allowable takeoff or landing times, i.e., the time window constraint, for the aircraft $Tcf_k$ for two constant $f_k^{\min}\leq f_k^{\max}$, $Y_{ij}$ represents the minimum separation time between an aircraft $Tcf_i$ and its trailing aircraft $Tcf_{j}$. For the sake of expression convenience, when no confusion arises, $F(\phi, Sr(\phi))$, $Sr(\phi)$ and
$S_i(\phi)$ can be abbreviated as $F(\phi)$, $Sr$ and $S_i$.

\begin{definition}\label{definition3}{\rm (Relevance) Consider an aircraft sequence $\phi=\langle Tcf_1, Tcf_2, \cdots, Tcf_n\rangle$. Let $S_{ij}=S_j-S_i$ represent the operation time interval between aircraft $Tcf_i$ and aircraft $Tcf_j$ for $i<j$. If $S_{ij}=Y_{{i}j}$ for $i<j$, it is said that aircraft $Tcf_j$ is relevant to the aircraft $Tcf_i$.}\end{definition}

In this paper, we will discuss the scheduling problems of aircraft under four typical scenarios including the landing/takeoff scheduling problem and the mixed landing and takeoff scheduling problems on a same runway and dual runways, which can be extended to the more general scheduling problems based on the analysis approach in this paper.

\section{Landing scheduling problem on a single runway}\label{seclanding}

Let $T_{ij}$ represent the minimum separation time between a landing aircraft of class $i$ and a trailing landing aircraft of class $j$ on a single runway without considering the influence of other aircraft. Let $T_0$ denote the minimum value of all possible separation times between aircraft, where $T_0$ is usually taken as $1$ minute.

Based on the minimum landing separation time standards at Heathrow Airport \cite{b30} and the understanding of the physical landing process, we propose the following assumptions and definitions.

\begin{assumption}\label{ass1.3.1}{\rm (1) For $i=1,2$, $T_{ii}=1.5T_0$.

(2) For $i=\rho_1, \rho_2$, $T_{ii}=T_0+\delta$, where $T_0/8<\delta<T_0/6$ is a positive number, $\rho_1=3$ and $\rho_2=5$.

(3) For $i\neq 1,2,\rho_1, \rho_2$, $T_{ii}=T_0$.}
\end{assumption}

        \begin{assumption}\label{ass1.3.3}{\rm
(1) For all $i<k<j$, $T_{ik}\leq T_{ij}\leq 3T_0$ and $T_{kj}<T_{ij}\leq 3T_0$.

(2)  For all $k\leq j\leq i$, $T_{ik}<T_{ij}+T_{jk}$ and $T_{ki}<T_{ji}+T_{kj}$.}\end{assumption}

\begin{assumption}\label{ass1.3.2}{\rm (1) $T_{21}=1.5T_0$.

(2) For $i>j$ and $i\neq 2$, $T_{ij}=T_0$.}\end{assumption}


\begin{remark}{\rm Though different airports might have different parameters, the analysis idea in this paper might be still valid except some special requirements of the airport.}\end{remark}

\begin{remark}{\rm For the wake impact caused by aircraft, the closer the class of the leading aircraft is, the smaller the difference in wake impact is, and the greater the difference in the leading aircraft class is, the greater the difference in wake impact is.}\end{remark}

\begin{remark}{\rm When observing the landing separation time standards at Heathrow Airport (See Table \ref{tab:Heathrow_separation} in Sec. \ref{simulations}), it was found that there are aircraft of adjacent classes with similar wake effects, such as medium and light medium aircraft.}\end{remark}

\begin{lemma}\label{lemma11s}{\rm For all $i,j, k\in \mathcal{I}$, $T_{ik}<T_{ij}+T_{jk}$.}\end{lemma}

\noindent{Proof:} When $i\geq k$, from Assumptions \ref{ass1.3.1} and \ref{ass1.3.2}, $T_{ik}\leq 1.5T_0<2T_0\leq T_{ij}+T_{jk}$. When $i\leq k \leq j$, from Assumption \ref{ass1.3.3}, $T_{ik}\leq T_{ij}<T_{ij}+T_{jk}$. When $j\leq i\leq k$, $T_{ik}\leq T_{jk}<  T_{ij}+T_{jk}$. When $i\leq j\leq k$, from Assumption \ref{ass1.3.3}, $T_{ik}<  T_{ij}+T_{jk}$.

\begin{definition}\label{definition1}{\rm (Breakpoint aircraft) Consider a landing/takeoff sequence $\phi=\langle Tcf_1, Tcf_2, \cdots, Tcf_{n}\rangle$. If the classes of two consecutive aircraft satisfy that $cl_{i}<cl_{i+1}$ where $cl_i$ and $cl_{i+1}$ represent the classes of aircraft $Tcf_i$ and $Tcf_{i+1}$, it is said that the aircraft $Tcf_i$ is a breakpoint aircraft of $\phi$.}\end{definition}

\begin{definition}\label{definition2}{\rm (Resident-point aircraft)  Consider a landing/takeoff sequence $\phi=\langle Tcf_1, Tcf_2, \cdots, Tcf_{n}\rangle$. If $S_1>t_0$, $Tcf_1$ is called  a resident-point aircraft of $\phi$, and $S_1-t_0$ is called the resident time of $Tcf_1$. If the aircraft $Tcf_i$ is not relevant to the aircraft $Tcf_{i-1}$, i.e., $S_{(i-1)i}-Y_{(i-1)i}>0$, for $i=2, 3, \cdots, n$,
it is said that $Tcf_i$ is a resident-point of $\phi$ and $S_{(i-1)i}-Y_{(i-1)i}$ is the resident time of $Tcf_i$.
 Consider a mixed landing and takeoff sequence $\phi=\langle Tcf_1, Tcf_2, \cdots, Tcf_{n}\rangle$. Let $\mu_1$ and $\mu_2$ be the largest integers smaller than $i$ for a landing aircraft and a takeoff aircraft $Tcf_{\mu_2}$. If the aircraft $Tcf_i$ is not relevant to the aircraft $Tcf_{\mu_1}$ and $Tcf_{\mu_2}$, $i=3, 4, \cdots, n$, it is said that $Tcf_i$ is a resident-point aircraft of $\phi$ and $\min\{S_{\mu_1i}(\phi)-Y_{\mu_1{i}}, S_{\mu_2i}(\phi)-Y_{\mu_2{i}}\}$ is the resident time of $Tcf_i$.}\end{definition}

\begin{remark}{\rm Breakpoint aircraft and resident-point aircraft often exist in the actual aircraft sequence. It should be noted that an aircraft $Tcf_i$ in a given sequence $\phi$ might be a breakpoint and a resident point at the same time.}\end{remark}

From the above definitions, it can be seen that the existence of breakpoint aircraft and resident-point aircraft might increase the value of the objective function $F(\cdot)$.

\begin{assumption}\label{ass1.3.4}{\rm (1) For $k=\rho_2$, $T_{(k-1)k}=T_0+\delta$. For $k\neq \rho_2$, $T_{(k-1)k}\geq 1.5T_0$.

(2) For all $i\leq k$, when $(i,k)\neq (\rho_2, \rho_2)$, $T_{(i-1)k}-T_{ik}>2\delta$.

(3) For $k=1$, $T_{k(k+1)}>2T_0,$ and for $k=2$, $T_{k(k+1)}>1.5T_0+2\delta$.

(4) For $k=1$, all $k+2\leq h \leq \eta$, and all $h\leq j\leq \eta$, $T_{kj}-T_{hj}>0.5T_0$.

(5) Let $E=\{\langle3,4\rangle, \langle 3,\eta\rangle\}$ be a sequence set such that $0.5T_0-\delta\leq T_{2j}-T_{kj}<0.5T_0$ for all $\langle k,j\rangle\in E$ and for all $\langle k,j\rangle\notin E$ with $2<k< j$, $T_{2j}-T_{kj}>0.5T_0$.}\end{assumption}

\begin{remark}{\rm Assumption \ref{ass1.3.4}(5) considers a typical scenario for the breakpoints which might yield local minimum points for the objective function $F(\cdot)$. More general scenarios can be studied based on class-sequence sets including $\Psi_0$, $\Psi_1$, $\cdots$, $\Psi_5$ defined later.}\end{remark}

In Theorem \ref{theorem1.1}, we study the aircraft relevance and show that when the landing times of the aircraft have no constraints, the occurrence of resident-point aircraft should be avoided to ensure the optimality of the sequence.

\begin{theorem}\label{theorem1.1}{\rm Consider the sequence of landing aircraft $\phi=\langle Tcf_1, Tcf_2, \cdots, Tcf_{n}\rangle$. Suppose that $Tf_k=[t_0, +\infty]$ for all $k$. Under Assumptions \ref{ass1.3.1}-\ref{ass1.3.4}, the following statements hold.

(1) For all $i=3, 4, \cdots, n$, the aircraft $Tcf_{i}$ is not relevant to $Tcf_1$, $Tcf_2, \cdots,$ $Tcf_{i-2}$.

(2) If the aircraft $Tcf_j$ is relevant to $Tcf_i$, then $j=i+1$.

(3) Suppose that the orders of all landing aircraft in $\phi$ are fixed.  
The optimization problem (\ref{optim1}) is solved if and only if the aircraft $Tcf_{i+1}$ is relevant to the aircraft $Tcf_i$, $i=1,2, \cdots, n-1$.}\end{theorem}
\noindent {Proof}: If $k<i-3$, $S_{ki}>3T_0\geq T_{cl_kcl_i}$ and hence the aircraft $Tcf_i$ is not relevant to $Tcf_k$. If $k=i-3$ or $k=i-2$, from Assumption \ref{ass1.3.3},
$S_{ki}=S_{k(k-1)}+S_{(k-1)i}\geq T_{cl_kcl_{k-1}}+T_{cl_{k-1}cl_i}> T_{cl_kcl_i}$. Therefore, if $k=i-3$ or $k=i-2$, the aircraft $Tcf_i$ is not relevant to the aircraft $Tcf_k$. From the statement (1), the statements (2) and (3) naturally hold.

\begin{assumption}\label{ass1.4}{\rm Consider a landing aircraft sequence $\phi=\langle Tcf_1, Tcf_2, \cdots, Tcf_{n}\rangle$. Suppose that $Tf_k=[t_0, +\infty]$ for all $k$, and for each pair of adjacent aircraft in landing sequence $\phi$, each aircraft is relevant to its leading aircraft.}\end{assumption}

In the following, we first give a calculation method for the objective function $F(\cdot)$ when the order of some aircraft is changed.

\begin{lemma}\label{lemmap131}{\rm Consider a landing aircraft sequence $\Phi_a=\langle\phi_1, Tcf_{s_2}, Tcf_{s_3},\phi_2\rangle$, where $\phi_1=\langle Tcf_1, Tcf_2, \cdots, Tcf_{s_1}\rangle$, $s_2=s_1+1$, $s_3=s_1+2$, $\phi_2=\langle Tcf_{s_3+1}, \cdots, Tcf_{n}\rangle$,  and $s_1, s_2, s_3$ are three positive integers. Suppose that Assumptions \ref{ass1.3.1}-\ref{ass1.4} hold for each landing aircraft sequence.

(1) Move $Tcf_{s_2}$ to be between $Tcf_{h_1}$ and $Tcf_{h_1+1}$ and convert $\phi_2$ into a new sequence $\phi_3$, where $Tcf_{h_1}, Tcf_{h_1+1}\in \phi_2$.
Let $\Phi_b=\langle \phi_1, Tcf_{s_3}, \phi_3\rangle$. It follows that $F(\Phi_a)-F(\Phi_b)=T_{cl_{s_1}cl_{s_2}}+T_{cl_{s_2}cl_{s_3}}+T_{cl_{h_1}cl_{h_1+1}}-T_{cl_{s_1}cl_{s_3}}-T_{cl_{h_1}cl_{s_2}}-T_{cl_{s_2}cl_{h_1+1}}$.

(2) Move $Tcf_{s_3}$ to be between $Tcf_{h_2}$ and $Tcf_{h_2+1}$ and convert $\phi_1$ into a new sequence $\phi_4$, where $Tcf_{h_2}, Tcf_{h_2+1}\in \phi_1$.
Let $\Phi_c=\langle \phi_4, Tcf_{s_2}, \phi_2\rangle$. It follows that $F(\Phi_a)-F(\Phi_c)=T_{cl_{s_2}cl_{s_3}}+T_{cl_{s_3}cl_{s_3+1}}+T_{cl_{h_2}cl_{h_2+1}}-T_{cl_{s_2}cl_{s_3+1}}-T_{cl_{h_2}cl_{s_3}}-T_{cl_{s_3}cl_{h_2+1}}$.
}\end{lemma}

When the aircraft are free of time window constraints, it is assumed by default that the operation time of the first aircraft of a sequence is $t_0$ if no otherwise specified.

Based on the calculation method given in Lemma \ref{lemmap131}, we give Lemma \ref{lemma1.2} to discuss the sum of the separation times between the aircraft of certain class and their trailing aircraft is discussed when the sequence containing a breakpoint aircraft (two class-monotonically-decreasing sequences) is merged into one class-monotonically-decreasing sequence. The sequences containing multiple breakpoint aircraft can be discussed in a similar way or by frequent use of Lemma \ref{lemma1.2}.

Let $f_{k_i}(\phi)$ represent the sum of the separation times between all the aircraft of class $k_i$ and their trailing aircraft in $\phi$. 
When an aircraft $Tcf_k$ is the end of the aircraft sequence, we use $T_{cl_{k}0}=0$ to denote that the aircraft $Tcf_k$ has no trailing aircraft.

\begin{lemma}{\rm \label{lemma1.2} Consider a landing aircraft sequence $\Phi_a=\langle \phi_1, \phi_2\rangle$, where $\phi_1=\langle Tcf_{1}, Tcf_{2}, \cdots, Tcf_{s_1}\rangle$, $\phi_2=\langle Tcf_{s_1+1}, Tcf_{s_1+2},\cdots, Tcf_{n}\rangle$, $cl_{1}\geq cl_{2}\geq \cdots \geq cl_{s_1}=k_c$, $cl_{s_1+1}\geq cl_{s_1+2}\geq \cdots \geq cl_{n}$, $s_1$ and $k_c$ are positive integers. Merge $\phi_1, \phi_2$ to form a class-monotonically-decreasing sequence $\Phi_b$. Suppose that Assumptions \ref{ass1.3.1}-\ref{ass1.4} hold true for each landing aircraft sequence, $\Theta_a=\{k_{1a},k_{2a},\cdots,k_{\eta_ca}\}$ is the set of all possible aircraft classes in $\phi_1$, and $\Theta_b=\{k_{1b},k_{2b},\cdots,k_{\eta_fb}\}$ is the set of all possible aircraft classes in $\phi_2$ for two positive integers $\eta_c$ and $\eta_f$, where $k_{1a}<k_{2a}<\cdots<k_{\eta_ca}\leq \eta$, $k_{1b}<k_{2b}<\cdots<k_{\eta_fb}\leq \eta$. The following statements hold.

(1) If $k_c\neq h_0=k_{ia}\in \Theta_a$ and $h_0=k_{jb}\in \Theta_b$, $f_{h_0}(\Phi_b)-f_{h_0}(\Phi_a)=T_{h_0h_0}+T_{h_0h_b}-T_{h_0k_{(i-1)a}}-T_{h_0k_{(j-1)b}}$, where $h_b=\max\{k_{(i-1)a}, k_{(j-1)b}\}$.

(2) If $k_c\neq h_0=k_{ia}\in \Theta_a$ and $h_0\notin \Theta_b$, $f_{h_0}(\Phi_b)-f_{h_0}(\Phi_a)=T_{h_0h_b}-T_{h_0k_{(i-1)a}}$, where $h_b$ is the largest integer smaller than $h_0$ in $\Theta_a\cup \Theta_b$.

(3) If $k_c\neq h_0=k_{jb}\in \Theta_b$ and $h_0\notin \Theta_a$, $f_{h_0}(\Phi_b)-f_{h_0}(\Phi_a)=T_{h_0h_b}-T_{h_0k_{(j-1)b}}$, where $h_b$ is the largest integer smaller than $h_0$ in $\Theta_a\cup \Theta_b$.}\end{lemma}

\begin{remark}{\rm It should be noted that in Lemma \ref{lemma1.2}, it is assumed that when $i-1=0$ or $j-1=0$, $T_{h_0(i-1)}=T_{h_00}=0$ and $T_{h_0(j-1)}=T_{h_00}=0$.}\end{remark}

\begin{corollary}{\rm Under the conditions in Lemma \ref{lemma1.2}, the following statements hold.

(1.1) If $k_c\neq h_0=2\in \Theta_a \cap\Theta_b$, and $1\in \Theta_a \cap\Theta_b$, $f_{h_0}(\Phi_b)-f_{h_0}(\Phi_a)=0$.

(1.2) If $k_c\neq h_0=2\in \Theta_a \cap\Theta_b$, and $1\notin \Theta_a \cap\Theta_b$, $f_{h_0}(\Phi_b)-f_{h_0}(\Phi_a)=1.5T_0$.

(2.1) If $k_c\neq h_0$, $2<h_0\in \Theta_a \cap\Theta_b$, and $\Theta_b-\{h_0, h_0+1, \cdots, \eta\}\neq \emptyset$, $f_{h_0}(\Phi_b)-f_{h_0}(\Phi_a)=T_{h_0h_0}-T_0$.

(2.2) If $k_c\neq h_0$, $2<h_0\in \Theta_a \cap\Theta_b$, and $\Theta_b-\{h_0, h_0+1, \cdots, \eta\}=\emptyset$, $f_{h_0}(\Phi_b)-f_{h_0}(\Phi_a)=T_{h_0h_0}$.}\end{corollary}

\begin{lemma}{\rm Consider two landing aircraft sequences $\phi_0=\langle Tcf_1, Tcf_2, Tcf_3, Tcf_4\rangle$ and $\phi_1=\langle Tcf_1, Tcf_3, Tcf_2,$ $ Tcf_4\rangle$. Suppose that Assumptions \ref{ass1.3.1}-\ref{ass1.4} hold for each aircraft sequence and $cl_1=cl_2=2$. If $(cl_3, cl_4)\in E$, $0<F(\phi_0)-F(\phi_1)<\delta$. If $(cl_3, cl_4)\notin E$, $F(\phi_0)-F(\phi_1)<0$.}\end{lemma}

\noindent{Proof:} By simple calculations, this lemma can be proved and hence its proof is omitted.
\begin{remark}{\rm This lemma gives examples to show how to deal with the scenarios in Assumption \ref{ass1.3.4}(5).}\end{remark}

In Lemmas \ref{lemmap1} and \ref{lemma1.31}, we discuss some typical cases for breakpoint aircraft which yields local minimum points. In Theorem \ref{lemma1.3s}, we give a rule to calculate the objective function $F(\cdot)$ based on a standard class-monotonically-decreasing sequence. In Theorem \ref{theorem1.2}, we gives a condition for the optimality of the aircraft sequence.

\begin{definition}{\rm (Class-monotonically-decreasing sequence) If a landing/takeoff sequence $\phi=\langle Tcf_1, Tcf_2, \cdots, Tcf_n\rangle$ satisfies $cl_1\geq cl_2\geq \cdots\geq cl_n$, then the aircraft sequence $\phi$ is called a class-monotonically-decreasing sequence.
}\end{definition}

\begin{lemma}\label{lemmap1}{\rm Consider a landing aircraft sequence $\Phi_a=\langle\phi_1, \phi_2\rangle$, where $\phi_1=\langle Tcf_1, Tcf_2, \cdots, Tcf_{s_1}\rangle$, $\phi_2=\langle Tcf_{s_1+1}, Tcf_{s_1+2}, \cdots, Tcf_{n}\rangle$, both of the sequences $\phi_1$ and $\phi_2$ are class-monotonically-decreasing sequences, $cl_{s_1}<cl_{s_1+1}$, and $s_1$ is a positive integer. Suppose that $\phi_2$ contains $s_2\geq1$ aircraft of class $h_1$, $Tcf_i\in \phi_2$, $cl_i=h_1<cl_{s_1}$, and $Tcf_{s_1+1}\neq Tcf_i$. Suppose that Assumptions \ref{ass1.3.1}-\ref{ass1.4} hold for each landing aircraft sequence. Generate a new sequence $\Phi_b$ by moving the aircraft $Tcf_i$ to be between $Tcf_{s_1}$ and $Tcf_{s_1+1}$ in $\Phi_a$. The following statements hold.

(1) If $h_1=1$, $F(\Phi_a)-F(\Phi_b)<0$.

(2) If $h_1=2$, $s_2>1$ and $\langle cl_{s_1}, cl_{s_1+1}\rangle\in E$, $0<F(\Phi_a)-F(\Phi_b)\leq \delta$.

(3) If $h_1=2$ and $\langle cl_{s_1}, cl_{s_1+1}\rangle\notin E$, $F(\Phi_a)-F(\Phi_b)<0$.

(4) If $h_1=2$, $s_2=1$ and $\langle cl_{s_1}, cl_{s_1+1}\rangle\in E$, $F(\Phi_a)-F(\Phi_b)<0$.

(5) If $h_1>2$, $F(\Phi_a)-F(\Phi_b)<0$.}\end{lemma}

\noindent{Proof}: (1) Since $h_1=1$, $cl_{s_1}\leq 2$. If $cl_{s_1}>2$ and the number of the aircraft of class no larger than 2 in $\phi_2$ is larger than $1$, from Assumption \ref{ass1.3.4}, $T_{h_1cl_{s_1+1}}-T_{cl_{s_1}cl_{s_1+1}}>0.5T_0$.
It follows from Assumptions \ref{ass1.3.1} and \ref{ass1.3.2} that \begin{eqnarray*}\begin{array}{lll}&F(\Phi_a)-F(\Phi_b)\\
&=T_{cl_{s_1}cl_{s_1+1}}-T_{h_1cl_{s_1+1}}-T_{cl_{s_1}h_1}\\&\quad+T_{h_1cl_{i+1}}+T_{cl_{i-1}h_1}-T_{cl_{i-1}cl_{i+1}}\\
&<-0.5T_0-T_0+1.5T_0<0,\end{array}\end{eqnarray*}
where $T_{cl_{i-1}h_1}=T_{cl_{i-1}cl_{i+1}}$ since $cl_{i+1}\leq h_1=1$.
If $cl_{s_1}=2$, noting that $Tcf_{s_1+1}\neq Tcf_i$, by a similar approach,
it can also be obtained that $F(\Phi_a)-F(\Phi_b)<0$.

(2) Since $s_2>1$, either $cl_{i-1}$ or $cl_{i+1}$ is equal to $2$.
We assume that $cl_{i+1}=2$ and the case of $cl_{i-1}=2$ can be discussed in the same way.
Since $h_1=2$ and $\langle cl_{s_1}, cl_{s_1+1}\rangle\in E$,
 \begin{eqnarray*}\begin{array}{lll}&F(\Phi_a)-F(\Phi_b)\\
&=T_{cl_{s_1}cl_{s_1+1}}-T_{h_1cl_{s_1+1}}-T_{cl_{s_1}h_1}+T_{h_1cl_{i+1}}\\&\quad+T_{cl_{i-1}h_1}-T_{cl_{i-1}cl_{i+1}}\\
&=T_{cl_{s_1}cl_{s_1+1}}-T_{h_1cl_{s_1+1}}-T_{cl_{s_1}h_1}+T_{h_1cl_{i+1}}\\
&=T_{cl_{s_1}cl_{s_1+1}}-T_{h_1cl_{s_1+1}}+0.5T_0,\end{array}\end{eqnarray*}
where $T_{cl_{s_1}h_1}=T_0$ since $h_1<cl_{s_1}$ and $T_{h_1cl_{i+1}}=1.5T_0$ since $h_1=cl_{i+1}=2$.
Note that $0.5T_0- \delta\leq T_{h_1cl_{s_1+1}} -T_{cl_{s_1}cl_{s_1+1}}$ $<0.5T_0$ since $\langle cl_{s_1}, cl_{s_1+1}\rangle\in E$. Hence, $0<F(\Phi_a)-F(\Phi_b)\leq \delta$.

(3)-(5) The proofs of statements (3)-(5) can be proved by a similar approach to the proofs of statements (1) and (2).

\begin{remark}{\rm In Lemma \ref{lemmap1}, we consider the transition from the sequence $\Phi_a$ to the sequence $\Phi_b$. Actually, the results in Lemma \ref{lemmap1} still hold for the transition from the sequence $\Phi_b$ to the sequence $\Phi_a$, which can be used to decrease the value of the objective function $F(\cdot)$.
}\end{remark}

\begin{remark}{\rm Lemma \ref{lemmap1}(2) shows the sum of the separation times might be reduced when consecutive aircraft of class $2$ are splitted to be in two class-monotonically-decreasing sequences.}\end{remark}

\begin{lemma}\label{lemma1.31}{\rm  Consider a landing sequence $\Phi_a=\langle\phi_1, \phi_2\rangle$, where $\phi_1=\langle Tcf_1, Tcf_2, \cdots, Tcf_{n_1}\rangle$, $\phi_2=\langle Tcf_{n_1+1},$ $ Tcf_{{n_1}+2}, \cdots, Tcf_{n_2}\rangle$, $\phi_1$ and $\phi_2$ are both class-monotonically-decreasing sequences, $cl_{n_1}<cl_{n_1+1}$, and
$n_1, n_2$ are two positive integers. Suppose that Assumptions \ref{ass1.3.1}-\ref{ass1.4} hold for each landing aircraft sequence. Merge the aircraft sequences $\phi_1$ and $\phi_2$ to form a new class-monotonically-decreasing sequence $\Phi_b=\langle Tcf_{s_1}, Tcf_{s_2}, \cdots, Tcf_{s_{n_1+n_2}}\rangle$.

(1) If $(cl_i, cl_{n_1}, cl_{n_1+1})=(\rho_2, \rho_2-1, \rho_2)$ for some $i\in \{1, 2, \cdots, n_1-1\}$, $F(\Phi_a)-F(\Phi_b)=0$.

(2) If $cl_{n_1}\notin\{1,2\}$ and $(cl_i, cl_{n_1}, cl_{n_1+1})\neq(\rho_2, \rho_2-1, \rho_2)$ for all $i\in \{1, 2, \cdots, n_1-1\}$, $F(\Phi_a)-F(\Phi_b)>0$.

(3) If $cl_{n_1}\in\{1,2\}$, $F(\Phi_a)-F(\Phi_b)>0$.}\end{lemma}

\noindent{Proof}: (1) Since $(cl_i, cl_{n_1})=(\rho_2, \rho_2-1)$, from Lemma \ref{lemma1.2}, under Assumptions \ref{ass1.3.2} and \ref{ass1.3.4},
it can be checked that $F(\Phi_a)-F(\Phi_b)=0$.

(2) When $(cl_{n_1}, cl_{n_1+1})=(\rho_2-1, \rho_2)$, $cl_i\neq \rho_2$ for all $i\in \{1, 2, \cdots, n_1-1\}$. It follows that $f_{\rho_2}(\Phi_a)+f_{\rho_2-1}(\Phi_a)-f_{\rho_2}(\Phi_b)-f_{\rho_2-1}(\Phi_b)=\delta$. Note that $f_{k}(\Phi_a)=f_{k}(\Phi_b)$ for all $k\neq \rho_2-1, \rho_2$. Therefore, $F(\Phi_a)-F(\Phi_b)>0$. When $(cl_{n_1}, cl_{n_1+1})\neq (\rho_2-1, \rho_2)$, from Assumptions \ref{ass1.3.3} and \ref{ass1.3.4}, $T_{cl_{n_1}cl_{(n_1+1)}}\geq 1.5T_0>T_0+2\delta$.
From Assumption \ref{ass1.3.1}, there are at most two classes of aircraft satisfying $T_{ii}=T_0+\delta$. Thus, during the transition from $\Phi_a$ to $\Phi_b$, there are at most two aircraft such that the separation times between them and their trailing aircraft become $T_0+\delta$ from $T_0$ or $T_{cl_{n_1}cl_{(n_1+1)}}$. It follows that $F(\Phi_a)-F(\Phi_b)\geq T_0+2\delta-T_0-2\delta>0$.

(3) Similar to the proof of statement (2), during the transition from $\Phi_a$ to $\Phi_b$, there are at most two aircraft of classes $\rho_1$ and $\rho_2$ such that the separation times between them and their trailing aircraft become $T_0+\delta$ from $T_0$.
When $\phi_2$ has aircraft of class no larger than $2$, we have  \begin{eqnarray*}\begin{array}{lll}&F(\Phi_a)-F(\Phi_b)\\
&\geq f_1(\Phi_a)+f_2(\Phi_a)-f_1(\Phi_b)-f_2(\Phi_b)-2\delta\\
&\geq T_{cl_{n_1}cl_{n_1+1}}-1.5T_0-2\delta.\end{array}\end{eqnarray*} When $\phi_2$ has no aircraft of class no larger than $2$, we have
\begin{eqnarray*}\begin{array}{lll}&F(\Phi_a)-F(\Phi_b)\\&\geq f_1(\Phi_a)\!+f_2(\Phi_a)\!-f_1(\Phi_b)\!-f_2(\Phi_b)-T_0-2\delta\\
&\geq T_{cl_{n_1}cl_{n_1+1}}-T_0-2\delta.\end{array}\end{eqnarray*}
Note from Assumption \ref{ass1.3.4} that $T_{cl_{n_1}cl_{n_1+1}}>1.5T_0+2\delta$ since $cl_{n_1}\in\{1,2\}$. Then it follows that $F(\Phi_a)-F(\Phi_b)>0$.

\begin{remark}{\rm Lemma \ref{lemma1.31} actually discusses the case of a breakpoint generated from a class-monotonically-decreasing sequence.
The existence of the breakpoint might increase the value of the objective function $F(\cdot)$. It is important to avoid the occurrence of breakpoints so as to decrease the value of the objective function $F(\cdot)$.}\end{remark}

\begin{remark}{\rm In Lemmas \ref{lemmap1} and \ref{lemma1.31}, it is shown that the value of the objective function $F(\cdot)$ can be decreased by interposing an aircraft to replace the original breakpoint or merging the breakpoint into a class-monotonically-decreasing sequence.}\end{remark}

Let $\mathrm{Exs}_i(\phi)$ be a function such that $\mathrm{Exs}_i(\phi)=1$ if the aircraft sequence $\phi$ contains aircraft of class $i$, and $\mathrm{Exs}_i(\phi)=0$ if the aircraft sequence $\phi$ contains no aircraft of class $i$, and let $\mathrm{sgn}(x)$ be a sign function such that $\mathrm{sgn}(x)=1$ if $x>0$ and $\mathrm{sgn}(x)=0$ if $x\leq 0$.

\begin{theorem}\label{lemma1.3s}{\rm  Consider a landing sequence $\Phi_a=\langle\phi_1, \phi_2, \cdots, \phi_s\rangle$ for a positive integer $s$, where $\phi_i=\langle Tcf_{i1}, Tcf_{i2}, \cdots,$ $ Tcf_{ic_i}\rangle$ is a class-monotonically-decreasing sequence for some positive integer $c_i$ and all $i\in \{1, 2, \cdots, s\}$, and $cl_{jc_j}<cl_{(j+1)1}$ for all $j\in \{1, 2, \cdots, s-1\}$. Suppose that Assumptions \ref{ass1.3.1}-\ref{ass1.4} hold for each landing aircraft sequence.
 Merge the aircraft sequence $\Phi_a$ to form a new class-monotonically-decreasing sequence $\Phi_b$. The following statements hold.


(1) $F(\Phi_a)=F(\Phi_b)+T_{d1}-T_{cl_{sc_s}1}-\mathrm{sgn}(\sum_{i=1}^s \mathrm{Exs}_{\rho_1}(\phi_i))(\sum_{i=1}^s \mathrm{Exs}_{\rho_1}(\phi_i)-1)\delta-$ $\mathrm{sgn}(\sum_{i=1}^s \mathrm{Exs}_{\rho_2}(\phi_i))(\sum_{i=1}^s \mathrm{Exs}_{\rho_2}(\phi_i) -1)\delta+$ $\sum_{i=1}^{s-1}[T_{cl_{ic_i}cl_{(i+1)1}}-T_{cl_{ic_i}1}]$, where $d$ denotes the class of the last aircraft of $\Phi_b$.

(2) Suppose that $(cl_j, cl_{kc_k}, cl_{(k+1)1})=(\rho_2, \rho_2-1, \rho_2)$ for some $j\in \{k1, k2, \cdots, k(c_k-1)\}$ and all $k\in\{1,2, \cdots, s-1\}$. It follows that $F(\Phi_a)-F(\Phi_b)=0$.

(3) Suppose that there is a positive integer $k_0\in \{1, 2, \cdots, s-1\}$ such that $(cl_j, cl_{k_0c_{k_0}}, cl_{(k_0+1)1})\neq(\rho_2, \rho_2-1, \rho_2)$ for all $j\in \{k1, k2, \cdots, k(c_k-1)\}$. It follows that $F(\Phi_a)-F(\Phi_b)>0$.}\end{theorem}

\noindent{Proof:} (1) Calculate $F(\Phi_a)-F(\Phi_b)$ in three parts: the separation times between the breakpoint aircraft and their trailing aircraft, the separation times between the aircraft of classes $\rho_1$ and $\rho_2$ and their trailing aircraft, and the separation times between the last aircraft in $\Phi_a$ and $\Phi_b$ and their trailing aircraft. It can be obtained that $F(\Phi_a)=F(\Phi_b)-T_{d1}+T_{cl_{sc_s}1}-\sum_{i=1}^s[\mathrm{Exs}_{\rho_1}(\phi_i)+\mathrm{Exs}_{\rho_2}(\phi_i)]\delta+\sum_{i=1}^{s-1}[T_{cl_{ic_i}cl_{(i+1)1}}-T_{cl_{ic_i}1}]$.

(2) Since $(cl_j, cl_{kc_k}, cl_{(k+1)1})=(\rho_2, \rho_2-1, \rho_2)$ for some $j\in \{k1, k2, \cdots, k(c_k-1)\}$ and all $k\in\{1,2, \cdots, s-1\}$, the subsequence $\phi_1$ contains at most the aircraft of classes $\rho_2-1$, $\rho_2$, $\cdots, \eta$, each of the subsequences $\phi_2, \phi_3, \cdots, \phi_{s-1}$ contains the aircraft of classes $\rho_2-1$ and $\rho_2$, and the subsequence $\phi_s$ contains at most the aircraft of classes $\rho_2$, $\rho_2-1$, $\cdots, 1$. Note that $T_{d1}-T_{cl_{sc_s}1}\geq 0$ since $d\leq cl_i$ for all $Tcf_i$ in $\Phi_b$. Thus, it can easily be checked that $F(\Phi_a)-F(\Phi_b)=0$ based on the statement (1).

(3) When $(cl_{ic_{i}}, cl_{(i+1)1})=(\rho_2-1, \rho_2)$ for some $i\in \{1, \cdots, s-1\}$, under the condition given in this statement, it follows that $T_{cl_{ic_{i}}cl_{(i+1)1}}-T_{cl_{ic_{i}}1}=\delta$ and $\mathrm{Exs}_{\rho_2}(\phi_{i})=\mathrm{Exs}_{\rho_1}(\phi_{i})=0$.
When $(cl_{ic_{i}}, cl_{(i+1)1})\neq (\rho_2-1, \rho_2)$ for some $i\in \{1,2, \cdots, s-1\}$, from Assumptions \ref{ass1.3.2}, \ref{ass1.3.3} and \ref{ass1.3.4}, if $cl_{ic_{i}}=1$, then $T_{cl_{ic_{i}}cl_{(i+1)1}}-T_{cl_{ic_i}1}>0.5T_0$, if $cl_{ic_{i}}=2$, then $T_{cl_{ic_{i}}cl_{(i+1)1}}-T_{cl_{ic_i}1}>2\delta$, and if $cl_{ic_{i}}\geq 3$, then $T_{cl_{ic_{i}}cl_{(i+1)1}}-T_{cl_{ic_i}1}\geq 0.5T_0$. Note that $T_{d1}-T_{cl_{sc_s}1}\geq 0$ since $d\leq cl_i$ for all $Tcf_i$ in $\Phi_b$.
Based on the statement (1), it can be obtained that $F(\Phi_a)-F(\Phi_b)>0$.

\begin{remark}{\rm In Theorem \ref{lemma1.3s}, the terms $\mathrm{sgn}(\sum_{i=1}^s \mathrm{Exs}_{\rho_1}(\phi_i))(\sum_{i=1}^s \mathrm{Exs}_{\rho_1}(\phi_i)-1)\delta+\mathrm{sgn}(\sum_{i=1}^s \mathrm{Exs}_{\rho_2}(\phi_i))(\sum_{i=1}^s \mathrm{Exs}_{\rho_2}(\phi_i)$ $-1)\delta$ mean that the dispersion of aircraft of classes $\rho_1$ and $\rho_2$ in different subsequences might decrease the value of the objective function $F(\cdot)$.}\end{remark}

\begin{remark}{\rm Theorem \ref{lemma1.3s} gives a calculation method for the objective function $F(\cdot)$ which can also be applied to the case when the landing or takeoff times of the aircraft are subject to different constraints, possibly resulting in the occurrence of the resident-point aircraft. Moreover, it should be noted that each subsequence $\phi_i$ might contain only one aircraft or the aircraft of the same class.}\end{remark}

\begin{theorem}\label{theorem1.2}{\rm Consider a group of landing aircraft $\{Tcf_1, Tcf_2, \cdots, Tcf_n\}$, where $cl_1 \geq cl_2 \geq cl_3 \geq \cdots \geq cl_n$. Suppose that Assumptions \ref{ass1.3.1}-\ref{ass1.4} hold for each landing aircraft sequence.
Then the optimization problem (\ref{optim1}) can be solved if the landing aircraft sequence is taken as $\phi_0=\langle Tcf_1, Tcf_2, \cdots, Tcf_n \rangle$.}\end{theorem}

\noindent{Proof}: This theorem can be regarded as a special case of Theorem \ref{lemma1.3s}.

 Theorem \ref{theorem1.2} shows that when the set of allowable landing times for each aircraft is $Tf_k=[t_0, +\infty]$, if all aircraft form a class-monotonically-decreasing sequence without breakpoints, the optimization problem (\ref{optim1}) can be solved.

From Theorems \ref{theorem1.1} and \ref{lemma1.3s}, it can be seen that the optimality of the aircraft sequence is heavily related to the resident-point aircraft, breakpoint aircraft and the aircraft with separation times larger than $T_0$. To solve optimization problem (\ref{optim1}), we can focus on these special aircraft so as to study the local minimum points of the objective function $F(\cdot)$. Due to the high complexity of the aircraft sequence and the constraints on the landing times of the aircraft, in the following, we only introduce some class-sequence sets, which can be calculated out offline, to study some typical scenarios and more general scenarios can be studied in a similar way.

Define class-sequence sets as $\Psi_0=\{\langle i_1, i_2, i_3, i_4, i_5\rangle\mid {i_1}\geq {i_2}\geq {i_3}, {i_4}\geq {i_5}, i_1, i_2, i_3, i_4, i_5\in \mathcal{I}\}$, $\Psi_1=\{\langle i_1, i_2, i_3, i_4, i_5\rangle\mid {i_1}\geq {i_2}, {i_2}<{i_3}, {i_4}\geq {i_5}, i_1, i_2, i_3, i_4, i_5\in \mathcal{I}\}$, $\Psi_2=\{\langle i_1, i_2, i_3, i_4, i_5\rangle\mid {i_1}<{i_2}, {i_2}\geq {i_3}, {i_4}\geq {i_5}, i_1, i_2, i_3, i_4, i_5\in \mathcal{I}\}$, $\Psi_3=\{\langle i_1, i_2, i_3, i_4, i_5\rangle\mid {i_1}\geq {i_2}, {i_2}<{i_3}, {i_4}< {i_5}, i_1, i_2, i_3, i_4, i_5\in \mathcal{I}\}$, $\Psi_4=\{\langle i_1, i_2, i_3, i_4, i_5\rangle\mid {i_1}<{i_2}, {i_2}\geq {i_3}, {i_4}<{i_5}, i_1, i_2, i_3, i_4, i_5\in \mathcal{I}\}$ and $\Psi_5=\{\langle i_1, i_2, i_3, i_4, i_5\rangle\mid {i_1}\geq {i_2}\geq {i_3}, {i_4}< {i_5}, i_1, i_2, i_3, i_4, i_5\in \mathcal{I}\}$. Let $\Gamma(\langle i_1, i_2, i_3, i_4, i_5\rangle)$ be a class-sequence function such that  $\Gamma(\langle i_1, i_2, i_3, i_4, i_5\rangle)=T_{{i_1}{i_2}}+T_{{i_2}{i_3}}+T_{{i_4}{i_5}}-T_{{i_1}{i_3}}-T_{{i_4}{i_2}}-T_{{i_2}{i_5}}$.

The role of the conditions in the definitions of $\Psi_0, \Psi_1, \Psi_2, \Psi_3, \Psi_4,\Psi_5$ is to describe relationship between the classes of the adjacent aircraft. For example, the condition that $i_1<i_2, i_2\geq i_3$ in $\Psi_2$ means that a breakpoint aircraft is formed at an aircraft of class $i_2$.

\begin{definition}{\rm (1) Let $\Omega_0\subseteq \Psi_0$ be a class-sequence set such that $\Gamma(\langle i_1, i_2, i_3, i_4, i_5\rangle)\geq 0$ for any $\langle i_1, i_2, i_3, i_4, i_5\rangle\in \Omega_0$, and $\Gamma(\langle i_1, i_2, i_3, i_4, i_5\rangle)<0$ for any $\langle i_1, i_2, i_3, i_4, i_5\rangle\in \Psi_0-\Omega_0$.

(2) Let $\Omega_1\subseteq \Psi_1$ be a class-sequence set such that $\Gamma(\langle i_1, i_2, i_3, i_4, i_5\rangle)\geq 0$ for any $\langle i_1, i_2, i_3, i_4, i_5\rangle\in \Omega_1$, and $\Gamma(\langle i_1, i_2, i_3, i_4, i_5\rangle)<0$ for any $\langle i_1, i_2, i_3, i_4, i_5\rangle\in \Psi_1-\Omega_1$.

(3) Let $\Omega_2\subseteq \Psi_2$ be a class-sequence set such that $\Gamma(\langle i_1, i_2, i_3, i_4, i_5\rangle)\geq 0$ for any $\langle i_1, i_2, i_3, i_4, i_5\rangle\in \Omega_2$, and $\Gamma(\langle i_1, i_2, i_3, i_4, i_5\rangle)<0$ for any $\langle i_1, i_2, i_3, i_4, i_5\rangle\in \Psi_2-\Omega_2$.

(4) Let $\Omega_3\subseteq \Psi_3$ be a class-sequence set such that $\Gamma(\langle i_1, i_2, i_3, i_4, i_5\rangle)\geq 0$ for any $\langle i_1, i_2, i_3, i_4, i_5\rangle\in \Omega_3$, and $\Gamma(\langle i_1, i_2, i_3, i_4, i_5\rangle)<0$ for any $\langle i_1, i_2, i_3, i_4, i_5\rangle\in \Psi_3-\Omega_3$.

(5) Let $\Omega_4\subseteq \Psi_4$ be a class-sequence set such that $\Gamma(\langle i_1, i_2, i_3, i_4, i_5\rangle)\geq 0$ for any $\langle i_1, i_2, i_3, i_4, i_5\rangle\in \Omega_4$, and $\Gamma(\langle i_1, i_2, i_3, i_4, i_5\rangle)<0$ for any $\langle i_1, i_2, i_3, i_4, i_5\rangle\in \Psi_4-\Omega_4$.

(6) Let $\Omega_5\subseteq \Psi_5$ be a class-sequence set such that $\Gamma(\langle i_1, i_2, i_3, i_4, i_5\rangle)\geq 0$ for any $\langle i_1, i_2, i_3, i_4, i_5\rangle\in \Omega_5$, and $\Gamma(\langle i_1, i_2, i_3, i_4, i_5\rangle)<0$ for any $\langle i_1, i_2, i_3, i_4, i_5\rangle\in \Psi_5-\Omega_5$.}\end{definition}

Though the forms of the sets $\Omega_0, \Omega_1,\Omega_2, \Omega_3, \Omega_4, \Omega_5$ seem to be complex, they can be obtained offline by a small amount of calculations since the aircraft are usually categorized into up to $6$ or $7$ classes in practical applications. Moreover, it is clear that $\Omega_i \cap \Omega_j=\emptyset$ for all $i\neq j$, but $\Omega_i\cap (\Psi_j-\Omega_j)$ might not be empty. For less calculation of the sets $\Omega_0, \Omega_1,\Omega_2, \Omega_3, \Omega_4, \Omega_5$, the relationship between $\Omega_i$ and $\Psi_j-\Omega_j$ can be further analyzed by splitting them into proper subsets.

To make the definitions of the sets $\Omega_0, \Omega_1,\Omega_2, \Omega_3, \Omega_4, \Omega_5$ more easily understood from a view point of class-monotonically-decreasing sequences, we give the following definition as an example.

\begin{definition}\label{ass1.5.1.1}{\rm Let $\Theta_{0}=\{\langle i, k, j, h\rangle\mid i< j<k, i< h<k, i,j,h,k\in \mathcal{I}\}$ be a class sequence set  such that $T_{i k}+T_{0}\geq T_{i j}+T_{jk}$ for any $\langle i, k, j, h\rangle\in \Theta_{0}$ and $T_{i k}+T_{0}<T_{i j}+T_{jk}$ for any $\langle i, k, j, h\rangle\notin \Theta_{0}$.   Let $\Theta_{1}=\{\langle i, k, j, h\rangle\mid k<i<j, k<h<j, i,j,h,k\in \mathcal{I}\}$ be a class sequence set  such that $T_{i k}+T_{0}\geq T_{i j}+T_{jk}$ for any $\langle i, k, j, h\rangle\in \Theta_{1}$ and $T_{i k}+T_{0}<T_{i j}+T_{jk}$ for any $\langle i, k, j, h\rangle\notin \Theta_{1}$.}\end{definition}

In Definition \ref{ass1.5.1.1}, the sets $\Theta_{0}$ and $\Theta_{1}$ are discussed for some scenarios when one breakpoint is split into two breakpoints or two breakpoints are merged into one breakpoint, and more scenarios, e.g., the split of one breakpoint into more than two breakpoints, can be analyzed in a similar way. All of these can be used to reduce the computation cost of our algorithms proposed later.

\begin{remark}{\rm Under Assumption \ref{ass1.4}, consider landing sequences $\phi_1=\langle Tcf_1, Tcf_2, Tcf_3, Tcf_4, Tcf_5\rangle$ and $\phi_2=\langle Tcf_1, Tcf_3, Tcf_4, Tcf_2, Tcf_5\rangle$, where $\langle cl_1, cl_2, cl_3, cl_4\rangle\in\Theta_{0}$ and $cl_5\leq cl_4<cl_2$. Since $\langle cl_1, cl_2, cl_3, cl_4\rangle\in\Theta_{0}$, $cl_1<cl_3<cl_2$, implying that $cl_2\geq 3$, and hence $T_{cl_1cl_3}+T_{cl_4cl_2}\leq T_{cl_1 cl_2}+T_0=T_{cl_1 cl_2}+T_{cl_2cl_3}$. Note that $T_{cl_4cl_5}\geq T_{cl_2cl_5}=T_0$. Therefore, $F(\phi_1)\geq F(\phi_2)$.
If the inequality $cl_5\leq cl_4<cl_2$ does not hold, it is hard to determine which of $F(\phi_1)$ and $F(\phi_2)$ is larger without more information. From this example, it can be observed that
the class sequence in $\Theta_0$ might be optimal in some situations but might not be in some other situations which is additionally related to the orders of other aircraft. It should be noted that when the conditions of the class sequence sets $\Theta_{0}$ and $\Theta_{1}$
do not hold, then the two breakpoints can be merged into one breakpoint to reduce the value of the objective function $F(\cdot)$.
}\end{remark}

\begin{proposition}\label{prop1}{\rm
Consider a sequence $\langle Tcf_1, Tcf_2, Tcf_3, Tcf_4, Tcf_5\rangle$, where $(cl_1, cl_2, cl_3, cl_4, cl_5)=(i_1, i_2, i_3, i_4, i_5)$. Suppose that Assumptions \ref{ass1.3.1}-\ref{ass1.4} hold for each landing sequence.

(1) Suppose that ${i_1}\geq {i_2}\geq {i_3}, {i_4}\geq {i_5},$ and $i_1, i_2, i_3, i_4, i_5\in \mathcal{I}$. If $(i_2, j, i_4)\neq (\rho_2, \rho_2, \rho_2-1)$ for all $j\in\{i_1,i_3\}$ with $i_2>i_4$ or $i_2<i_5$, $\langle i_1, i_2, i_3, i_4, i_5\rangle\notin \Omega_0$.

(2) Suppose that ${i_1}\geq {i_3}, {i_4}\geq {i_2}\geq {i_5},$ and $i_1, i_2, i_3, i_4, i_5\in \mathcal{I}$. If $(i_2, j, i_1)\neq (\rho_2, \rho_2, \rho_2-1)$ for all $j\in\{i_4,i_5\}$ with $i_2>i_1$ or $i_2<i_3$, $\langle i_1, i_2, i_3, i_4, i_5\rangle\in \Omega_0$.

(3) Suppose that ${i_1}\geq {i_2}, {i_2}<{i_3}, {i_4}\geq {i_5},$ and $ i_1, i_2, i_3, i_4, i_5\in \mathcal{I}$. If $(i_2, j)\neq (\rho_2, \rho_2)$ for all $j\in\{i_4,i_5\}$ with $i_1> i_3$ and $i_4\geq i_2\geq i_5$, $\langle i_1, i_2, i_3, i_4, i_5\rangle\in \Omega_1$.

(4) Suppose that ${i_1}<{i_2}, {i_2}\geq {i_3}, {i_4}\geq {i_5},$ and $i_1, i_2, i_3, i_4, i_5\in \mathcal{I}$. If $(i_2, j)\neq (\rho_2, \rho_2)$ for all $j\in\{i_4,i_5\}$ with $i_1\geq i_3$ and $i_4\geq i_2\geq i_5$, $\langle i_1, i_2, i_3, i_4, i_5\rangle\in \Omega_2$.

(5) Suppose that ${i_1}\geq {i_2}\geq {i_3}, {i_4}< {i_5}, i_1, i_2, i_3, i_4, i_5\in \mathcal{I}$. If $\langle i_4,i_5\rangle\in E$ and $i_2=2$, $\langle i_1, i_2, i_3, i_4, i_5\rangle\in \Omega_5$.

}\end{proposition}

\noindent{Proof:} Using the previous lemmas and theorems, this proposition can be proved and its proof is omitted.

It should be noted that the transitions in Proposition \ref{prop1} (1)(2) can be regarded as reciprocal inverse transformations of each other.

\begin{proposition}\label{prop2}{\rm
Consider two sequences $\phi_1=\langle Tcf_1, Tcf_2, Tcf_3, Tcf_4, Tcf_5\rangle$ and $\phi_2=\langle Tcf_1, Tcf_3, Tcf_4, Tcf_2,$ $ Tcf_5\rangle$, where $(cl_1, cl_2, cl_3, cl_4, cl_5)=(i_1, i_2, i_3, i_4, i_5)$. Suppose that Assumptions \ref{ass1.3.1}-\ref{ass1.4} hold for all landing aircraft sequences, ${i_1}\geq {i_2}\geq {i_3}, {i_4}\geq {i_5},$ and $i_1, i_2, i_3, i_4, i_5\in \mathcal{I}$. If $i_2=j=k$ for some $j\in\{i_1, i_3\}$ and some $k\in \{\rho_4, \rho_5\}$, then $F(\phi_1)=F(\phi_2)$.
}\end{proposition}

In Proposition \ref{prop1}, we only discuss some typical scenarios for the elements of the class-sequence sets $\Omega_0-\Omega_5$, and more discussions can be made according to Assumptions \ref{ass1.3.1}-\ref{ass1.4} and the practical situations.

It should be noted that the class-sequence sets $\Omega_1$, $\Omega_2$, $\Psi_3-\Omega_3$ and $\Psi_4-\Omega_4$ might contain elements for the scenarios of the split of one breakpoint into two breakpoints whereas
 the class-sequence sets $\Omega_3$, $\Omega_4$, $\Psi_1-\Omega_1$ and $\Psi_2-\Omega_2$ might contain elements for the scenarios for the merger of two breakpoints into one breakpoint, which might be able to decrease the value of the objective function $F(\cdot)$.

\begin{remark}{\rm For given aircraft sequences, the class sequences might not belong to $\Omega_i$ or $\Psi_i-\Omega_i$, $i=0,1,2, \cdots, 5$. But in special situations, aircraft sequences can be adjusted to satisfy the conditions of the class-sequence sets of $\Omega_i$ or $\Psi_i-\Omega_i$, $i=0,1,2, \cdots, 5$, by changing the aircraft orders.
For example, consider the sequences $\phi_1=\langle Tcf_1, Tcf_2, Tcf_3, Tcf_4, Tcf_5\rangle$ and $\phi_2=\langle Tcf_4, Tcf_2, Tcf_3, Tcf_1, Tcf_5\rangle$, where $\langle cl_4, cl_2, cl_3, cl_1, cl_5\rangle\in \Omega_0$.
 It is clear that the class sequence of the aircraft in
$\phi_1$ does not satisfy the set $\Omega_0$ but can be converted into $\phi_2$ to make the class sequence of the aircraft satisfies the class-sequence set $\Omega_0$ by exchanging the orders of the aircraft $Tcf_1$ and $Tcf_4$.}\end{remark}

\begin{proposition}\label{theorem4}{\rm Consider two sequences $\Phi_1=\langle \phi_{1}, \phi_{2}, \phi_{3}, Tcf_{k_0}, \phi_{4}\rangle$  and $\Phi_2=\langle \phi_{1}, \phi_{3}, \phi_{2},  Tcf_{k_0}, \phi_{4}\rangle$, where the aircraft of each $\phi_{i}$ have the same class $k_i$ for $i=1, 2, 3, 4$. Suppose that Assumptions \ref{ass1.3.1}-\ref{ass1.4} hold for the sequences $\Phi_1$ and $\Phi_2$. If $k_2<k_1<k_3$, $cl_{k_0}=k_3$ and $(k_1, k_3)\neq (\rho_2-1,\rho_2)$, then $F(\Phi_1)<F(\Phi_2)$.
}\end{proposition}

\noindent{Proof:} Calculating $F(\Phi_2)-F(\Phi_1)$, we have   $F(\Phi_2)-F(\Phi_1)=T_{k_1k_3}+T_{k_3k_2}+(m_3-1)T_{k_3k_3}-T_{k_1k_2}-m_3T_{k_3k_3}=T_{k_1k_3}+T_{k_3k_2}-T_{k_3k_3}-T_{k_1k_2}$, where $m_3$ denotes the number of aircraft in $\phi_3$. Note that $k_2<k_1<k_3$. It follows that
$k_3\geq 3$ and $k_3-k_2\geq 2$. From Assumption \ref{ass1.3.4}, when $k_1\leq 2$, $T_{k_1k_3}> 1.5T_0+2\delta=T_{k_1k_2}+2\delta$, $T_{k_3k_2}=T_0$, and $T_{k_3k_3}\leq T_0+\delta$. Hence, $F(\Phi_1)<F(\Phi_2)$ when $k_1\leq 2$.
  When $k_1\geq 3$, it follows that  $T_{k_1k_3}\geq 1.5T_0$ since $(k_1, k_3)\neq (\rho_2-1,\rho_2)$, $T_{k_3k_2}=T_{k_1k_2}=T_0$, and $T_{k_3k_3}\leq T_0+\delta$. Therefore, $F(\Phi_1)<F(\Phi_2)$ when $k_1\geq 3$. This proposition is proved.

\begin{remark}{\rm Proposition \ref{theorem4} shows that if the consecutive aircraft are of the same class and one of them needs to be moved, all the consecutive aircraft with the same class are usually needed to be moved together simultaneously.}\end{remark}

\section{Takeoff scheduling problem on a single runway}\label{sectak}

In this section, we study the takeoff scheduling problem on a single runway. The main analysis idea is similar to the landing scheduling problem in Sec. \ref{seclanding}.

Let $D_{ij}$ represent the minimum separation time between an aircraft of class $i$ and a trailing aircraft of class $j$ on a single runway without considering the influence of other aircraft.

\begin{assumption}\label{ass2.3.1}{\rm (1) When $i=1,2, 3, \eta$, $D_{ii}=(1+1/3)T_0$.

(2) When $i\neq 1, 2, 3, \eta$, $D_{ii}=T_0$.}\end{assumption}

\begin{assumption}\label{ass2.3.2}{\rm (1) $D_{21}=(1+1/3)T_0$.

(2) When $i>j$ and $i\geq 3$, $D_{ij}=T_0$.}\end{assumption}

\begin{assumption}\label{ass2.3.3}{\rm
(1) For all $k\leq i\leq j$, $D_{ij}\leq D_{kj}\leq 3T_0$ and $D_{ki}\leq D_{kj}\leq 3T_0$.

(2) For all $k\leq j\leq i$, $D_{ik}<D_{ij}+D_{jk}$, $D_{ki}< D_{ji}+D_{kj}$.}\end{assumption}

\begin{lemma}\label{lemma22s}{\rm For all $i,j, k\in \mathcal{I}$, $D_{ik}<D_{ij}+D_{jk}$.}\end{lemma}

\begin{assumption}\label{ass2.3.41}{\rm
(1) For $k=1, 2$, $D_{k(k+1)}=D_{kk}+T_0/6$.

(2) For $k=3, \rho_2-1$, $D_{k(k+1)}=D_{kk}$.

(3) For $k=1$, $D_{k3}=D_{23}+T_0/3$. For $k=1$ and all $k+2\leq j \leq \eta$, $D_{kj}-D_{(k+1)j}=2T_0/3$.

(4) For $k=3$, $D_{k\eta}=D_{(k+1)\eta}$, and for $k=\rho_2-1$, $D_{k\eta}=D_{\rho_2\rho_2}+T_0$.

(5) For all $2\leq k<j$ and $j\geq 3$ such that $(k,j)\neq (3,\eta)$, $(k,j)\neq (\rho_2-1,\rho_2)$ and $(k,j)\neq (\rho_2-2,\rho_2)$, $D_{kj}=D_{(k+1)j}+T_0/3$.

}\end{assumption}

\begin{lemma}\label{lemma2.111}{\rm Under Assumptions \ref{ass2.3.1}-\ref{ass2.3.41}, the following statements hold.

(1) Let $E_1=\{\langle 3,j\rangle, 3\leq j\leq \eta, \langle 4,\eta\rangle\}$ be a sequence set. Then
$D_{2j}-D_{kj}=T_0/3$ for all $\langle k,j\rangle\in E_1$ and for all $\langle k,j\rangle\notin E_1$ with $2<k\leq j$, $D_{2j}-D_{kj}>0.5T_0$.

(2) Let $E_{20}=\{\langle 4,\eta\rangle\}$ and $E_{21}=\{\langle4,j\rangle, 4<j<\eta, j\neq \rho_2, \langle 5,\eta\rangle\}$ be two sequence sets. Then, $D_{3j}-D_{kj}=0$ for $\langle k,j\rangle\in E_{20}$, $D_{3j}-D_{kj}=T_0/3$ for $\langle k,j\rangle\in E_{21}$, and $D_{3j}-D_{kj}>T_0/3$ for all $\langle k,j\rangle\notin E_{20}\cup E_{21}$ with $3<k<j$.

(3) Let $E_{30}=\{\langle 1,\rho_2\rangle, \langle 2,\rho_2\rangle, \langle 3,\rho_2\rangle\}$ and $E_{31}=\{\langle \rho_2-1,\eta\rangle\}$ be a sequence set. Then $D_{kj}-D_{k(\rho_2-1)}=T_0/3$ for $\langle k,j\rangle\in E_{30}$ and $D_{kj}-D_{\rho_2j}=T_0/3$ for $\langle k,j\rangle\in E_{31}$.

(4) Let $E_4=\{\langle \eta,\eta\rangle\}$ be a sequence set. Then $D_{(\eta-1)j}-D_{kj}=T_0/3$ for $\langle k,j\rangle\in E_4$.

}\end{lemma}

\begin{remark}{\rm In Lemma \ref{lemma2.111}, we only discuss some typical scenarios which might result in local minimum point and more general scenarios can be discussed according to the class-sequence sets defined later.}\end{remark}

As discussed in Theorem \ref{theorem1.1}, in Theorem \ref{theorempp2.1}, we study the aircraft relevance and show that when the takeoff times of the aircraft have no constraints, the occurrence of resident-point aircraft should be avoided to ensure the optimality of the sequence.

\begin{theorem}\label{theorempp2.1}{\rm Consider a takeoff aircraft sequence $\phi=\langle Tcf_1, Tcf_2, \cdots, Tcf_{n}\rangle$. The following statements hold.

(1) For all $i=3, 4, \cdots, n$, the aircraft $Tcf_{i}$ is not relevant to $Tcf_1$, $Tcf_2, \cdots,$ $Tcf_{i-2}$.

(2) If the aircraft $Tcf_j$ is relevant to the aircraft $Tcf_i$, $j=i+1$.

(3) Suppose that the orders of all the takeoff aircraft in $\phi$ is fixed.
The optimization problem (\ref{optim1}) is solved if and only if $Tcf_{i+1}$ is relevant to $Tcf_i$, $i=1,2, \cdots, n-1$.}\end{theorem}

\begin{assumption}\label{ass2.4}{\rm Consider a takeoff aircraft sequence $\phi=\langle Tcf_1, Tcf_2, \cdots, Tcf_{n}\rangle$. Suppose that $Tf_k=[t_0, +\infty]$ for all $k$, and for each pair of adjacent aircraft in the takeoff sequence $\phi$, each aircraft is relevant to its leading aircraft.}\end{assumption}

\begin{lemma}\label{lemma2p}{\rm Consider two takeoff aircraft sequences $\phi_0=\langle Tcf_1, Tcf_2, Tcf_3, Tcf_4\rangle$ and $\phi_1=\langle Tcf_1, Tcf_3, Tcf_2,$ $ Tcf_4\rangle$. Under Assumptions \ref{ass2.3.1}-\ref{ass2.4}, the following statements hold.

(1) Suppose that $cl_1=cl_2=2$. If $\langle cl_3, cl_4\rangle\in E_1$, $F(\phi_0)-F(\phi_1)=0$. If $\langle cl_3, cl_4 \rangle\notin E_1$, {$F(\phi_0)-F(\phi_1)<0$.}

(2) Suppose that $cl_1=cl_2=3$. If $\langle cl_3, cl_4 \rangle\in E_{20}$, $F(\phi_0)-F(\phi_1)=T_0/3$. If $\langle cl_3, cl_4\rangle\in E_{21}$, $F(\phi_0)-F(\phi_1)=0$. If $\langle cl_3, cl_4\rangle\notin E_{20}\cup E_{21}$ with $3<cl_3<cl_4$, $F(\phi_0)-F(\phi_1)<0$.

(3) Suppose that $cl_1=\eta, cl_2=\rho_2-1$. If $\langle cl_3, cl_4\rangle\in E_{30}$, $F(\phi_0)-F(\phi_1)=T_0/3$. Suppose that $cl_1=\eta, cl_2=\rho_2$. If $\langle cl_3, cl_4\rangle\in E_{31}$,
$F(\phi_0)-F(\phi_1)=T_0/3$.

(4) Suppose that $cl_1=cl_2=\eta-1$. If $\langle cl_3, cl_4\rangle\in E_{4}$, $F(\phi_0)-F(\phi_1)=0$.}\end{lemma}

\noindent{Proof:} By simple calculations, this lemma can be proved and hence its proof is omitted.

Lemma \ref{lemma2p} gives examples to illustrate the scenarios in Lemma \ref{lemma2.111}.

In the following, we first give a calculation method for the objective function $F(\cdot)$ when the order of some aircraft is changed as discussed for landing sequences.

\begin{lemma}\label{lemmap1314}{\rm Consider a takeoff aircraft sequence $\Phi_a=\langle\phi_1, Tcf_{s_2}, Tcf_{s_3},\phi_2\rangle$, where $\phi_1=\langle Tcf_1, Tcf_2, \cdots, Tcf_{s_1}\rangle$, $s_2=s_1+1$, $s_3=s_1+2$, $\phi_2=\langle Tcf_{s_3}, Tcf_{s_3+1}, \cdots, Tcf_{n}\rangle$,  and $s_1, s_2, s_3$ are three positive integers. Suppose that Assumptions \ref{ass2.3.1}-\ref{ass2.4} hold for each landing aircraft sequence.

(1) Move $Tcf_{s_2}$ to be between $Tcf_{h_1}$ and $Tcf_{h_1+1}$ and convert $\phi_2$ into a new sequence $\phi_3$, where $Tcf_{h_1}, Tcf_{h_1+1}\in \phi_2$.
Let $\Phi_b=\langle \phi_1, Tcf_{s_3}, \phi_3\rangle$. It follows that $F(\Phi_a)-F(\Phi_b)=\Gamma(\langle cl_{s_1},cl_{s_2},cl_{s_3},cl_{h_1},cl_{h_1+1}\rangle)=D_{cl_{s_1}cl_{s_2}}+D_{cl_{s_2}cl_{s_3}}+D_{cl_{h_1}cl_{h_1+1}}-D_{cl_{s_1}cl_{s_3}}-D_{cl_{h_1}cl_{s_2}}-D_{cl_{s_2}cl_{h_1+1}}$.

(2) Move $Tcf_{s_3}$ to be between $Tcf_{h_2}$ and $Tcf_{h_2+1}$ and convert $\phi_1$ into a new sequence $\phi_4$, where $Tcf_{h_2}, Tcf_{h_2+1}\in \phi_1$.
Let $\Phi_c=\langle \phi_4, Tcf_{s_2}, \phi_2\rangle$. It follows that $F(\Phi_a)-F(\Phi_c)=\Gamma(\langle cl_{s_2},cl_{s_3}, cl_{s_3+1},cl_{h_2},cl_{h_2+1}\rangle)=D_{cl_{s_2}cl_{s_3}}+D_{cl_{s_3}cl_{s_3+1}}+D_{cl_{h_2}cl_{h_2+1}}-D_{cl_{s_2}cl_{s_3+1}}-D_{cl_{h_2}cl_{s_3}}-D_{cl_{s_3}cl_{h_2+1}}$.
}\end{lemma}

\begin{lemma}{\rm \label{lemma2.2} Consider a takeoff aircraft sequence $\Phi_a=\langle \phi_1, \phi_2\rangle$, where $\phi_1=\langle Tcf_{1}, Tcf_{2}, \cdots, Tcf_{s_1}\rangle$, $\phi_2=\langle Tcf_{s_1+1}, Tcf_{s_1+2},\cdots, Tcf_{n}\rangle$, $cl_{1}\geq cl_{2}\geq \cdots \geq cl_{s_1}=k_c$, $cl_{s_1+1}\geq cl_{s_1+2}\geq \cdots \geq cl_{n}$, and $s_1$ and $k_c$ are positive integers. Merge $\phi_1, \phi_2$ to form a class-monotonically-decreasing sequence $\Phi_b$. Suppose that Assumptions \ref{ass2.3.1}-\ref{ass2.4} hold for each takeoff aircraft sequence, $\Theta_a=\{k_{1a},k_{2a},\cdots,k_{\eta_ca}\}$ is the set of all possible aircraft classes in $\phi_1$, and $\Theta_b=\{k_{1b},k_{2b},\cdots,k_{\eta_fb}\}$ is the set of all possible aircraft classes in $\phi_2$ for two positive integers $\eta_c$ and $\eta_f$, where $k_{1a}<k_{2a}<\cdots<k_{\eta_ca}\leq \eta$, $k_{1b}<k_{2b}<\cdots<k_{\eta_fb}\leq \eta$. The following statements hold.

(1) If $k_c\neq h_0=k_{ia}\in \Theta_a$ and $h_0=k_{jb}\in \Theta_b$, $f_{h_0}(\Phi_b)-f_{h_0}(\Phi_a)=D_{h_0h_0}+D_{h_0h_b}-D_{h_0k_{(i-1)a}}-D_{h_0k_{(j-1)b}}$, where $h_b=\max\{k_{(i-1)a}, k_{(j-1)b}\}$.

(2) If $k_c\neq h_0=k_{ia}\in \Theta_a$ and $h_0\notin \Theta_b$, $f_{h_0}(\Phi_b)-f_{h_0}(\Phi_a)=D_{h_0h_b}-D_{h_0k_{(i-1)a}}$, where $h_b$ is the largest integer smaller than $h_0$ in $\Theta_a\cup \Theta_b$.

(3) If $k_c\neq h_0=k_{jb}\in \Theta_b$ and $h_0\notin \Theta_a$, $f_{h_0}(\Phi_b)-f_{h_0}(\Phi_a)=D_{h_0h_b}-D_{h_0k_{(j-1)b}}$, where $h_b$ is the largest integer smaller than $h_0$ in $\Theta_a\cup \Theta_b$.}\end{lemma}

\begin{corollary}{\rm Under the conditions in Lemma \ref{lemma2.2}, the following statements hold.

(1.1) If $k_c\neq h_0=2\in \Theta_a \cap\Theta_b$, and $1\in \Theta_a \cap\Theta_b$, $f_{h_0}(\Phi_b)-f_{h_0}(\Phi_a)=0$.

(1.2) If $k_c\neq h_0=2\in \Theta_a \cap\Theta_b$, and $1\notin \Theta_a \cap\Theta_b$, $f_{h_0}(\Phi_b)-f_{h_0}(\Phi_a)={4T_0/3}$.

(2.1) If $k_c\neq h_0$, $2<h_0\in \Theta_a \cap\Theta_b$, and $\Theta_b-\{h_0, h_0+1, \cdots, \eta\}\neq \emptyset$, $f_{h_0}(\Phi_b)-f_{h_0}(\Phi_a)=D_{h_0h_0}-T_0$.

(2.2) If $k_c\neq h_0$, $2<h_0\in \Theta_a \cap\Theta_b$, and $\Theta_b-\{h_0, h_0+1, \cdots, \eta\}=\emptyset$, $f_{h_0}(\Phi_b)-f_{h_0}(\Phi_a)=D_{h_0h_0}$.}\end{corollary}

In Lemmas \ref{lemma2.31}-\ref{lemma2.1}, we discuss some typical cases for breakpoint aircraft which yields local minimum points.

\begin{lemma}\label{lemma2.31}{\rm  Consider a  takeoff sequence $\Phi_a=\langle\phi_1, \phi_2\rangle$, where $\phi_1=\langle Tcf_1, Tcf_2, \cdots, Tcf_{n_1}\rangle$, $\phi_2=\langle Tcf_{n_1+1},$ $ Tcf_{{n_1}+2}, \cdots, Tcf_{n_2}\rangle$, $\phi_1$ and $\phi_2$ are both class-monotonically-decreasing sequences, $cl_{n_1}<cl_{n_1+1}$, and
$n_1, n_2$ are two positive integers. Merge the aircraft sequences $\phi_1$ and $\phi_2$ to form a new class-monotonically-decreasing sequence $\Phi_b$. Suppose that Assumptions \ref{ass2.3.1}-\ref{ass2.4} hold for each takeoff aircraft sequence. The following statements hold.


(1) Suppose that $(cl_{n_1}, cl_{n_1+1})=(\rho_2-1, \rho_2)$. Then $F(\Phi_a)-F(\Phi_b)=0$.

(2) Suppose that $(cl_{n_1}, cl_{n_1+1}, j)=(3,4,3)$ for some $j\in\{n_1+2,n_1+3,\cdots, n_1+n_2\}$. Then $F(\Phi_a)-F(\Phi_b)=0$.

(3) Suppose that $(cl_{n_1}, cl_{n_1+1}, j)=(\eta-1,\eta,\eta)$ for some $j\in\{1, 2, \cdots, n_1-1\}$. Then $F(\Phi_a)-F(\Phi_b)=0$.

(4) Suppose that $cl_{n_2}\leq 2$ and $(cl_{n_1}, cl_{n_1+1},j)=(2,3,3)$ for some $j\in\{1, 2, \cdots, n_1-1\}$. Then {$F(\Phi_a)-F(\Phi_b)=0$}.

(5) Suppose that all of the supposition conditions in (1)-(4) do not hold. Then $F(\Phi_a)-F(\Phi_b)>0$.
}\end{lemma}

\noindent{Proof}: (1)-(4) The results are obvious and hence the proofs are omitted.

(5) From Assumptions \ref{ass2.3.1} and \ref{ass2.3.2}, $D_{11}=D_{22}=D_{21}$, $D_{33}=4T_0/3$ and $D_{\eta\eta}=4T_0/3$.

If $cl_{n_1}=1$, from Assumption \ref{ass2.3.41}, by simple calculations for the cases when $cl_{n_1+1}=2$, when $2<cl_{n_1+1}<\eta$ and when $cl_{n_1+1}=\eta$, it follows that $f_1(\Phi_a)+f_2(\Phi_a)+f_3(\Phi_a)+f_{\eta}(\Phi_a)-[f_1(\Phi_b)+f_2(\Phi_b)+f_3(\Phi_b)+f_{\eta}(\Phi_b)]>0$. Note that $f_i(\Phi_a)=f_i(\Phi_a)$ for $i\neq 1, 2, 3, \eta$. Hence $F(\Phi_a)-F(\Phi_b)>0$.

If $cl_{n_1}=2$, from Assumption \ref{ass2.3.41}, by simple calculations for the cases when $2<cl_{n_1+1}<\eta$ and when $cl_{n_1+1}=\eta$, it follows that $f_2(\Phi_a)+f_3(\Phi_a)+f_{\eta}(\Phi_a)-[f_2(\Phi_b)+f_3(\Phi_b)+f_{\eta}(\Phi_b)]>0$ and hence $F(\Phi_a)-F(\Phi_b)>0$.

If $cl_{n_1}=3$, from Assumption \ref{ass2.3.41}, by simple calculations for the cases when $3<cl_{n_1+1}<\eta$ and when $cl_{n_1+1}=\eta$, it follows that $f_3(\Phi_a)+f_{\eta}(\Phi_a)-[f_3(\Phi_b)+f_{\eta}(\Phi_b)]>0$ and hence $F(\Phi_a)-F(\Phi_b)>0$.

If $cl_{n_1}<3$, from Assumptions \ref{ass2.3.41}, by simple calculations for the cases when $cl_{n_1}<cl_{n_1+1}<\eta$ and when $cl_{n_1+1}=\eta$, it follows that $f_{cl_{n_1}}(\Phi_a)+f_{\eta}(\Phi_a)-[f_{cl_{n_1}}(\Phi_b)+f_{\eta}(\Phi_b)]>0$ and hence $F(\Phi_a)-F(\Phi_b)>0$.

Summarizing the above analysis, this lemma is proved.

\begin{lemma}\label{lemma2.1}{\rm Consider a takeoff aircraft sequence $\Phi_a=\langle\phi_1, \phi_2\rangle$, where $\phi_1=\langle Tcf_1, Tcf_2, $ $\cdots, Tcf_{s_1}\rangle$, $\phi_2=\langle Tcf_{s_1+1}, Tcf_{s_1+2}, \cdots, Tcf_{n}\rangle$, both of the sequences $\phi_1$ and $\phi_2$ are class-monotonically-decreasing sequences, $cl_{s_1}<cl_{s_1+1}$, and $s_1<n-1$ is a positive integer. Suppose that $\phi_2$ contains $s_2$ aircraft of class $h_1$, $Tcf_i\in \phi_2$, and $cl_i=h_1<cl_{s_1}$. Generate a new sequence $\Phi_b$ by moving the aircraft $Tcf_i$ to be between $Tcf_{s_1}$ and $Tcf_{s_1+1}$. Suppose that Assumptions \ref{ass2.3.1}-\ref{ass2.4} hold for each aircraft sequence. The following statements hold.

(1) If $h_1=1$, $F(\Phi_a)-F(\Phi_b)\leq0$.

{(2) If $h_1=2$, $cl_{n-1}\leq 2$ and $\langle cl_{s_1}, cl_{s_1+1}\rangle\in E_1$, $F(\Phi_a)-F(\Phi_b)=0$.}

(3) If $h_1=2$, $cl_{n-1}>2$ and $\langle cl_{s_1}, cl_{s_1+1}\rangle\in E_1$, $F(\Phi_a)-F(\Phi_b)=-T_0/3$. If $h_1=2$ and $\langle cl_{s_1}, cl_{s_1+1}\rangle\notin E_1$,  $F(\Phi_a)-F(\Phi_b)<0$.

(4) If $h_1=3$, $s_2>1$ and $\langle cl_{s_1}, cl_{s_1+1}\rangle\in E_{20}$, $F(\Phi_a)-F(\Phi_b)=T_0/3$.
 If $h_1=3$, $s_2=1$ and $\langle cl_{s_1}, cl_{s_1+1}\rangle\in E_{20}$, $F(\Phi_a)-F(\Phi_b)=0$.

(5) If $h_1=3$, $s_2>1$ and $\langle cl_{s_1}, cl_{s_1+1}\rangle\in E_{21}$, $F(\Phi_a)-F(\Phi_b)=0$.
 If $h_1=3$, $s_2=1$ and $\langle cl_{s_1}, cl_{s_1+1}\rangle\in E_{21}$, $F(\Phi_a)-F(\Phi_b)=T_0/3$.

(6) If $h_1=3$ and $\langle cl_{s_1}, cl_{s_1+1}\rangle\notin E_{20}\cup E_{21}$, $F(\Phi_a)-F(\Phi_b)<0$.

(7) If $h_1>3$, $F(\Phi_a)-F(\Phi_b)<0$.

(8) If $h_1=\rho_2-1$ and $\langle cl_{s_1}, cl_{s_1+1}\rangle\in E_{30}$, $F(\Phi_a)-F(\Phi_b)=T_0/3$.

(9) If $h_1=\rho_2$ and $\langle cl_{s_1}, cl_{s_1+1}\rangle \in E_{31}$, $F(\Phi_a)-F(\Phi_b)=T_0/3$.
}\end{lemma}
\noindent{Proof:} (1) Note that $h_1<cl_{s_1}<cl_{s_1+1}$. From Lemma \ref{lemma2.111}, $D_{h_1cl_{s_1+1}}-D_{cl_{s_1}cl_{s_1+1}}>0.5T_0$. By simple calculations for the cases when $D_{cl_{s_1}h_1}=T_0$ and when $D_{cl_{s_1}h_1}=1.5T_0$, it follows that $F(\Phi_a)-F(\Phi_b)<0$.

(2) When $h_1=2$ and $\langle cl_{s_1}, cl_{s_1+1}\rangle\in E_1$, from Lemma \ref{lemma2.111}, $D_{2cl_{s_1+1}}-D_{cl_{s_1}cl_{s_1+1}}=T_0/3$. Note that $cl_{n-1}\leq 2$. Then there are at least two aircraft of class no larger than $2$ in $\phi_2$. By simple calculations, $F(\Phi_a)-F(\Phi_b)=0$.

(3) When $h_1=2$ and $\langle cl_{s_1}, cl_{s_1+1}\rangle\in E_1$, from Lemma \ref{lemma2.111}, $D_{2cl_{s_1+1}}-D_{cl_{s_1}cl_{s_1+1}}=T_0/3$. Note that $cl_{n-1}> 2$. Then there is at most one aircraft of class no larger than $2$ in $\phi_2$. By simple calculations, $F(\Phi_a)-F(\Phi_b)=\Gamma(cl_{i}, h_1, cl_{i+1},cl_{s_1},cl_{s_1+1})=-T_0/3$. When $h_1=2$ and $\langle cl_{s_1}, cl_{s_1+1}\rangle\in E_1$, from Lemma \ref{lemma2.111}, $D_{2cl_{s_1+1}}-D_{cl_{s_1}cl_{s_1+1}}>0.5T_0$. By simple calculations, $F(\Phi_a)-F(\Phi_b)<0$.

The statements (4)-(9) can be proved by using similar approaches in (1)-(3).

In Theorem \ref{lemma2.3s}, we give a rule to calculate the objective function $F(\cdot)$ based on a standard class-monotonically-decreasing sequence. As a matter of fact, combing Lemmas \ref{lemma2.2}-\ref{lemma2.1}, Theorem \ref{lemma2.3s} and other special properties/constraints, the optimal sequence can be studied for the optimization problem (\ref{optim1}).

\begin{theorem}\label{lemma2.3s}{\rm Consider a takeoff sequence $\Phi_a=\langle\phi_1, \phi_2, \cdots, \phi_s\rangle$ for a positive integer $s$, where $\phi_i=\langle Tcf_{i1}, Tcf_{i2}, \cdots,$ $ Tcf_{ic_i}\rangle$ is a class-monotonically-decreasing sequence for some positive integer $c_i$ and all $i\in \{1, 2, \cdots, s\}$, and $cl_{jc_j}<cl_{(j+1)1}$ for all $j\in \{1, 2, \cdots, s-1\}$. Merge the aircraft sequences $\Phi_a$ to form a new class-monotonically-decreasing sequence $\Phi_b$. Suppose that Assumptions \ref{ass2.3.1}-\ref{ass2.4} hold for each sequence. Then, $F(\Phi_a)=F(\Phi_b)+D_{d1}-D_{cl_{sc_s}1}-\mathrm{sgn}(\sum_{i=1}^s \mathrm{Exs}_{\rho_1}(\phi_i))[\sum_{i=1}^s \mathrm{Exs}_{\rho_1}(\phi_i)-1]T_0/3-\mathrm{sgn}(\sum_{i=1}^s \mathrm{Exs}_{\rho_2}(\phi_i))[\sum_{i=1}^s\mathrm{Exs}_{\rho_2}(\phi_i)-1]T_0/3+\sum_{i=1}^{s-1}[D_{cl_{ic_i}cl_{(i+1)1}}-D_{cl_{ic_i}1}]$, where $d$ denotes the class of the last aircraft of $\Phi_b$.
}\end{theorem}

\begin{theorem}\label{theorem2.211bb}{\rm Consider a group of takeoff aircraft $\{Tcf_1, Tcf_2, \cdots, Tcf_n\}$, where $cl_1 \geq cl_2 \geq cl_3 \geq \cdots \geq cl_n$. Suppose that Assumptions \ref{ass2.3.1}-\ref{ass2.4} hold for each takeoff aircraft sequence.
Then the optimization problem (\ref{optim1}) can be solved if the takeoff aircraft sequence is taken as $\phi_0=\langle Tcf_1, Tcf_2, \cdots, Tcf_n \rangle$.}\end{theorem}

As discussed in Sec. \ref{seclanding}, here we also introduce some class-sequence sets as an example to study the local minimum points for the optimization problem (\ref{optim1}).

Let $\Upsilon_0=\{\langle i_1, i_2, i_3, i_4, i_5\rangle\mid {i_1}\geq {i_2}\geq {i_3}, {i_4}\geq {i_5}, i_1, i_2, i_3, i_4, i_5\in \mathcal{I}\}$, $\Upsilon_1=\{\langle i_1, i_2, i_3, i_4, i_5\rangle\mid {i_1}\geq {i_2}, {i_2}<{i_3}, {i_4}\geq {i_5}, i_1, i_2, i_3, i_4, i_5\in \mathcal{I}\}$, $\Upsilon_2=\{\langle i_1, i_2, i_3, i_4, i_5\rangle\mid {i_1}<{i_2}, {i_2}\geq {i_3}, {i_4}\geq {i_5}, i_1, i_2, i_3, i_4, i_5\in \mathcal{I}\}$, $\Upsilon_3=\{\langle i_1, i_2, i_3, i_4, i_5\rangle\mid {i_1}\geq {i_2}, {i_2}<{i_3}, {i_4}< {i_5}, i_1, i_2, i_3, i_4, i_5\in \mathcal{I}\}$, $\Upsilon_4=\{\langle i_1, i_2, i_3, i_4, i_5\rangle\mid {i_1}<{i_2}, {i_2}\geq {i_3}, {i_4}<{i_5}, i_1, i_2, i_3, i_4, i_5\in \mathcal{I}\}$ and $\Upsilon_5=\{\langle i_1, i_2, i_3, i_4, i_5\rangle\mid {i_1}\geq {i_2}\geq {i_3}, {i_4}< {i_5}, i_1, i_2, i_3, i_4, i_5\in \mathcal{I}\}$.

\begin{definition}{\rm (1) Let $\Lambda_0\subseteq \Upsilon_0$ be a class-sequence set such that $\Gamma(\langle i_1, i_2, i_3, i_4, i_5\rangle)\geq 0$ for any $\langle i_1, i_2, i_3, i_4, i_5\rangle\in \Lambda_0$, and $\Gamma(\langle i_1, i_2, i_3, i_4, i_5\rangle)<0$ for any $\langle i_1, i_2, i_3, i_4, i_5\rangle\in \Upsilon_0-\Lambda_0$.

(2) Let $\Lambda_1\subseteq \Upsilon_1$ be a class-sequence set such that $\Gamma(\langle i_1, i_2, i_3, i_4, i_5\rangle)\geq 0$ for any $\langle i_1, i_2, i_3, i_4, i_5\rangle\in \Lambda_1$, and $\Gamma(\langle i_1, i_2, i_3, i_4, i_5\rangle)<0$ for any $\langle i_1, i_2, i_3, i_4, i_5\rangle\in \Upsilon_1-\Lambda_1$.

(3) Let $\Lambda_2\subseteq \Upsilon_2$ be a class-sequence set such that $\Gamma(\langle i_1, i_2, i_3, i_4, i_5\rangle)\geq 0$ for any $\langle i_1, i_2, i_3, i_4, i_5\rangle\in \Lambda_2$, and $\Gamma(\langle i_1, i_2, i_3, i_4, i_5\rangle)<0$ for any $\langle i_1, i_2, i_3, i_4, i_5\rangle\in \Upsilon_2-\Lambda_2$.

(4) Let $\Lambda_3\subseteq \Upsilon_3$ be a class-sequence set such that $\Gamma(\langle i_1, i_2, i_3, i_4, i_5\rangle)\geq 0$ for any $\langle i_1, i_2, i_3, i_4, i_5\rangle\in \Lambda_3$, and $\Gamma(\langle i_1, i_2, i_3, i_4, i_5\rangle)<0$ for any $\langle i_1, i_2, i_3, i_4, i_5\rangle\in \Upsilon_3-\Lambda_3$.

(5) Let $\Lambda_4\subseteq \Upsilon_4$ be a class-sequence set such that $\Gamma(\langle i_1, i_2, i_3, i_4, i_5\rangle)\geq 0$ for any $\langle i_1, i_2, i_3, i_4, i_5\rangle\in \Lambda_4$, and $\Gamma(\langle i_1, i_2, i_3, i_4, i_5\rangle)<0$ for any $\langle i_1, i_2, i_3, i_4, i_5\rangle\in \Upsilon_4-\Lambda_4$.

(6) Let $\Lambda_5\subseteq \Upsilon_5$ be a class-sequence set such that $\Gamma(\langle i_1, i_2, i_3, i_4, i_5\rangle)\geq 0$ for any $\langle i_1, i_2, i_3, i_4, i_5\rangle\in \Lambda_5$, and $\Gamma(\langle i_1, i_2, i_3, i_4, i_5\rangle)<0$ for any $\langle i_1, i_2, i_3, i_4, i_5\rangle\in \Upsilon_5-\Lambda_5$.}\end{definition}

\section{Scheduling problem of landing and  takeoff aircraft on a same runway}

In this section, the scheduling problem of mixed takeoff and landing aircraft on a same way is discussed. Let $D_T$ denote the minimum  separation time between a landing aircraft and a leading takeoff aircraft. Let $T_D$ denote the minimum separation time between a takeoff aircraft and a leading landing aircraft.
As defined previously, $Y_{ij}$ is used to unifiedly denote the minimum separation time between any given two aircraft $Tcf_i$ and $Tcf_j$. Specifically, when a landing aircraft $Tcf_j$ and a leading takeoff aircraft $Tcf_i$ are considered, $Y_{ij}=D_T$; when a takeoff aircraft $Tcf_j$ and a leading landing aircraft $Tcf_i$ are considered, $Y_{ij}=T_D$; when two consecutive landing aircraft $Tcf_i$ and $Tcf_j$ are considered, $Y_{ij}=T_{cl_icl_j}$; and when two consecutive takeoff aircraft $Tcf_i$ and $Tcf_j$ are considered, $Y_{ij}=D_{cl_icl_j}$.

\begin{assumption}\label{ass4.11}{\rm Suppose that $T_0\leq T_D<1.5T_0$ and $T_0\leq D_T<1.5T_0$.}\end{assumption}

\begin{assumption}\label{ass4.21a}{\rm Consider an aircraft sequence $\phi=\langle Tcf_1, Tcf_2, \cdots, Tcf_{n}\rangle$. Suppose that the aircraft $Tcf_{j_1}$ is relevant to $Tcf_{j_0}$, the aircraft $Tcf_{j_2}$ is relevant to $Tcf_{j_1}$, $Y_{j_0j_1}\geq T_D+D_T$ and $Y_{j_1j_2}\geq T_D+D_T$. If the aircraft $Tcf_{j_3}$ is relevant to $Tcf_{j_2}$, then $Y_{j_2j_3}<T_D+D_T$.}\end{assumption}

From the conditions that $Y_{j_0j_1}\geq T_D+D_T\geq 2T_0$ and $Y_{j_1j_2}\geq T_D+D_T\geq 2T_0$, it can be obtained the aircraft $Tcf_{j_0}$, $Tcf_{j_1}$ and $Tcf_{j_2}$ are all landing or takeoff aircraft and $cl_{j_0}<cl_{j_1}<cl_{j_2}$. Assumption \ref{ass4.21a} means that the separation time between $Tcf_{j_2}$ and $Tcf_{j_3}$ is smaller than
$T_D+D_T$, which is consistent with the minimum separation time standards at Heathrow Airport and for the RECAT-EU system (See, e.g., Table \ref{tab:Taking-off_separation}).

From the previous sections, when the leading and the trailing aircraft are both landing or takeoff aircraft, the minimum separation time between an aircraft and its leading aircraft is only relevant to their own classes. In contrast, when the takeoff and landing aircraft are simultaneously considered, the minimum separation time between an aircraft and its leading aircraft might be related to not only their own classes, but also the classes of the aircraft ahead of them.

In Theorem \ref{theorem4.11b}, we study the aircraft relevance and show that when the operation times of the aircraft have no constraints, each trailing aircraft should be relevant to at least one of its nearest landing and takeoff aircraft ahead to ensure the optimality of the sequence. Theorem \ref{theorem3.4} shows that the optimal value of the objective function $F(\cdot)$ can be calculated along the path from the last aircraft to the first aircraft. 

\begin{definition} {\rm(Path) Consider a sequence $\phi=\langle Tcf_1, Tcf_2, \cdots, Tcf_{n}\rangle$. For any given aircraft $Tcf_i$ and $Tcf_j$, if there exists an aircraft subsequence $\langle Tcf_i^0, Tcf_i^1, \cdots, Tcf_i^{\rho}\rangle$ for some positive integer $\rho>0$ such that $Tcf_i^0=Tcf_i$, $Tcf_i^\rho=Tcf_j$ and each aircraft $Tcf^h_i$ is relevant to
aircraft $Tcf^{h-1}_i$, $h=1, 2, \cdots, \rho$, then the sequence $\langle Tcf_i^0, Tcf_i^1, \cdots, Tcf_i^{\rho}\rangle$ is said to be a path from the aircraft $Tcf_j$ to the aircraft $Tcf_i$. It is assumed by default that each aircraft has a path to itself.}\end{definition}

\begin{theorem}\label{theorem4.11b}{\rm  Consider an aircraft sequence $\phi=\langle Tcf_1, Tcf_2, \cdots, Tcf_{n}\rangle$. Suppose that $Tf_j=[t_0, +\infty]$ for all $j$.
The following statements hold.

(1) Suppose that aircraft $Tcf_i$ is a takeoff (landing) aircraft and aircraft $Tcf_j$ is a landing (takeoff) aircraft, and $0<j-i\leq 3$. If aircraft $Tcf_j$ is relevant to aircraft $Tcf_i$, then $j=i+1$.

(2) Suppose that aircraft $Tcf_i$ and $Tcf_j$ are both landing (takeoff) aircraft. If $j=i+1$ and aircraft $Tcf_j$ is relevant to aircraft $Tcf_k$, then $k=i$.

(3) Suppose that aircraft $Tcf_j$ is relevant to aircraft $Tcf_i$. If $Tcf_i$ is a landing (takeoff) aircraft, then aircraft $Tcf_{i+1}$, $Tcf_{i+2}$, $\cdots$, $Tcf_{j-1}$ are all takeoff (landing) aircraft.

(4) Suppose that $j-i>3$. The aircraft $Tcf_j$ is not relevant to aircraft $Tcf_i$.

(5) Suppose that $j-i=3$. If the aircraft $Tcf_j$ is relevant to aircraft $Tcf_i$, the aircraft $Tcf_j$ is also relevant to aircraft $Tcf_{j-1}$.
}\end{theorem}

\noindent{Proof:} (1) Since $Y_{ij}=D_T<1.5T_0$ from Assumption \ref{ass4.11}, if $j\geq i+2$, then $S_{ij}\geq 2T_0>Y_{ij}$ and aircraft $Tcf_j$ is not relevant to aircraft $Tcf_i$. If aircraft $Tcf_j$ is relevant to aircraft $Tcf_i$, then $j=i+1$.

(2) Note that $D_T<1.5T_0$ and $T_D<1.5T_0$ from Assumption \ref{ass4.11}. If $k\neq i$, it follows that $k<i$, and $Tcf_k$ is a landing (takeoff) aircraft.  From Lemma \ref{lemma11s}, $S_{kj}\geq S_{ij}+S_{ki}\geq Y_{i{j}}+Y_{k{i}}>Y_{k{j}}$, implying that aircraft $Tcf_j$ is not relevant to aircraft $Tcf_k$. Therefore, $k=i$.

(3)-(4) These two is obvious and hence their proofs are omitted.

(5) It is clear that the aircraft $Tcf_i$ and $Tcf_j$ are both landing (or takeoff) aircraft, and $Y_{ij}=3T_0$. This implies that $Y_{(j-1)j}=T_0$ and $Y_{(j-2)j}=2T_0$. From (3),
the aircraft $Tcf_{j-1}$ and $Tcf_{j-2}$ are both takeoff (or landing) aircraft. Thus, the aircraft $Tcf_j$ is also relevant to aircraft $Tcf_{j-1}$ and not relevant to aircraft $Tcf_{j-2}$.

\begin{remark}{\rm When the scheduling problem of mixed landing and takeoff aircraft on a same runway, each aircraft might be relevant to two aircraft ahead: one takeoff aircraft and one landing aircraft, and is not relevant to other aircraft.}\end{remark}

\begin{theorem}\label{theorem3.4}{\rm Consider an aircraft sequence $\phi=\langle Tcf_1, Tcf_2, \cdots, Tcf_{n}\rangle$. Suppose that the orders of all the aircraft in $\phi$ is fixed.                 

(1) The optimization problem (\ref{optim1}) is solved if and only if there is a path from the aircraft $Tcf_n$ to the aircraft $Tcf_1$.

(2) If the aircraft $Tcf_2$ is relevant to the aircraft $Tcf_1$ and the aircraft $Tcf_{i+2}$ is relevant to the aircraft $Tcf_i$ or $Tcf_{i+1}$, $i=1,2, \cdots, n-2$, there is a path from the aircraft $Tcf_n$ to the aircraft $Tcf_1$.
}\end{theorem}

\noindent{Proof:} (1) Sufficiency is obvious and we only discuss the necessity of this theorem. Suppose that there is no path
from the aircraft $Tcf_n$ to the aircraft $Tcf_1$,
when the optimization problem (\ref{optim1}) is solved. Construct an aircraft set $\Theta_1$ such that there is a path from the aircraft $Tcf_n$ to the aircraft $Tcf_i$ for any $Tcf_i\in \Theta_1$ and there is no path from the aircraft $Tcf_n$ to the aircraft $Tcf_i$ for any $Tcf_i\in \{Tcf_1, Tcf_2, \cdots, Tcf_{n}\}-\Theta_1$. Clearly, $Tcf_1\notin \Theta_1$. Shift the takeoff or landing times of all the aircraft in $\Theta_1$ as a whole to an earlier time instant with their relative landing or takeoff times unchanged while keeping necessary minimum separation times with the aircraft in $\{Tcf_1, Tcf_2, \cdots, Tcf_{n}\}-\Theta_1$. Such an operation can be realized because the aircraft in $\Theta_1$ are not relevant to any aircraft in $\{Tcf_1, Tcf_2, \cdots, Tcf_{n}\}-\Theta_1$. After this operation, the value of the objective function $F(\phi, Sr(\phi))$ becomes smaller, which yields a contraction. Therefore, there is a path from the aircraft $Tcf_n$ to the aircraft $Tcf_1$.

 (2) From Theorem \ref{theorem4.11b} and the definition of path, if the aircraft $Tcf_2$ is relevant to the aircraft $Tcf_1$ and the aircraft $Tcf_{i+2}$ is relevant to the aircraft $Tcf_i$ or $Tcf_{i+1}$, $i=1,2, \cdots, n-2$, there is a path from the aircraft $Tcf_n$ to the aircraft $Tcf_1$.

\begin{assumption}\label{ass3.4}{\rm Consider the aircraft sequence $\phi=\langle Tcf_1, Tcf_2, \cdots, Tcf_{n}\rangle$. Suppose that the aircraft $Tcf_2$ is relevant to aircraft $Tcf_1$, and aircraft $Tcf_{i+2}$ is relevant to aircraft $Tcf_{i+1}$ or aircraft $Tcf_i$, $i=1,2, \cdots, n-2$.}\end{assumption}

When there is a path from the last aircraft $Tcf_n$ to the first aircraft $Tcf_1$, there might be some aircraft that do not belong to the path and have no relevant aircraft. We can adjust the landing or takeoff times of these aircraft without affecting other aircraft so as to make the aircraft sequence satisfy the condition in Assumption \ref{ass3.4}.

\begin{assumption}\label{ass3.411}{\rm Suppose that $Tf_j=[t_0, +\infty]$ for all $j$, and Assumptions \ref{ass1.3.1}-\ref{ass1.3.4}, \ref{ass2.3.1}-\ref{ass2.3.41} and \ref{ass4.11}-\ref{ass3.4} hold.}\end{assumption}

When the separation time of two consecutive takeoff/landing aircraft is large, by adding one aircraft with different operation task, the total operation time of the aircraft can be decreased. In the following, we make discussions about this issue.

\begin{theorem}\label{theorem4.3}{\rm Consider an aircraft sequence $\phi_0=\langle Tcf_1, Tcf_2, \cdots, Tcf_{n}\rangle$ and the sequence $\phi_1$ generated by moving the aircraft $Tcf_{i_0}$ to be between
$Tcf_{j_0}$ and $Tcf_{j_0+1}$, where aircraft $Tcf_{j_0}$ and $Tcf_{j_0+1}$ are both takeoff (landing) aircraft, and aircraft $Tcf_{i_0}$ is a landing (takeoff) aircraft. Suppose that Assumption \ref{ass3.411} holds. The following statements hold.

(1)$ F(\phi_1)-F(\phi_0)=S_{j_0i_0}(\phi_1)+S_{i_0(j_0+1)}(\phi_1)+S_{(i_0-1)(i_0+1)}(\phi_1)
-S_{(i_0-1)i_0}(\phi_0)-S_{i_0(i_0+1)}(\phi_0)-S_{j_0(j_0+1)}(\phi_0).$

(2) Suppose that in $\phi_0,\phi_1$, the aircraft $Tcf_{i_0}$, $Tcf_{i_0+1}$ are relevant to their leading aircraft, and the aircraft $Tcf_{j_0+1}$ is relevant to aircraft $Tcf_{j_0}$. It follows that $Y_{{j_0}{(j_0+1)}}\geq S_{j_0i_0}(\phi_1)+S_{i_0(j_0+1)}(\phi_1)=T_D+D_T$, and
 \begin{eqnarray*}\begin{array}{lll}F(\phi_1)-F(\phi_0)=Y_{(i_0-1)(i_0+1)}-Y_{(i_0-1)i_0}-Y_{i_0(i_0+1)}.\end{array}\end{eqnarray*}}\end{theorem}

Let $\mathrm{Det}(\phi)$ be a function such that $\mathrm{Det}(\phi)=1$ when the aircraft in $\phi$ are all landing aircraft, and $\mathrm{Det}(\phi)=0$ when the aircraft in $\phi$ are all takeoff aircraft.
\begin{theorem}\label{lemma3.31s}{\rm  Consider two sequences $\Phi_a=\langle\phi^a_1, \phi^a_2, \cdots, \phi^a_s\rangle$ and $\Phi_b=\langle\phi^b_1, \phi^b_2\rangle$ formed by the same group of aircraft for a positive integer $s$, where $\phi^a_i=\langle Tcf_{i1}, Tcf_{i2}, \cdots, Tcf_{ic_i}\rangle$ is a class-monotonically-decreasing takeoff or landing sequence for some positive integer $c_i$ and all $i\in \{1, 2, \cdots, s\}$, each aircraft $Tcf_{jc_j}$ is a takeoff-landing or landing-takeoff transition aircraft for all $j\in \{1, 2, \cdots, s-1\}$, $\phi^b_1$ is a class-monotonically-decreasing landing (takeoff) sequence with $Tcf_{n_1}$ as its last aircraft, and $\phi^b_2$ is a class-monotonically-decreasing takeoff (landing) sequence with $Tcf_{n_2}$ and $Tcf_{n_3}$ as its first and last aircraft.
Suppose that Assumption \ref{ass3.411} holds for each sequence.
  Then, $F(\Phi_a)=F(\Phi_b)+D_{cl_{n_1}1}-Y_{n_1n_2}+D_{cl_{n_3}1}-D_{cl_{sc_s}1}-\sum_{k=\rho_1, \rho_2} \mathrm{sgn}(p_k)
  (p_k-1)T_0/3-\sum_{k=\rho_1, \rho_2}\mathrm{sgn}(q_k)(q_k-1)\delta+\sum_{i=1}^{s-1}[Y_{(ic_i)((i+1)1)}-Y_{(ic_i)1}]$,
  where $p_k=\sum_{i=1}^s [1-\mathrm{Det}(\phi^a_i)]\mathrm{Exs}_{k}(\phi^a_i)$ and $q_k=\sum_{i=1}^s \mathrm{Det}(\phi^a_i)\mathrm{Exs}_{k}(\phi^a_i)$.
}\end{theorem}

\begin{definition}{\rm Consider the aircraft sequence $\langle Tcf_i, Tcf_j\rangle$. If the aircraft $Tcf_i$ is a takeoff aircraft, and the aircraft $Tcf_j$ is a landing aircraft, it is said that the aircraft sequence forms a takeoff-landing transition at the aircraft $Tcf_i$. If the aircraft $Tcf_i$ is a landing aircraft, and the aircraft $Tcf_j$ is a takeoff aircraft, it is said that the aircraft sequence forms a landing-takeoff transition at the aircraft $Tcf_i$.}\end{definition}

\begin{remark}{\rm In Theorem \ref{lemma3.31s}, what is different from the landing sequence and the takeoff sequence is that there is a landing-takeoff transition or a takeoff-landing transition at aircraft $Tcf_{n_1}$ in $\Phi_a$. Due to this reason, the term $D_{cl_{n_1}1}-Y_{n_1n_2}$ is included in the relationship between $F(\Phi_a)$ and $F(\Phi_b)$.}\end{remark}

\begin{remark}{\rm Theorem \ref{lemma3.31s} gives a rule to calculate the value of the objective function $F(\cdot)$ for scheduling problem of mixed takeoff and landing aircraft on a same runway.}\end{remark}

\begin{lemma}\label{theorem4.41s}{\rm Consider an aircraft sequence $\phi_1=\langle\phi_{11}, \phi_{12}\rangle$, where $\phi_{11}=\langle Tcf_{1}, Tcf_{2}, Tcf_{3}, Tcf_{4},Tcf_{5}, Tcf_{6}\rangle$, all the aircraft of $\phi_{11}$ are takeoff (landing) aircraft, $\phi_{12}=\langle Tcf_{7}, Tcf_{8}, Tcf_{9}, Tcf_{10}\rangle$, all the aircraft of $\phi_{12}$ are landing (takeoff) aircraft. Suppose that Assumption \ref{ass3.411} holds for each sequence.

Generate a new sequence $\phi_2$ by moving the aircraft $Tcf_{9}$ and $Tcf_{10}$ in $\phi_1$ to be between the aircraft $Tcf_{2}$ and $Tcf_{3}$ and to be between $Tcf_{3}$ and $Tcf_{4}$. Generate another new sequence $\phi_3$ by moving the aircraft $Tcf_{10}$ and $Tcf_{9}$ in $\phi_1$ to be between the aircraft $Tcf_{2}$ and $Tcf_{3}$ and to be between $Tcf_{3}$ and $Tcf_{4}$.

(1) Suppose that $Y_{23},Y_{34}\leq T_D+D_T$. If $Y_{9,{10}}\leq T_D+D_T$ and $Y_{{10},9}\leq T_D+D_T$, then $F(\phi_2)=F(\phi_3)$.

(2) Suppose that $Y_{23},Y_{34}\leq T_D+D_T$. If $Y_{9,{10}}>T_D+D_T$ and $Y_{{10},9}\leq T_D+D_T$, then $F(\phi_2)>F(\phi_3)$.

Generate a new sequence $\phi_4$ by moving the aircraft $Tcf_{9}$ and $Tcf_{10}$ in $\phi_1$ to be between the aircraft $Tcf_{3}$ and $Tcf_{4}$ in the order of $\langle Tcf_{9}, Tcf_{10}\rangle$. Generate another new sequence $\phi_5$ by moving the aircraft $Tcf_{10}$ and $Tcf_{9}$ in $\phi_1$ to be between the aircraft $Tcf_{3}$ and $Tcf_{4}$ in the order of $\langle Tcf_{10}, Tcf_{9}\rangle$.

(3) Suppose that $Y_{34}<3T_0$. If $Y_{9,{10}}<Y_{{10},{9}}$, then $F(\phi_4)<F(\phi_5)$.

Generate a new sequence $\phi_6$ by moving the aircraft $Tcf_{8}$ in $\phi_1$ to be between the aircraft $Tcf_{2}$ and $Tcf_{3}$ and moving the aircraft $Tcf_{9}$ and $Tcf_{10}$  to be between the aircraft $Tcf_{3}$ and $Tcf_{4}$ in the order of $\langle Tcf_{9}, Tcf_{10}\rangle$. Generate a new sequence $\phi_7$ by moving the aircraft $Tcf_{9}$ in $\phi_1$ to be between the aircraft $Tcf_{2}$ and $Tcf_{3}$ and moving the aircraft $Tcf_{8}$ and $Tcf_{10}$  to be between the aircraft $Tcf_{3}$ and $Tcf_{4}$ in the order of $\langle Tcf_{8}, Tcf_{10}\rangle$.

(4) Suppose that $cl_8 > cl_{10}$ and $cl_9 > cl_{10}$. If $Y_{8{9}}\leq T_D+D_T$ and $Y_{{9}8}\leq T_D+D_T$, then $F(\phi_6)=F(\phi_7)$.

(5) Suppose that $cl_8 > cl_{10}$ and $cl_9 > cl_{10}$. If $Y_{8{9}}>T_D+D_T$ and $Y_{{9}8}\leq T_D+D_T$, then $F(\phi_6)>F(\phi_7)$. }\end{lemma}

\begin{remark}{\rm From Lemma \ref{theorem4.41s}, the order of the aircraft classes has a direct impact on the value of the objective function $F(\cdot)$. Generally, compared with class-monotonically-increasing sequence, class-monotonically-decreasing sequence can distinctly reduce the mutual influence between takeoff-landing and landing-takeoff transitions.}\end{remark}

In Theorem \ref{theorem4.5}, we will make an estimation about the impact range when the aircraft sequence is changed. To this end, we first give the following lemma.

\begin{lemma}\label{lemma4.3}{\rm Consider an aircraft sequence $\phi^1_0=\phi^2_0=\langle Tcf_i, Tcf_{i+1}, Tcf_{i+2}\rangle$, where the aircraft $Tcf_i$ and $Tcf_{i+2}$ are both takeoff (landing) aircraft,
the aircraft $Tcf_{i+1}$ is a landing (takeoff) aircraft, and the separation times between the aircraft in $\phi^1_0$ and $\phi^2_0$ are different. Suppose that in $\phi_0^1$ and $\phi_0^2$, $S_i(\phi^1_0)=S_i(\phi^2_0)$ and the aircraft $Tcf_{i+2}$ is relevant to either the aircraft $Tcf_i$ or the aircraft $Tcf_{i+1}$. When
$S_{i(i+1)}(\phi_{0}^1)>S_{i(i+1)}(\phi_{0}^2)$, then $S_{(i+1)(i+2)}(\phi_{0}^1)=S_{(i+1)(i+2)}(\phi_{0}^2)$ or $S_{(i+1)(i+2)}(\phi_{0}^1)<S_{(i+1)(i+2)}(\phi_{0}^2)$. Further, when $S_{i(i+1)}(\phi_{0}^1)>S_{i(i+1)}(\phi_{0}^2)$ and $S_{(i+1)(i+2)}(\phi_{0}^1)<S_{(i+1)(i+2)}(\phi_{0}^2)$, $S_{i(i+2)}(\phi_{0}^2)=Y_{i(i+2)}$, i.e., in $\phi_0^2$, the aircraft $Tcf_{i+2}$ is relevant to the aircraft $Tcf_i$.}\end{lemma}

\noindent{Proof:} If $\max\{S_{(i+1)(i+2)}(\phi_{0}^1), S_{(i+1)(i+2)}(\phi_{0}^2)\}=Y_{{(i+1)}{(i+2)}}$, then in $\phi_0^1$ and $\phi_0^2$, the aircraft $Tcf_{i+2}$ is relevant to the aircraft $Tcf_{i+1}$, and $S_{(i+1)(i+2)}(\phi_{0}^1)=S_{(i+1)(i+2)}(\phi_{0}^2)$.
Suppose that $\max\{S_{(i+1)(i+2)}(\phi_{0}^1), S_{(i+1)(i+2)}(\phi_{0}^2)\}$ $>Y_{{(i+1)}{(i+2)}}$, and $S_{(i+1)(i+2)}(\phi_{0}^1)\geq S_{(i+1)(i+2)}(\phi_{0}^2)$. Clearly, $S_{(i+1)(i+2)}(\phi_{0}^1)> Y_{({i+1})({i+2})}$.
In $\phi_0^1$ and $\phi_0^2$, the aircraft $Tcf_{i+2}$ is relevant to the aircraft $Tcf_i$ or the aircraft $Tcf_{i+1}$. It follows that $S_{i(i+2)}(\phi_{0}^1)=Y_{{i}({i+2})}$. On the other hand, note that $S_{i(i+2)}(\phi_{0}^1)=Y_{{i}({i+2})}>S_{(i+1)(i+2)}(\phi_{0}^2)+S_{i(i+1)}(\phi_{0}^2)$, which contradicts with the definition of relevance. Therefore, $\max\{S_{(i+1)(i+2)}(\phi_{0}^1), S_{(i+1)(i+2)}(\phi_{0}^2)\}=Y_{({i+1})({i+2})}$, or $S_{(i+1)(i+2)}(\phi_{0}^1)<S_{(i+1)(i+2)}(\phi_{0}^2)$.

When $S_{(i+1)(i+2)}(\phi_{0}^1)<S_{(i+1)(i+2)}(\phi_{0}^2)$, if $S_{(i+1)(i+2)}(\phi_{0}^2)=Y_{(i+1)(i+2)}$, then $Y_{(i+1)(i+2)}>S_{(i+1)(i+2)}(\phi_{0}^1)$, which also contradicts with the definition of relevance. Note that in $\phi_0^1$ and $\phi_0^2$, the aircraft $Tcf_{i+2}$ is relevant to the aircraft $Tcf_i$ or the aircraft $Tcf_{i+1}$. It follows that $S_{i(i+2)}(\phi_{0}^2)=Y_{i(i+2)}$.

\begin{theorem}\label{theorem4.5}{\rm Consider two aircraft sequences $\phi_0=\langle \phi_{01}, \phi_{02}\rangle$ and $\phi_1=\langle \phi_{11}, \phi_{12}\rangle$, where $\phi_{01}, \phi_{02}, \phi_{11}, \phi_{12}$, respectively, denote the subsequences of $\phi_0$ and $\phi_1$, $\phi_{01}\neq \phi_{11}$, and $\phi_{02}=\phi_{12}=\langle Tcf_{b_0}, Tcf_{b_0+1}, \cdots, Tcf_{b_0+m}\rangle$ for two positive integers $b_0>0$ and $m>0$.  Suppose that Assumption \ref{ass3.411} holds.

 (1) If there is an integer $0\leq j_0<m$ such that $S_{(b_0+j_0)(b_0+j_0+1)}(\phi_{0})=S_{(b_0+j_0)(b_0+j_0+1)}(\phi_{1})$, then $S_{(b_0+k)(b_0+k+1)}(\phi_{0})$ $=S_{(b_0+k)(b_0+k+1)}(\phi_{1})$ for all $k=j_0, j_0+1, \cdots, m-1$.

 (2) If the aircraft $Tcf_{b_0\!+j_0}\!$ and $Tcf_{b_0\!+j_0\!+1}\!$ are both takeoff or landing aircraft, $S_{(b_0\!+k)(b_0\!+k+1)}(\phi_{0})\!=\!S_{(b_0\!+k)(b_0\!+k\!+1)}(\phi_{1})$ for all $k=j_0, j_0+1, \cdots, m-1$.

 (3) Suppose that the aircraft $Tcf_{b_0}, Tcf_{b_0+2}, \cdots, Tcf_{b_0+2m_1}$ are all takeoff (landing) aircraft, and the aircraft $Tcf_{b_0+1}, Tcf_{b_0+3}, \cdots, Tcf_{b_0+2m_1+1}$ are all landing (takeoff) aircraft, where $m_1\geq0$ is an integer such that $2m_1+1\leq m$. If $S_{b_0(b_0+1)}(\phi_{0})\neq S_{b_0(b_0+1)}(\phi_{1})$, then $S_{(b_0+k)(b_0+k+1)}(\phi_{0})=S_{(b_0+k)(b_0+k+1)}(\phi_{1})$ for all $k=4, 5, \cdots, m-1$.}\end{theorem}

\noindent{Proof:} (1) From Assumption \ref{ass3.4}, the statement (1) holds.

 (2) If the aircraft $Tcf_{b_0}$ and $Tcf_{b_0+1}$ are both takeoff or landing aircraft, from Theorem \ref{theorem4.11b}, in $\phi_0$ and $\phi_1$, the aircraft $Tcf_{b_0+1}$ is relevant to the aircraft $Tcf_{b_0}$, and $S_{b_0(b_0+1)}(\phi_{0})=S_{b_0(b_0+1)}(\phi_{1})$. From the statement (1), the statement (2) holds.

 (3) In view of the condition that $S_{b_0(b_0+1)}(\phi_{0})\neq S_{b_0(b_0+1)}(\phi_{1})$, without loss of generality, suppose that $S_{b_0(b_0+1)}(\phi_{0})>S_{b_0(b_0+1)}(\phi_{1})$. It follows that $S_{b_0(b_0+1)}(\phi_{0})>Y_{b_0(b_0+1)}$. That is, the aircraft $Tcf_{b_0+1}$ is not relevant to the aircraft $Tcf_{b_0}$ in $\phi_0$. From Assumption \ref{ass3.4}, the aircraft $Tcf_{b_0+1}$ is relevant to the last aircraft of $\phi_{01}$. If $S_{(b_0+1)(b_0+2)}(\phi_{0})=S_{(b_0+1)(b_0+2)}(\phi_{1})$, from the statement (1), $S_{(b_0+k)(b_0+k+1)}(\phi_{0})=S_{(b_0+k)(b_0+k+1)}(\phi_{1})$ for all $k=1, 2, \cdots, m-1$. Suppose that $S_{(b_0+1)(b_0+2)}(\phi_{0})\neq S_{(b_0+1)(b_0+2)}(\phi_{1})$, from Lemma \ref{lemma4.3}, $S_{(b_0+1)(b_0+2)}(\phi_{0})<S_{(b_0+1)(b_0+2)}(\phi_{1})$, and the aircraft $Tcf_{b_0+2}$ is relevant to $Tcf_{b_0}$ in $\phi_{1}$. By analogy, if $S_{(b_0+2)(b_0+3)}(\phi_{0})=S_{(b_0+2)(b_0+3)}(\phi_{1})$, then $S_{(b_0+k)(b_0+k+1)}(\phi_{0})=S_{(b_0+k)(b_0+k+1)}(\phi_{1})$ for all $k=2, 3, \cdots, m-1$. When $S_{(b_0+2)(b_0+3)}(\phi_{0})\neq S_{(b_0+2)(b_0+3)}(\phi_{1})$, from Lemma \ref{lemma4.3}, $S_{(b_0+2)(b_0+3)}(\phi_{0})>S_{(b_0+2)(b_0+3)}(\phi_{1})$, and the aircraft $Tcf_{b_0+3}$ is relevant to $Tcf_{b_0+1}$ in $\phi_{0}$.

 Note that the aircraft $Tcf_{b_0+1}$ is relevant to the last aircraft of $\phi_{01}$ in $\phi_0$ and the aircraft $Tcf_{b_0+3}$ is relevant to the aircraft $Tcf_{b_0+1}$ in $\phi_{0}$. From Assumption \ref{ass4.21a}, the aircraft $Tcf_{b_0+5}$ is relevant to the aircraft $Tcf_{b_0+4}$ in $\phi_0$ and $Y_{(b_0+3)(b_0+5)}<T_D+D_T$. Since $Tcf_{b_0+4}$ is a takeoff (landing) aircraft and $Tcf_{b_0+3}, Tcf_{b_0+5}$ are landing (takeoff) aircraft, $S_{(b_0+3)(b_0+5)}(\phi_1)\geq T_D+D_T>Y_{(b_0+3)(b_0+5)}$. Therefore, the aircraft $Tcf_{b_0+5}$ is also relevant to the aircraft $Tcf_{b_0+4}$ in $\phi_1$. From statement (1), $S_{(b_0+k)(b_0+k+1)}(\phi_{0})=S_{(b_0+k)(b_0+k+1)}(\phi_{1})$ for all $k=4, 5, \cdots, m-1$.

\begin{remark}{\rm {Theorem \ref{theorem4.5} shows that the adjustments of the aircraft sequence at some point might affect at most $4$ aircraft and have no impact on other separation times between aircraft.}}\end{remark}

\section{Scheduling problem of landing and takeoff aircraft on dual runways}\label{secdual}

In this section, scheduling problem of mixed landing and takeoff aircraft on dual runways whose spacing is no larger than $760$ $m$, where  all of the landing aircraft lands on one runway and all of the takeoff aircraft take off from the other runway. Let $P_D$ denote the minimum  separation time between a takeoff aircraft and a leading landing aircraft. Let $D_P$ denote the minimum separation time between a landing aircraft and a leading takeoff aircraft.

\begin{assumption}\label{ass4.1}{\rm Suppose that $D_P=T_0$, $P_D=0$, and for any two aircraft $Tcf_i$ and $Tcf_j$ with the same operation tasks, $Y_{ij}=T_{cl_icl_j}$ or $Y_{ij}=D_{cl_icl_j}$.}\end{assumption}

\begin{remark}\label{remakrr}{\rm If one aircraft $Tcf_i$ lands and one aircraft $Tcf_j$ takes off at the same time, it is assumed by default that aircraft $Tcf_i$ is ahead of aircraft $Tcf_j$, i.e., $S_i\leq S_j$, the aircraft  $Tcf_i$ is the leading aircraft and $Tcf_j$ is the trailing aircraft, and it is said that the aircraft $Tcf_j$ is relevant to the aircraft $Tcf_i$.}\end{remark}

Consider a group of landing aircraft and takeoff aircraft operating on dual runways. Let $\Phi=\langle Tcf_1, Tcf_2, \cdots, $ $Tcf_{n+m}\rangle$ denote the whole mixed landing and takeoff sequence on dual runways, $\Phi_0=\langle Tcf^0_1, Tcf^0_2, \cdots, Tcf^0_n\rangle$ denote the landing aircraft sequence in $\Phi$, $\Phi_1=\langle Tcf^1_{1}, Tcf^1_{2}, \cdots, Tcf^1_{m}\rangle$ denote the takeoff aircraft sequence in $\Phi$. From the definition of aircraft sequence, it follows that $S_1(\Phi)\leq S_2(\Phi)\leq \cdots \leq S_{n+m}(\Phi)$, $S_1(\Phi_0)\leq S_2(\Phi_0)\leq \cdots \leq S_{n}(\Phi_0)$ and $S_1(\Phi_1)\leq S_2(\Phi_1)\leq \cdots \leq S_{m}(\Phi_1)$.

In the following, we first discuss the aircraft relevance in Lemma \ref{lemma111s} and Theorem \ref{theorem4.11}, and based on Theorem \ref{theorem4.11}, we give an optimal conditions for the optimization problem (\ref{optim1}) in Theorem \ref{theorem4.4}.

\begin{lemma}\label{lemma111s}{\rm Consider the aircraft sequence $\Phi$. Suppose that $Tf_k=[t_0, +\infty]$ for all $k$ and Assumptions \ref{ass1.3.1}-\ref{ass1.3.4}, \ref{ass2.3.1}-\ref{ass2.3.41} and \ref{ass4.1} hold.
 For a landing aircraft $Tcf_{i}$ and a takeoff aircraft $Tcf_{k}$ in $\Phi$, if aircraft $Tcf_{k}$ is relevant to aircraft $Tcf_{i}$, then $k=i+1$.
For a takeoff aircraft $Tcf_{i}$ and a landing aircraft $Tcf_{k}$ in $\Phi$, if aircraft $Tcf_{k}$ is relevant to aircraft $Tcf_{i}$, then $k=i+1$.}\end{lemma}

\noindent{Proof:} For a landing aircraft $Tcf_{i}$ and a takeoff aircraft $Tcf_{k}$ in $\Phi$, if aircraft $Tcf_{k}$ is relevant to aircraft $Tcf_{i}$, then $S_i(\Phi)=S_k(\Phi)$. Note that $|S_i-S_j|\geq T_0$ for each $Tcf_j\in \Phi_1-\{Tcf_i\}$ and $|S_k-S_h|\geq T_0$ for each $Tcf_h\in \Phi_0-\{Tcf_i\}$, and $S_{h}\leq S_{j}$ for all $0<h< j$ from the sequence definition. Note that if one aircraft $Tcf_a$ lands and one aircraft $Tcf_b$ takes off at the same time, it is assumed by default that aircraft $Tcf_a$ is ahead of aircraft $Tcf_b$. Thus, $k=i+1$. For a takeoff aircraft $Tcf_{i}$ and a landing aircraft $Tcf_{k}$ in $\Phi$, if aircraft $Tcf_{k}$ is relevant to aircraft $Tcf_{i}$, $S_k(\Phi)-S_i(\Phi)=T_0$.
Note that $|S_i-S_j|\geq T_0$ for each $Tcf_j\in \Phi_1-\{Tcf_i\}$ and $|S_k-S_h|\geq T_0$ for each $Tcf_h\in \Phi_0-\{Tcf_i\}$,  and $S_{h}\leq S_{j}$ for all $0<h< j$ from the sequence definition.
 Also, note that if one aircraft $Tcf_a$ lands and one aircraft $Tcf_b$ takes off at the same time, it is assumed by default that aircraft $Tcf_a$ is ahead of aircraft $Tcf_b$. Since the takeoff aircraft $Tcf_{k}$ is relevant to the landing aircraft $Tcf_{i}$, then $k=i+1$.

\begin{theorem}\label{theorem4.11}{\rm  Consider the aircraft sequence $\Phi$. Suppose that $Tf_k=[t_0, +\infty]$ for all $k$. Under Assumptions \ref{ass1.3.1}-\ref{ass1.3.4}, \ref{ass2.3.1}-\ref{ass2.3.41} and \ref{ass4.1}, the following statements hold.

(1) Suppose that aircraft $Tcf_i$ is a takeoff (landing) aircraft and aircraft $Tcf_j$ is a landing (takeoff) aircraft. If aircraft $Tcf_j$ is relevant to aircraft $Tcf_i$, then $j=i+1$.

(2) Suppose that aircraft $Tcf_i$ and $Tcf_{i+1}$ are both landing (takeoff) aircraft. If aircraft $Tcf_{i+1}$ is relevant to some aircraft $Tcf_k$, then $k=i$.

(3) Suppose that aircraft $Tcf_i$ and $Tcf_j$ are both landing (takeoff) aircraft, and $j-i>1$. If aircraft $Tcf_j$ is relevant to aircraft $Tcf_i$, aircraft $Tcf_{i+1}$ and $Tcf_{j-1}$ are takeoff (landing) aircraft.

(4) Suppose that aircraft $Tcf_j$ is relevant to aircraft $Tcf_i$. If $Tcf_i$ is a landing (takeoff) aircraft, then aircraft $Tcf_{i+1}$, $Tcf_{i+2}$, $\cdots$, $Tcf_{j-1}$ are all takeoff (landing) aircraft.

(5) Suppose that the aircraft $Tcf_j$ is relevant to aircraft $Tcf_i$. Then $0 < j-i \leq 4$.

}\end{theorem}

\noindent{Proof:} (1) From Lemma \ref{lemma111s}, this statement holds.

(2) From the sequence definition, $S_{h(i+1)}\geq S_{i(i+1)}\geq T_0$ for all $0<h< i$. If $k\neq i$, then $0<k<i$. Since $D_P=T_0$, if the landing aircraft $Tcf_{i+1}$ is relevant to a takeoff aircraft $Tcf_k$, then $Y_{k(i+1)}=T_0$. It is assumed that an aircraft lands and an aircraft takes off at the same time, it is assumed by default that the landing aircraft is ahead of the takeoff aircraft. Then the takeoff time of aircraft $Tcf_k$ is later than the landing time of aircraft $Tcf_i$, i.e., $k>i$.
This yields a contradiction. Therefore, the landing aircraft $Tcf_{i+1}$ is relevant to a landing aircraft $Tcf_k$. From Theorem \ref{theorem1.1}, $k=i$. If aircraft $Tcf_i$ and $Tcf_{i+1}$ are both takeoff aircraft, it can be similarly proved that $k=i$.

(3) If one of the aircraft $Tcf_{i+1}$ and $Tcf_{j-1}$, denoted by $Tcf_h$, is a landing (takeoff) aircraft, from Lemmas \ref{lemma11s} and \ref{lemma22s}, $S_{ij}\geq Y_{ih}+Y_{hj}>Y_{ij}$, implying that aircraft $Tcf_j$ is not relevant to aircraft $Tcf_i$, which is a contradiction.

(4) From statements (1) and (3), this statement holds naturally.

(5) Since the aircraft $Tcf_j$ is relevant to aircraft $Tcf_i$, $Y_{ij}\leq 3T_0$. If aircraft $Tcf_i$ and $Tcf_j$ are both landing aircraft, $Tcf_{i+1}, Tcf_{i+2}, \cdots, Tcf_{j-1}$ are all takeoff aircraft. If $j-i>4$, it follows that $S_{ij} \ge S_{(i+1)(j-1)} + S_{(j-1)j} \geq 3T_0 + T_0 > Y_{ij}$, implying that aircraft $Tcf_j$ is not relevant to aircraft $Tcf_i$, which yields a contradiction.
Suppose that one of aircraft $Tcf_i$ and $Tcf_j$ is a landing aircraft and one is a takeoff aircraft. If $j-i>4$, it is easy to see that $S_{ij}>T_0$, and aircraft $Tcf_j$ cannot be relevant to aircraft $Tcf_i$. If aircraft $Tcf_i$ and $Tcf_j$ are both takeoff aircraft and $j-i>4$, aircraft $Tcf_j$ cannot be relevant to aircraft $Tcf_i$.
Therefore, if the aircraft $Tcf_j$ is relevant to aircraft $Tcf_i$, $0<j-i\leq4$.

\begin{remark}{\rm Theorem \ref{theorem4.11} shows that aircraft might be relevant to its nearest landing and takeoff aircraft ahead and is not relevant to other aircraft.}\end{remark}

\begin{remark}{\rm When the scheduling problem of mixed landing and takeoff aircraft on dual runways whose spacing is no larger than $760$ $m$ is discussed, each aircraft might be relevant to two aircraft.}\end{remark}

\begin{theorem}\label{theorem4.4}{\rm Consider the sequence $\Phi$. Suppose that the orders of all aircraft in $\Phi$ is fixed and $Tf_k=[t_0, +\infty]$ for all $k$. Under Assumptions \ref{ass1.3.1}-\ref{ass1.3.4}, \ref{ass2.3.1}-\ref{ass2.3.41} and \ref{ass4.1}, the following statements hold.

(1) The optimization problem (\ref{optim1}) can be solved if and only if there is a path from the aircraft $Tcf_{n+m}$ to the aircraft $Tcf_1$.

(2) If the aircraft $Tcf_2$ is relevant to the aircraft $Tcf_1$ and for each $i\in \{3,4, \cdots, n+m\}$, and the aircraft $Tcf_{i}$ is relevant to an aircraft $Tcf_k$ for $k\in \{Tcf_{i-4}, Tcf_{i-3}, Tcf_{i-2}, Tcf_{i-1}\}$, then the optimization problem (\ref{optim1}) can be solved.
}\end{theorem}

\noindent{Proof:} The proof is similar to that of Theorem \ref{theorem3.4}.

\begin{remark}{\rm Theorem \ref{theorem4.4}(2) means that when each aircraft $Tcf_i$ is relevant to the nearest landing or takeoff aircraft, the optimization problem (\ref{optim1}) is solved. Moreover, it should be noted that when the subscript of a given aircraft $Tcf_j$ is not positive, the aircraft $Tcf_j$ is not taken into account by default. }\end{remark}

\begin{assumption}\label{ass4.21}{\rm Suppose that for each aircraft $Tcf_i\in \{2,3, \cdots, n+m\}$ in $\Phi$, the aircraft $Tcf_{i}$ is relevant to an aircraft $Tcf_k$ for $Tcf_k\in \{Tcf_{i-4}, Tcf_{i-3}, Tcf_{i-2}, Tcf_{i-1}\}$. If aircraft $Tcf_{i}$ and $Tcf_{i-1}$ have different (takeoff/landing) operation tasks and aircraft $Tcf_i$ is relevant to $Tcf_{i-1}$ for $i\in \{2,3, \cdots, n+m\}$, there is at least an aircraft $Tcf_{j}\in \{Tcf_1, Tcf_2, \cdots, Tcf_{k-1}\}$ such that $S_{ji}(\Phi)-T_0<Y_{ji}$.}\end{assumption}

To illustrate this assumption, consider the sequence $\phi_1=\langle Tcf_1, Tcf_2, Tcf_3, Tcf_4\rangle$ and $\phi_2=\langle Tcf_1, Tcf_2, Tcf_4,$ $ Tcf_3\rangle$, where $Tcf_1$ and $Tcf_4$ are landing aircraft such that $Y_{14}=T_0$, and $Tcf_2$ and $Tcf_3$ are takeoff aircraft such that $Y_{23}=T_0$. Suppose that $S_{12}(\phi_1)=S_{12}(\phi_2)=S_{43}(\phi_2)=0$ and $S_{23}(\phi_1)=S_{34}(\phi_1)=S_{24}(\phi_2)=T_0$. It can be checked that $S_{14}(\phi_1)-T_0\geq Y_{14}$,
$S_{24}(\phi_1)-T_0\geq Y_{24}$. Thus, there is no aircraft $Tcf_{j}\in \{Tcf_1, Tcf_2\}$ such that $S_{j4}(\phi_1)-T_0<Y_{j4}$. That is, Assumption \ref{ass4.21} does not hold for $\phi_1$. Moreover, it can be checked that $S_{14}(\phi_2)-T_0<Y_{14}$.
Comparing these two sequences, it can be obtained that $F(\phi_1)-F(\phi_2)>0$, showing that Assumption \ref{ass4.21} can be used to exclude some non-optimal cases.

When there is a path from the last aircraft $Tcf_{n+m}$ to the first aircraft $Tcf_1$ in $\Phi$, there might be some aircraft that do not satisfy Assumption \ref{ass4.21}. We can adjust the landing or takeoff times of these aircraft to make the aircraft sequence satisfy Assumption \ref{ass4.21}.

\begin{assumption}\label{ass411}{\rm Consider the sequence $\Phi$. Suppose that $Tf_j=[t_0, +\infty]$ for all $j$ and Assumptions \ref{ass1.3.1}-\ref{ass1.3.4}, \ref{ass2.3.1}-\ref{ass2.3.41} and \ref{ass4.1}, \ref{ass4.21} hold.}\end{assumption}

\begin{theorem}\label{theorem5.366}{\rm Consider the sequence $\Phi$. Suppose that Assumption \ref{ass411} holds, and there exists a path $\phi_0$ from the aircraft $Tcf_{n+m}$ to the aircraft $Tcf_1$ in $\Phi$. Suppose that $(\Phi_0, Sr(\Phi_0))$ is an optimal solution of the optimization problem (\ref{optim1}) for the aircraft $Tcf^0_1, Tcf^0_2, \cdots, Tcf^0_n$ and $(\Phi_1, Sr(\Phi_1))$ is an optimal solution of the optimization problem (\ref{optim1}) for the aircraft $Tcf^1_{1}, Tcf^1_{2}, \cdots, Tcf^1_{m}$.

(1) If $\phi_0$ is formed completely by the landing aircraft of $\Phi$, $(\Phi, Sr(\Phi))$ is an optimal solution of the optimization problem (\ref{optim1}).

(2) If $\phi_0$ is formed completely by the takeoff aircraft of $\Phi$, $(\Phi, Sr(\Phi))$ is an optimal solution of the optimization problem (\ref{optim1}).
}\end{theorem}

\begin{remark}{\rm Theorem \ref{theorem5.366} gives a rule for the case of two runways with spacing no larger than $760$ $m$ when there exists a landing or takeoff aircraft subsequence playing a leading role in the whole landing and takeoff aircraft sequence.
}\end{remark}

In the following, we will present important definitions of block subsequences,  which are very useful to significantly reduce the computation to find the optimal sequence for the optimization problem (\ref{optim1}). Before that, we need to first introduce a definition of semi-resident-point aircraft, which are key aircrafts of the block subsequences, and study the impact of the change of semi-resident-point aircraft.

{\begin{definition}{\rm (Semi-resident-point aircraft) Consider the sequence $\Phi$. If a landing (takeoff) aircraft $Tcf_i$ is a resident-point aircraft in $\Phi_0$ ($\Phi_1$) and relevant to a takeoff (landing) aircraft in $\Phi$, the aircraft $Tcf_i$ is said to be a semi-resident-point aircraft in $\Phi$.}\end{definition}}

\begin{theorem}\label{thep1}{\rm Consider the aircraft sequence $\Phi$. Let $\Phi_a$ be a sequence generated by exchanging the orders of aircraft $Tcf_k$ and $Tcf_{k+1}$ for some $2<k<n+m-1$. Suppose that Assumption \ref{ass411} holds for the sequences $\Phi$ and $\Phi_a$, and $Tcf_{k+1}$ is relevant to only $Tcf_{k}$ in $\Phi$.

(1) Suppose that aircraft $Tcf_{k}$ is a landing aircraft and $Tcf_{k+1}$ is a takeoff aircraft, and $S_k(\Phi_a)>S_{k+1}(\Phi_a)>S_{k-1}(\Phi_a)$. Then $S_{k}(\Phi)=S_{k+1}(\Phi)$, $S_{k}(\Phi_a)-S_{k+1}(\Phi_a)=T_0$, and
$S_{k+1}(\Phi)-S_{k+1}(\Phi_a)+S_k(\Phi_a)-S_k(\Phi)=T_0$.

(2)  Suppose that aircraft $Tcf_{k}$ is a takeoff aircraft and $Tcf_{k+1}$ is a landing aircraft, and $S_{k+1}(\Phi)>S_{k}(\Phi)>S_{k-1}(\Phi)$. Then $S_{k}(\Phi_a)=S_{k+1}(\Phi_a)$, $S_{k+1}(\Phi)-S_{k}(\Phi)=T_0$, and
$S_{k+1}(\Phi)-S_{k+1}(\Phi_a)+S_k(\Phi_a)-S_k(\Phi)=T_0$. }\end{theorem}

\noindent{Proof:} (1)
In $\Phi$, $Tcf_{k+1}$ is located after $Tcf_k$. Thus, in $\Phi_a$, $Tcf_{k}$ is located after $Tcf_{k+1}$. From Remark \ref{remakrr} and Assumptions \ref{ass4.1}, \ref{ass4.21},  $S_k(\Phi_a)=S_{k+1}(\Phi_a)+T_0$. Since $Tcf_{k+1}$ is relevant to only $Tcf_{k}$ in $\Phi$, $S_{k+1}(\Phi)=S_{k}(\Phi)$. It follows that $S_{k+1}(\Phi)-S_{k+1}(\Phi_a)+S_k(\Phi_a)-S_k(\Phi)=T_0$.

(2) Similar to the proof of statement (1), this statement can be proved and the proof is omitted.

\begin{remark}{\rm In Theorem \ref{thep1}, the aircraft $Tcf_{k+1}$ is actually a semi-resident-point aircraft. From Theorem \ref{thep1}, the relative time increment is $T_0$ when the orders of aircraft $Tcf_k$ and $Tcf_{k+1}$ are exchanged, which is a significant property for the algorithm proposed later.}\end{remark}

\begin{theorem}\label{thep2}{\rm Consider the aircraft sequence $\Phi$. Let $\Phi_a$ be a sequence generated by exchanging the orders of aircraft $Tcf_k$ and $Tcf_{k+1}$ for some $2<k<n+m-1$. Suppose that Assumption \ref{ass411} holds for the sequences $\Phi$ and $\Phi_a$, and $Tcf_{k+1}$ is only relevant to $Tcf_{k}$ in $\Phi$. 

 (1) Suppose that $S_k(\Phi_a)>S_{k+1}(\Phi_a)>S_{k-1}(\Phi_a)$, aircraft $Tcf_k$ is $Tcf_{i_0}^0$ in $\Phi_0$ and aircraft $Tcf_{k+1}$ is $Tcf_{i_1}^1$ in $\Phi_1$, and $S_{i_0+h+1}(\Phi_0)-S_{i_0+h}(\Phi_0)=S^a_{i_0+h+1}(\Phi_0)-S^a_{i_0+h}(\Phi_0)$ for all $h\in \{0, 1, \cdots, n-i_0-1\}$, where $S_{h}$ denotes the aircraft landing/takeoff time in $\Phi$ and $S^a_{h}$ denotes the aircraft landing/takeoff time in $\Phi_a$.

(1.1) If $S_{i_1+h}(\Phi_1)-S_{i_1+h-1}(\Phi_1)=S_{i_0+h}(\Phi_0)-S_{i_0+h-1}(\Phi_0)=T_0$ and $S_{i_1+h_0}(\Phi_1)-S_{i_1+h_0-1}(\Phi_1)=S_{i_0+h_0}(\Phi_0)-S_{i_0+h_0-1}(\Phi_0)\triangleq Q_1>T_0$ for all $h\in\{1,2,\cdots,h_0-1\}$ and some integer $0<h_0\leq\min\{n-i_0,m-i_1\}$, $S^a_{i_1+h}(\Phi_1)=S^a_{i_1+h-1}(\Phi_1)+T_0$ for all $h\in \{1, 2, \cdots, h_0-1\}$ and $S^a_{i_1+h_0}(\Phi_1)=S^a_{i_0+h_0-1}(\Phi_0)+Q_1-T_0$.

(1.2)  If $S_{i_1+h}(\Phi_1)-S_{i_1+h-1}(\Phi_1)=S_{i_0+h}(\Phi_0)-S_{i_0+h-1}(\Phi_0)=T_0$ and $2T_0\geq S_{i_1+h_0}(\Phi_1)-S_{i_1+h_0-1}(\Phi_1)>S_{i_0+h_0}(\Phi_0)-S_{i_0+h_0-1}(\Phi_0)\geq T_0$ for all $h\in\{1,2,\cdots,h_0-1\}$ and some integer $0<h_0<\min\{n-i_0,m-i_1\}$, $S^a_{i_1+h}(\Phi_1)=S^a_{i_1+h-1}(\Phi_1)+T_0$ for all $h\in \{1, 2, \cdots, h_0-1\}$ and $S^a_{i_1+h_0}(\Phi_1)=S^a_{i_0+h_0}(\Phi_0)$.

(2)  Suppose that $S_{k+1}(\Phi)>S_{k}(\Phi)>S_{k-1}(\Phi)$, aircraft $Tcf_k$ is $Tcf_{i_1}^1$ in $\Phi_1$ and aircraft $Tcf_{k+1}$ is $Tcf_{i_0}^0$ in $\Phi_0$, and $S_{i_1+h+1}(\Phi_1)-S_{i_1+h}(\Phi_1)=S^a_{i_1+h+1}(\Phi_1)-S^a_{i_1+h}(\Phi_1)$ for all $h\in \{0, 1, \cdots, m-i_1-1\}$.

(2.1) If $S_{i_0+h}(\Phi_0)-S_{i_0+h-1}(\Phi_0)=S_{i_1+1}(\Phi_1)-S_{i_1}(\Phi_1)=S_{i_1+h+1}(\Phi_1)-S_{i_1+h}(\Phi_1)=T_0$ and $2T_0\geq S_{i_0+h_0}(\Phi_0)-S_{i_0+h_0-1}(\Phi_0)=S_{i_1+h_0+1}(\Phi_1)-S_{i_1+h_0}(\Phi_1)> T_0$ for all $h\in\{1,2,\cdots,h_0-1\}$ and some integer $0<h_0<\min\{n-i_0,m-i_1\}$, $S^a_{i_0+h}(\Phi_0)=S^a_{i_0+h-1}(\Phi_0)+T_0$ for all $h\in \{1, 2, \cdots, h_0-1\}$ and $S^a_{i_0+h_0}(\Phi_0)=S^a_{i_1+h_0}(\Phi_1)+T_0$.

(2.2)  If $S_{i_0+h}(\Phi_0)-S_{i_0+h-1}(\Phi_0)=S_{i_1+1}(\Phi_1)-S_{i_1}(\Phi_1)=S_{i_1+h+1}(\Phi_1)-S_{i_1+h}(\Phi_1)=T_0$ and $T_0\leq S_{i_0+h_0}(\Phi_0)-S_{i_0+h_0-1}(\Phi_0)-[S_{i_1+h_0+1}(\Phi_0)-S_{i_1+h_0}(\Phi_0)]\leq 2T_0$ for all $h\in\{1,2,\cdots,h_0-1\}$ and some integer $0<h_0<\min\{n-i_0,m-i_1\}$, $S^a_{i_0+h}(\Phi_0)=S^a_{i_0+h-1}(\Phi_0)+T_0$ for all $h\in \{1, 2, \cdots, h_0-1\}$ and $S^a_{i_0+h_0}(\Phi_0)=S^a_{i_1+h_0+1}(\Phi_1)+T_0$.

}\end{theorem}

\noindent{Proof:} From Theorem \ref{thep1}(1), $S_{i_1}(\Phi_1)=S_{i_0}(\Phi_0)$ and $S^a_{i_0}(\Phi_0)-S^a_{i_1}(\Phi_1)=T_0$, while from Theorem \ref{thep1}(2),
$S_{i_1}(\Phi_1)+T_0=S_{i_0}(\Phi_0)$ and $S^a_{i_0}(\Phi_0)=S^a_{i_1}(\Phi_1)$. By simple calculations, this theorem can be proved.

In Theorems \ref{thep1} and \ref{thep2}, different orders of the aircraft $Tcf_k$ and $Tcf_{k+1}$ discussed might have different effects on the aircraft after them. Theorem \ref{thep1} shows the relative changes of the aircraft after them and Theorem \ref{thep2} shows the weakening process of the effects of the order changes of aircraft $Tcf_k$ and $Tcf_{k+1}$. These two theorems are also very useful in our algorithm proposed later.

\begin{theorem}\label{theorem5.5}{\rm  Consider the sequence $\Phi$. Suppose that Assumption \ref{ass411} holds.

(1) Suppose that $S_i(\Phi_0)=S_{j}(\Phi_1)$ and $Y_{j(j+1)}(\Phi_1)=Y_{(j+1)(j+2)}(\Phi_1)=Y_{(j+2)(j+3)}(\Phi_1)=T_0$. There exists an integer $k_0\in \{j+1, j+2, j+3\}$ such that
$S_{i+1}(\Phi_0)=S_{k_0}(\Phi_1)$.

(2) Suppose that $S_i(\Phi_0)=S_{j}(\Phi_1)$ and $Y_{i(i+1)}(\Phi_0)=Y_{(i+1)(i+2)}(\Phi_0)=Y_{(i+2)(i+3)}(\Phi_0)=Y_{(j+1)(j+2)}(\Phi_1)=T_0$. There exist two integers $k_1\in \{i+1, i+2, i+3\}$ and $k_2\in\{j+1,j+2\}$ such that $S_{k_1}(\Phi_0)=S_{k_2}(\Phi_1)$.

(3) Suppose that $S_i(\Phi_0)<S_{j}(\Phi_1)< S_{i+1}(\Phi_0)$ and $Y_{j(j+1)}(\Phi_1)=Y_{(j+1)(j+2)}(\Phi_1)=Y_{(j+2)(j+3)}(\Phi_1)=T_0$. There exists an integer $k_0\in \{j+1, j+2, j+3\}$ such that
$S_{i+1}(\Phi_0)=S_{k_0}(\Phi_1)$.

(4) Suppose that $S_{j}(\Phi_1)<S_i(\Phi_0)<S_{j+1}(\Phi_1)$ and $Y_{j(j+1)}(\Phi_1)=Y_{(j+1)(j+2)}(\Phi_1)=Y_{(j+2)(j+3)}(\Phi_1)=Y_{(j+3)(j+4)}(\Phi_1)=T_0$. There exists an integer $k_0\in \{j+1, j+2, j+3, j+4\}$ such that
$S_{i+1}(\Phi_0)=S_{k_0}(\Phi_1)$.

(5) Suppose that $S_i(\Phi_0)< S_{j}(\Phi_1)<S_{i+1}(\Phi_0)$ and $Y_{i(i+1)}(\Phi_0)=Y_{(i+1)(i+2)}(\Phi_0)=Y_{(i+2)(i+3)}(\Phi_0)=Y_{(i+3)(i+4)}(\Phi_0)=Y_{(j+1)(j+2)}(\Phi_1)=T_0$. There exist two integers $k_1\in \{i+1, i+2, i+3, i+4\}$ and $k_2\in\{j+1,j+2\}$ such that $S_{k_1}(\Phi_0)=S_{k_2}(\Phi_1)$.

(6) Suppose that $S_{j}(\Phi_1)<S_i(\Phi_0)<S_{j+1}(\Phi_1)$ and $Y_{i(i+1)}(\Phi_0)=Y_{(i+1)(i+2)}(\Phi_0)=Y_{(i+2)(i+3)}(\Phi_0)=Y_{(j+1)(j+2)}(\Phi_1)=T_0$. There exist two integers $k_1\in \{i+1, i+2, i+3\}$ and $k_2\in\{j+1,j+2\}$ such that $S_{k_1}(\Phi_0)=S_{k_2}(\Phi_1)$.
}\end{theorem}

\noindent{Proof:} (1) To prove statement (1), we perform the analysis in two cases: (1.1) the aircraft $Tcf^0_{i+1}$ is not relevant to the aircraft $Tcf^0_{i}$ and (1.2) the aircraft $Tcf^0_{i+1}$ is relevant to the aircraft $Tcf^0_{i}$.

Case (1.1). Since the aircraft $Tcf^0_{i+1}$ is not relevant to the aircraft $Tcf^0_{i}$, $S_{j(j+3)}(\Phi_1)\geq 3T_0$ and $D_P=T_0$, the aircraft $Tcf^0_{i+1}$ is relevant to the aircraft $Tcf^1_{k}$, $k\geq j+1$. If $k\geq j+3$, it can be checked that $S_{i(i+1)}(\Phi_0)-T_0\geq Y_{i(i+1)}$, i.e., Assumption \ref{ass4.21} does not hold.
Thus, $k<j+3$. Note that $Y_{j(j+1)}(\Phi_1)=Y_{(j+1)(j+2)}(\Phi_1)=Y_{(j+2)(j+3)}(\Phi_1)=T_0$. To ensure Assumption \ref{ass4.21} to be valid, the aircraft should satisfy that $S_{(j+h)(j+h+1)}(\Phi_1)=T_0$ for $h=0, 1, \cdots, k$ and $S_{i+1}(\Phi_0)-S_{k}(\Phi_1)=T_0$. Clearly, $S_{i+1}(\Phi_0)=S_{k+1}(\Phi_1)$ where $k+1\in \{j+1, j+2, j+3\}$.

Case (1.2). Since $D_P=T_0$, for any integer $k\geq j$ if $S_{k}(\Phi_1)\leq S_{i+1}(\Phi_0)$, then $S_{i+1}(\Phi_0)-S_{k}(\Phi_1)\geq T_0$. Since the aircraft $Tcf^0_{i+1}$ is relevant to the aircraft $Tcf^0_{i}$, $T_0\leq S_{i(i+1)}(\Phi_0)\leq 3T_0$. Thus, if $S_{k}(\Phi_1)\leq S_{i+1}(\Phi_0)$ for some integer $k\geq j$, under Assumption \ref{ass4.21}, it follows that $j\leq k\leq j+2$. Note that $S_{i+1}(\Phi_0)-S_{j+k}(\Phi_1)\geq T_0$ and $Y_{j(j+1)}(\Phi_1)=Y_{(j+1)(j+2)}(\Phi_1)=Y_{(j+2)(j+3)}(\Phi_1)=T_0$. Thus, under Assumption \ref{ass4.21}, $S_{i+1}(\Phi_0)=S_{k+1}(\Phi_1)$ where $j\leq k\leq j+2$.

Summarizing the above analysis of both cases, there exists an integer $k_0\in \{j+1, j+2, j+3\}$ such that
$S_{i+1}(\Phi_0)=S_{k_0}(\Phi_1)$.

(2)-(6) The proofs are similar to that for the statement (1) and hence omitted.

\begin{remark}{\rm From Theorem \ref{theorem5.5}, it can be observed that the unequal separation times caused possibly by the existence of breakpoint aircraft and resident-point aircraft  might result in mismatches between the landing sequence $\Phi_0$ and
the takeoff sequence $\Phi_1$, but the effects of the unequal separation times would disappear after $1$ to $4$ aircraft in the absence of the influence of other
breakpoint aircraft and resident-point aircraft. Moreover, it should be noted here that the forms of the minimum separation times such as $Y_{j(j+1)}(\Phi_1)$ and $Y_{(j+1)(j+2)}(\Phi_1)$ are used instead of the forms $Y_{j(j+1)}$ and $Y_{(j+1)(j+2)}$ to distinguish different separation times between different sequences. }\end{remark}

Now, we present an important definition of block subsequences, which is very useful to significantly reduce the computation to find the optimal sequence for the optimization problem (\ref{optim1}).

\begin{definition}{\rm (T-block, D-block, TD-block subsequences) Consider a subsequence of the sequence $\Phi$, denoted by $\phi=\langle Tcf_{i_0}, Tcf_{i_0+1}, \cdots, Tcf_{i_1}\rangle$ for two integers $0<i_0<i_1\leq n+m$. Let $\phi_0$ denote the landing aircraft subsequence of $\phi$, $\phi_1$ denote the takeoff aircraft subsequence of $\phi$, and $|\phi_0|$ and $|\phi_1|$ denote the aircraft numbers in $\phi_0$ and $\phi_1$. Suppose that there is a path from aircraft $Tcf_{i_1}$ to aircraft $Tcf_{i_0}$,  the operation (landing/takeoff) tasks of $Tcf_{i_0}, Tcf_{i_0+1}$ are different, and $Tcf_{i_0+1}$ is relevant to $Tcf_{i_0}$.

(1) Suppose that $Tcf_{i_1-1}\in \phi_0$, $Tcf_{i_1}\in \phi_1$, for the sequence $\phi_0$, there is a path from its last aircraft $Tcf_{i_1-1}$ to its first aircraft, and for the sequence $\phi_1$, there is a path from its second last aircraft to its first aircraft.
If $Tcf_{i_1}$ is relevant to only $Tcf_{i_1-1}$, it is said that the sequence $\phi$ is a T-block subsequence of $\Phi$, and $(|\phi_1|-1)/(|\phi_0|-1)$, $S_{i_1}(\phi)-S_{i_0}(\phi)$ and  $S_{i_1}(\phi)-S_{i_2}(\phi)-Y_{i_2i_1}$ are the capacity, the time length and the takeoff time increment of T-block subsequence, where $Tcf_{i_2}$ is the second last takeoff aircraft in $\phi$.
Let $\langle Tcf_{k_0}, Tcf_{k_1}, \cdots, Tcf_{k_s}\rangle$ denote the takeoff sequence in the T-block subsequence $\phi$. Consider a subsequence $\bar{\phi}$ such that the landing times of the landing aircraft in $\bar{\phi}$ are consistent with those in $\phi$ and the takeoff times of the takeoff aircraft satisfy that
$S_{k_0}(\bar{\phi})-S_{k_0}(\phi)=S_{k_1}(\bar{\phi})-S_{k_1}(\phi)=\cdots=S_{k_{s-1}}(\bar{\phi})-S_{k_{s-1}}(\phi)=\mu$ and $S_{k_{s}}(\bar{\phi})=S_{k_{s}}(\phi)$ for some constant $\mu\geq 0$. Suppose that for a given constant $\mu_0$, if $0\leq \mu\leq \mu_0$, the minimum separation requirements can be satisfied in $\bar{\phi}$ and if $\mu>\mu_0$, the minimum separation requirements cannot be satisfied in $\bar{\phi}$. Then the constant $\mu_0$ is said to be the initial redundant time of the T-block subsequence $\phi$.

(2) Suppose that $Tcf_{i_1-1}\in \phi_1$, $Tcf_{i_1}\in \phi_0$, for the sequence $\phi_1$, there is a path from its last aircraft $Tcf_{i_1-1}$ to its first aircraft, and for the sequence $\phi_0$, there is a path from the second last aircraft to the first aircraft in $\phi_0$.
If $Tcf_{i_1}$ is relevant to only $Tcf_{i_1-1}$, it is said that the sequence $\phi$ is a D-block subsequence of $\Phi$, and $(|\phi_1|-1)/(|\phi_0|-1)$, $S_{i_1}(\phi)-S_{i_0}(\phi)$ and  $S_{i_1}(\phi)-S_{i_2}(\phi)-Y_{i_2i_1}$ are the capacity, the time length and the landing time increment of D-block subsequence, where $Tcf_{i_2}$ is the second last landing aircraft in $\phi$.
Let $\langle Tcf_{k_0}, Tcf_{k_1}, \cdots, Tcf_{k_s}\rangle$ denote the landing sequence in the D-block subsequence $\phi$.
Consider a subsequence $\bar{\phi}$ such that the takeoff times of the takeoff aircraft in $\bar{\phi}$ are consistent with those in $\phi$ and
$S_{k_0}(\bar{\phi})-S_{k_0}(\phi)=S_{k_1}(\bar{\phi})-S_{k_1}(\phi)=\cdots=S_{k_{s-1}}(\bar{\phi})-S_{k_{s-1}}(\phi)=\mu$ and $S_{k_{s}}(\bar{\phi})=S_{k_{s}}(\phi)$ for some constant $\mu\geq 0$. Suppose that for a given constant $\mu_0$, if $0\leq \mu\leq \mu_0$, the minimum separation requirements can be satisfied in $\bar{\phi}$ and if $\mu>\mu_0$, the minimum separation requirements cannot be satisfied in $\bar{\phi}$. Then the constant $\mu_0$ is said to be the initial redundant time of the D-block subsequence $\phi$.

(3) Suppose that the operation (landing/takeoff) tasks of $Tcf_{i_1}, Tcf_{i_1-1}$ are different, for the sequence $\phi_0$, there is a path from its last aircraft $Tcf_{i_1-1}$ to its first aircraft, and for the sequence $\phi_1$, there is a path from its last aircraft to its first aircraft. If $Tcf_{i_1}$ is relevant to $Tcf_{i_1-1}$, it is said that the sequence $\phi$ is a TD-block subsequence of $\Phi$.
}\end{definition}

Note that in the definitions of T-block subsequences and D-subsequences, the last aircraft $Tcf_{i_1}$ is a semi-resident-point aircraft.

 In fact, under Assumption \ref{ass4.21}, the sequence $\Phi$ can be expressed as a sequence of block subsequences in the form of, e.g., $\langle TB_1, TB_2, DB_1, TB_3, DB_2, DB_3, TB_4, \cdots\rangle$, where each $TB_i$ denotes a T-block subsequence and each $DB_i$ denotes a D-block subsequence. It should be noted that for any two consecutive block subsequences, the last two aircraft of the first block subsequence are the first two aircraft of the second block subsequences. Based on the definitions of block subsequences, we can focus on the switching between T-block subsequences and D-block subsequences to study the time increments of the landing sequence and the takeoff sequence of the whole sequence $\Phi$.

According to the combination of the first two aircraft, each T-block/D-block subsequence has two different forms. One form is that the first aircraft is a landing aircraft, the second aircraft is a takeoff aircraft and the separation time between these two aircraft is $0$. The other form is that the first aircraft is a takeoff aircraft, the second aircraft is a landing aircraft and the separation time between these two aircraft is $T_0$.

To illustrate the definitions of block subsequences, we give the following proposition.

\begin{proposition}\label{pprr1}{\rm Consider the sequence $\Phi$. Suppose that Assumption \ref{ass411} holds.

(1) Suppose that $S_i(\Phi_0)=S_j(\Phi_1)$, $S_{i+1}(\Phi_0)-S_i(\Phi_0)=S_{i+2}(\Phi_0)-S_{i+1}(\Phi_0)=T_0$ and $Y_{j(j+1)}(\Phi_1)=4T_0/3$ for two integers $i,j$. If $S_j(\Phi_1)<S_{j+1}(\Phi_1)\leq S_{i+2}(\Phi_0)$, it follows that $S_{j+1}(\Phi_1)=S_{i+2}(\Phi_0)$ and aircraft $Tcf^0_{i}, Tcf^1_{j}, Tcf^0_{i+1}, Tcf^0_{i+2}, Tcf^1_{j+1}$ form a T-block subsequence in $\Phi$.

(2) Suppose that $S_i(\Phi_0)=S_j(\Phi_1)$, $S_{j+1}(\Phi_j)-S_j(\Phi_1)=S_{j+2}(\Phi_1)-S_{j+1}(\Phi_1)=T_0$ and $Y_{i(i+1)}(\Phi_0)=1.5T_0$ for two integers $i,j$. If $S_j(\Phi_1)<S_{i+1}(\Phi_0)\leq S_{j+2}(\Phi_1)$, it follows that $S_{i+1}(\Phi_0)=S_{j+2}(\Phi_1)$ and aircraft $Tcf^0_{i}, Tcf^1_{j}, Tcf^0_{j+1}, Tcf^0_{i+1}, Tcf^1_{j+2}$ form a D-block subsequence in $\Phi$.
}\end{proposition}

\begin{table}[h]
    \centering
       \caption{Takeoff separation times. (T-Blocks) }
    \begin{tabular}{|c|c|c|c|c|c|c|c|}
        \hline
\diagbox{LSTV}{MTST}     & 60 & 80  & 100 & 120 & 140 & 160 & 180\\
        \hline
        (60,60,60) & 60 & 120 & 120 & 120 & 180 & 180 & 180\\
        \hline
        (60,60,90) & 60 & 120 & 120 & 120 & \textcolor{red}{140} & 210 & 210 \\
        \hline
        (60,90,90) & 60 & \textcolor{red}{80}  & 150 & 150 & 150 & \textcolor{red}{160} & 180 \\
        \hline
        (90,90,90) & 90 & 90  & \textcolor{red}{100} & \textcolor{red}{120} & 180 & 180 & 180 \\
        \hline
    \end{tabular}
    \label{tab1}
\end{table}

\begin{table}[h]
    \centering
      \caption{Takeoff separation times. (T-Blocks)}
    \begin{tabular}{|c|c|c|c|c|c|c|c|}
        \hline
\diagbox{LSTV}{MTST}     & 60 & 80  & 100 & 120 & 140 & 160 & 180\\
        \hline
        (60,60,60) & 60 & 120 & 120 & 120 & 180 & 180 & 180\\
        \hline
        (60,60,90) & 60 & 120 & 120 & 120 & 180 & 180 & 180 \\
        \hline
        (60,90,90) & 60 & 120 & 120 & 120 & \textcolor{red}{140} & 210 & 210 \\
        \hline
        (90,90,90) & 60 & \textcolor{red}{80}  & 150 & 150 & 150 & \textcolor{red}{160} & 180 \\
        \hline
    \end{tabular}
    \label{tab2}
\end{table}

\begin{table}[h]
    \centering
    \caption{Landing separation times. (D-Blocks)}
    \begin{tabular}{|c|c|c|c|c|c|c|c|}
        \hline
\diagbox{TSTV}{MLST}     & 60 & 68  & 90 & 113 & 135 & 158 & 180\\
        \hline
        (60,60,60) & 60 & 120 & 120 &120  & 180 & 180 & 180\\
        \hline
        (60,60,80) & 60 & 120 & 120 &120  & 180 & 180 & 180 \\
        \hline
        (60,80,80) & 60 & 120 & 120 &120  & \textcolor{red}{135} & 200 & 200 \\
        \hline
        (80,80,80) & \textcolor{red}{60} & \textcolor{red}{68}  & 140 &140  & 140 & \textcolor{red}{158} & 220 \\
        \hline
    \end{tabular}

    \label{tab3}
\end{table}

\begin{table}[h]
    \centering
    \caption{Landing separation times. (D-Blocks)}
    \begin{tabular}{|c|c|c|c|c|c|c|c|}
        \hline
        {\diagbox{TSTV}{MLST}}        & 60 & 68  & 90 & 113 & 135 & 158 & 180\\

        \hline
        (60,60,60) & 60 & 120 & 120 &120  & 180 & 180 & 180\\
        \hline
        (60,60,80) & 60 & 120 & 120 &120  & \textcolor{red}{135} & 200 & 200 \\
        \hline
        (60,80,80) & 60 & \textcolor{red}{68} & 140 &140  & 140 & \textcolor{red}{158} & 220 \\
        \hline
        (80,80,80) & 80 & 80 & \textcolor{red}{90} &160  & 160 & 160 & 180 \\
        \hline
    \end{tabular}

    \label{tab4}
\end{table}

In Proposition \ref{pprr1}, we only discuss the simplest block subsequences in $\Phi$. In the following tables, we show the actual separation times in block subsequences
when the separation times between the landing or takeoff aircraft are fixed, where the specific data are given based on the minimum separation time standards at Heathrow Airport and the RECAT-EU system.

In Tables \ref{tab1}, \ref{tab2}, \ref{tab5} and \ref{tab6}, the separation times between aircraft in the takeoff subsequence $\phi_1=\langle Tcf_1^1, Tcf_2^1\rangle$ and the takeoff time increments of T-block subsequences are shown, where
the rows represent the landing separation time vector between aircraft (LSTV in short) in the landing subsequence $\phi_0=\langle Tcf^0_1, Tcf^0_2, Tcf^0_3, Tcf^0_4\rangle$, and the columns represent the minimum takeoff separation time (MTST in short) between aircraft in $\phi_1$.

In Tables \ref{tab3}, \ref{tab4}, \ref{tab7} and \ref{tab8}, the separation times between aircraft in the landing subsequence $\phi_0=\langle Tcf_1^0, Tcf_2^0\rangle$ and the landing time increments of D-block subsequences are shown, where
the rows represent the takeoff separation time vector (TSTV in short) between aircraft in the takeoff subsequence $\phi_1=\langle Tcf^1_1, Tcf^1_2, Tcf^1_3, Tcf^1_4\rangle$, and the columns represent the minimum landing separation time (MLST) between aircraft in $\phi_0$.

Moreover, it is assumed that $S_1(\phi_0)=S_1(\phi_1)$ in Tables \ref{tab1}, \ref{tab3}, \ref{tab5} and \ref{tab7}, and $S_1(\phi_0)=S_1(\phi_1)+60$ in Tables \ref{tab2}, \ref{tab4}, \ref{tab6} and \ref{tab8}.

\begin{table}[!h]
    \centering
       \caption{Takeoff time increments. (T-Blocks)}\vspace{0.1cm}
    \begin{tabular}{|c|c|c|c|c|c|c|c|}
        \hline
\diagbox{LSTV}{MTST}     & 60 & 80  & 100 & 120 & 140 & 160 & 180\\
        \hline
        (60,60,60) & 0 & 40 & 20 & 0 & 40 & 20 & 0\\
        \hline
        (60,60,90) & 0 & 40 & 20 & 0 & {\color{red}0} & 50 & 30 \\
        \hline
        (60,90,90) & 0 & {\color{red}0}  & 50 & 30 & 10& {\color{red}0} & 0 \\
        \hline
        (90,90,90) & 30 & 10  & {\color{red}0} & {\color{red}0} & 40 & 20 & 0 \\
        \hline
    \end{tabular}
    \label{tab5}
\end{table}

\begin{table}[!h]
    \centering
      \caption{Takeoff time increments. (T-Blocks)}\vspace{0.1cm}
    \begin{tabular}{|c|c|c|c|c|c|c|c|}
        \hline
\diagbox{LSTV}{MTST}     & 60 & 80  & 100 & 120 & 140 & 160 & 180\\
        \hline
        (60,60,60) & 0 & 40 & 20 & 0 & 40 & 20 & 0\\
        \hline
        (60,60,90) & 0 & 40 & 20 & 0 & 40 & 20 & 0 \\
        \hline
        (60,90,90) & 0 & 40 & 20 & 0 & {\color{red}0} & 50 & 30 \\
        \hline
        (90,90,90) & 0 & {\color{red}0}  & 50 & 30 & 10 & {\color{red}0} & 0 \\
        \hline
    \end{tabular}
    \label{tab6}
\end{table}

\begin{table}[!h]
    \centering
    \caption{Landing time increments. (D-Blocks)}\vspace{0.1cm}
    \begin{tabular}{|c|c|c|c|c|c|c|c|}
        \hline
\diagbox{TSTV}{MLST}     & 60 & 68  & 90 & 113 & 135 & 158 & 180\\
        \hline
        (60,60,60) & 0 & 52 & 30 &7  & 45 & 22 & 0\\
        \hline
        (60,60,80) & 0 & 52 & 30 &7  & 45 & 22 & 0 \\
        \hline
        (60,80,80) & 0 & 52 & 30 &7  & {\color{red}0} & 42 & 20 \\
        \hline
        (80,80,80) & {\color{red}0} & {\color{red}0}  & 50 &27  & 5 & {\color{red}0} & 40 \\
        \hline
    \end{tabular}

    \label{tab7}
\end{table}

\begin{table}[!h]
    \centering
    \caption{Landing time increments. (D-Blocks)}\vspace{0.1cm}
    \begin{tabular}{|c|c|c|c|c|c|c|c|}
        \hline
        {\diagbox{TSTV}{MLST}}        & 60 & 68  & 90 & 113 & 135 & 158 & 180\\

        \hline
        (60,60,60) & 0 & 52 & 30 &7  & 45 & 22 & 0\\
        \hline
        (60,60,80) & 0 & 52 & 30 &7  & {\color{red}0} & 42 & 20 \\
        \hline
        (60,80,80) & 0 & {\color{red}0} & 50 &27  & 5 & {\color{red}0} & 40 \\
        \hline
        (80,80,80) & 20 & 12 & {\color{red}0} &47 & 25 & 2 & 0 \\
        \hline
    \end{tabular}

    \label{tab8}
\end{table}

From Tables \ref{tab1}-\ref{tab4}, it can be seen that once the separation times between the landing/takeoff aircraft are fixed, in most of the scenarios, each of the takeoff/landing aircraft is relevant to an aircraft with different operation task, which form block subsequences of capacity $1/(*)$. There are only a very small number of scenarios marked in red such that each of the takeoff/landing aircraft is relevant to only an aircraft with the same operation task and block subsequences are not formed. For example, consider an optimal matching problem between two sequences, $\phi_0=\langle Tcf_1, Tcf_2, Tcf_3\rangle$ and $\phi_1=\langle Tcf_4, Tcf_5, Tcf_6\rangle$, where $S_1(\phi_0)=S_4(\phi_1)$, $Y_{12}=Y_{23}=90$, $Y_{45}=60$, and $Y_{56}=80$. If we let $S_{12}(\phi_0)=S_{23}(\phi_0)=90$, $S_{45}(\phi_1)$ and $S_{56}(\phi_1)$ should be taken as $S_{45}(\phi_1)=S_{56}(\phi_1)=90$, whereas if we let $S_{45}(\phi_1)=60$ and $S_{56}(\phi_1)=80$, $S_{12}(\phi_0)$ should be taken as $S_{12}(\phi_0)=140$, and $S_{23}(\phi_0)$ can be taken as $S_{23}(\phi_0)=90$ if there are no other aircraft. In this example, it can be seen that the matching between a landing sequence and a takeoff sequence is not very complicated and the number of all possible block subsequences is $2$.

From Tables \ref{tab1}-\ref{tab8}, it can also be seen that the ranges of the time length and the time increments of block subsequences are both not large as well, in particular when the T-block subsequences are considered.
As a matter of fact, based on the calculations for more general scenarios including those listed in Tables \ref{tab1}-\ref{tab8}, it is found that if the number of breakpoint aircraft in a block subsequence is fixed, the number of all possible block subsequences of different class combinations and the ranges of the time lengths and the takeoff/landing time increments of block subsequences are all not large under the minimum separation time standards at Heathrow Airport and for the RECAT-EU system in the absence of resident-point aircraft. This property can be fully explored to find the optimal sequence when a landing sequence and a takeoff sequence are matched.

\section{Algorithms}
In the previous sections, the properties of the aircraft sequence and some typical scenarios were discussed. In this section, we will propose algorithms to find the optimal solution for the optimization problem (\ref{optim1}).

\subsection{Some necessary results}

In this subsection, we present some necessary definitions, lemmas and theorems for the algorithms.

\begin{definition}{\rm (Insertion operation, extraction operation and transformation) For two adjacent aircraft $Tcf_i$ and $Tcf_j$ in an aircraft sequence $\phi$, if an aircraft $Tcf_k$ is inserted between aircraft $Tcf_i$ and $Tcf_j$, it is said that an insertion operation is performed on aircraft $Tcf_k$ between aircraft $Tcf_i$ and $Tcf_j$. For three consecutive aircraft $Tcf_i, Tcf_j$ and $Tcf_k$ in an aircraft sequence $\phi$, if aircraft $Tcf_j$ is extracted from the aircraft sequence $\phi$, it is said that an extraction operation is performed on aircraft $Tcf_k$ from the sequence $\phi$. For two aircraft sequences $\phi_1$ and $\phi_2$, if $\phi_1$ is converted into $\phi_2$ through a series of insertion and extraction operations, it is said that $\phi_2$ is a transformation of $\phi_1$. Further, if $F(\phi_1, Sr(\phi_1))=F(\phi_2, Sr(\phi_2))$ and $\phi_2$ is a transformation of $\phi_1$, it is said that $\phi_2$ is an equivalent transformation of $\phi_1$. }\end{definition}

From the above definitions, when an aircraft sequence is a transformation of another aircraft sequence, then the two sequences are composed of the same group of aircraft.

\begin{definition}{\rm (Drift operation) Consider an aircraft sequence $\Phi_a=\langle Tcf_1, Tcf_2, \cdots, Tcf_n\rangle$. Generate a new sequence $\Phi_b$ by moving a block of aircraft $Tcf_k, Tcf_{k+1}, \cdots, Tcf_{k+h}$ to be between aircraft $Tcf_{i}$ and $Tcf_{i+1}$. If the time intervals between any two aircraft of $Tcf_k, Tcf_{k+1}, \cdots, Tcf_{k+h}$ are the same in $\Phi_a$ and $\Phi_b$, it is said that a drift operation is performed on the aircraft $Tcf_k, Tcf_{k+1}, \cdots, Tcf_{k+h}$.
}\end{definition}

\begin{definition}{\rm (Minimum insertion time increment) Consider an insertion of an aircraft $Tcf_i$ between two adjacent aircraft $Tcf_{k}$ and $Tcf_j$. The quantity $Y_{ki}+Y_{ij}-Y_{kj}$ is said to be the minimum insertion time increment of aircraft $Tcf_i$ with respect to aircraft $Tcf_{k}$ and $Tcf_j$.
}\end{definition}

\begin{definition}{\rm (Minimum extraction time increment)  Consider an extraction of an aircraft $Tcf_i$ from a sequence $\langle Tcf_k, Tcf_i, Tcf_j \rangle$.  The quantity $Y_{ki}+Y_{ij}-Y_{kj}$ is said to be the minimum extraction time increment of aircraft $Tcf_i$ with respect to aircraft $Tcf_{k}$ and $Tcf_j$.}\end{definition}

It can be easily seen that the minimum insertion/extraction time increment is related to only the aircraft classes. Since the number of the aircraft classes is $\eta$, all possible minimum insertion/extraction time increments can easily be calculated out.

\begin{definition}\label{defipp}{\rm (The breakpoint-drift equivalent transformation) Consider two landing/takeoff sequences $\Phi_a=\langle\phi_{a1}, \phi_{a2}, \cdots, \phi_{as}\rangle$ and $\Phi_b=\langle\phi_{b1}, \phi_{b2}, \cdots, \phi_{bs}\rangle$ for a positive integer $s$, where each $\phi_{ai}=\langle Tcf_{ai1}, Tcf_{ai2}, \cdots, Tcf_{aic_i}\rangle$ and each $\phi_{bi}=\langle Tcf_{bi1}, Tcf_{bi2}, \cdots, Tcf_{bih_i}\rangle$ are both class-monotonically-decreasing sequences for positive integers $c_i$, $h_i$, and all $i\in \{1, 2, \cdots, s\}$, such that $cl_{ajc_j}<cl_{a(j+1)1}$ and $cl_{bjh_j}<cl_{b(j+1)1}$ for all $j\in \{1, 2, \cdots, s-1\}$. Suppose that Assumptions  \ref{ass1.3.1}-\ref{ass2.4} hold for both aircraft sequences $\Phi_a$ and $\Phi_b$.
If $Y_{(ajc_j)(a(j+1)1)}=Y_{(bjc_j)(b(j+1)1)}$ for all $j\in \{1, 2, \cdots, s-1\}$ and $F(\Phi_a)=F(\Phi_b)$, $\Phi_b$ is one breakpoint equivalent transformation of $\Phi_a$. The set of all the breakpoint equivalent transformations of $\Phi_a$ is said to be the breakpoint-drift equivalent transformation set of $\Phi_a$, denoted by $\mathrm{Bet}(\Phi_a)$, and $\Phi_a$ is said to be a sequence basis of $\mathrm{Bet}(\Phi_a)$.
} \end{definition}

For the optimization problem (\ref{optim1}), the set of the optimal sequences might be composed of the breakpoint-drift equivalent transformation sets spanned by multiple optimal sequence basis. For example, consider a landing sequence $\phi_a=\langle Tcf_1, Tcf_2, Tcf_3, Tcf_4, Tcf_5\rangle$ and $\phi_b=\langle Tcf_1, Tcf_3, Tcf_4, Tcf_2, Tcf_5\rangle$ under Assumptions \ref{ass1.3.1}-\ref{ass1.4}, where $Y_{12}=180, Y_{13}=160, Y_{23}=Y_{25}=Y_{34}=60$, $Y_{45}=120$, and $Y_{42}=140$. It is clear that $F(\phi_a)=F(\phi_b)$ but the separation times between the breakpoint aircraft and their trailing aircraft in $\phi_a$ and $\phi_b$ are different.

From this example, it can also be seen that
the bases of all the breakpoint-drift equivalent transformation sets can be interconverted by focusing on the orders of the breakpoint aircraft and their trailing aircraft.

\begin{lemma}\label{lemmadd}{\rm Consider a landing/takeoff sequence $\Phi_a=\langle Tcf_1, Tcf_2, \cdots, Tcf_n\rangle$. Generate a sequence $\Phi_b$ by moving $Tcf_i$ to be between $Tcf_{k}$ and $Tcf_{k+1}$ for two integers $0<k,i\leq n$. Suppose that Assumptions \ref{ass1.3.1}-\ref{ass2.4} hold for all sequences.
If $k+1<i$, then $S_j(\Phi_a)\leq S_j(\Phi_b)$ for $k+1<j<i$ and if $k>i$, then $S_j(\Phi_a)\geq S_j(\Phi_b)$ for $i<j<k$.
}\end{lemma}
\noindent{Proof:} This lemma is obvious and its proof is omitted.

Lemma \ref{lemmadd} shows the impact of aircraft order adjustments on the other aircraft. This lemma can be used to easily obtain the minimum and maximum adjustment ranges of aircraft by first making adjustments that have little impact on other aircrafts, which can be applied to find all breakpoint-drift equivalent transformations of an optimal sequence.
For example, based on $\Phi_a$ in Lemma \ref{lemmadd}, suppose that a sequence $\Phi_{b1}$ containing a subsequence $\langle Tcf_i, Tcf_{i+3}\rangle$ can be generated by moving the aircraft $Tcf_{i+1}, Tcf_{i+2}$ to be between $Tcf_{k_0}$ and $Tcf_{k_0+1}$ and between $Tcf_{k_1}$ and $Tcf_{k_1+1}$ for $k_0<k_1<i$. From Lemma \ref{lemmadd}, we can first move $Tcf_{i+1}$ and then move $Tcf_{i+2}$ since the movement of $Tcf_{i+1}$ has less effects on $Tcf_{i+2}$. If we first perform the movement of the aircraft $Tcf_{i+2}$ and keep the order of $Tcf_{i+1}$ unchanged, the generated sequence $\Phi_c$ might not exist because a possible case satisfying that $S_{i+2}(\Phi_c)\notin [f^{\min}_{i+2}, f^{\max}_{i+2}]$ and $S_{i+2}(\Phi_{b1})\in [f^{\min}_{i+2}, f^{\max}_{i+2}]$ might occur.

{In the following, we give an example to show how to obtain the maximum number of aircraft that can be moved to be before one aircraft based on Lemma \ref{lemmadd} and the analysis of the above example. Consider a sequence $\phi=\langle Tcf_1, Tcf_2, \cdots, Tcf_{s_0}, Tcf_{s_0+1}, \cdots, Tcf_n\rangle$, where $\phi_1=\langle Tcf_1, Tcf_2, \cdots, Tcf_{s_0}\rangle$ is a class-monotonically-decreasing sequence. The objective is to move the aircraft $Tcf_{s_0+1}, \cdots, Tcf_n$ as many as possible to be before $Tcf_{s_0}$ without generating new breakpoint aircraft. To this end, we first rearrange the aircraft $Tcf_{s_0+1}, \cdots, Tcf_n$ in descending order of their aircraft classes, where the aircraft of the same class are rearranged in ascending order of their earliest landing/takeoff times. Let $\phi_2=\langle Tcf^a_{s_0+1}, \cdots, Tcf^a_n\rangle$ denote the generated sequence of the aircraft $Tcf_{s_0+1}, \cdots, Tcf_n$. Attempt to insert the aircraft in the order of $\phi_2$ into the subsequence $\phi_1$ and we can roughly find the maximum number of aircraft that can be moved to be before aircraft $Tcf_{s_0}$.}

\begin{lemma}\label{llemma}{\rm Consider a landing/takeoff sequence $\Phi_a=\langle Tcf_1, Tcf_2, \cdots, Tcf_n\rangle$. Suppose that $(\Phi_a, Sr(\Phi_a))$ is an optimal solution for the optimization problem (\ref{optim1}), $\Phi_b$ is a sequence generated by moving an aircraft $Tcf_i$ to be between $Tcf_{k}$ and $Tcf_{k+1}$ in $\Phi_a$, and $\Phi_c$ is a sequence generated by moving an aircraft $Tcf_j$ to be between $Tcf_{h}$ and $Tcf_{h+1}$ in $\Phi_b$ for $i<j<k<h$ and Assumptions \ref{ass1.3.1}-\ref{ass2.4} hold for all sequences. The following statements hold.

(1) $F(\Phi_b)\geq F(\Phi_a)$ and $F(\Phi_c)\geq F(\Phi_a)$.

(2) If $j-i>1$, $k-j>1$, $h-k>1$ and $F(\Phi_c)=F(\Phi_a)$, then $F(\Phi_b)=F(\Phi_a)$.

}\end{lemma}
\noindent{Proof:} The statement (1) holds naturally from the optimality of the sequence $\Phi_a$. Suppose that $j-i>1$, $k-j>1$, and $h-k>1$. From Theorems \ref{theorem1.1} and \ref{theorempp2.1}, the movement of $Tcf_i$ from $\Phi_a$ to $\Phi_b$ and the movement of $Tcf_j$ from $\Phi_b$ to $\Phi_c$ have no effects with each other. Generate a new sequence $\Phi_d$ by moving an aircraft $Tcf_j$ to be between $Tcf_{h}$ and $Tcf_{h+1}$ in $\Phi_a$ under Assumptions \ref{ass1.3.1}-\ref{ass2.4}. It follows that $F(\Phi_c)-F(\Phi_a)=F(\Phi_b)-F(\Phi_a)+F(\Phi_d)-F(\Phi_a)$. If $F(\Phi_b)>F(\Phi_a)$ and $F(\Phi_c)=F(\Phi_a)$, it follows that $F(\Phi_d)<F(\Phi_a)$, which is a contradiction. Therefore, $F(\Phi_b)=F(\Phi_a)$.

\begin{remark}{\rm From Lemmas \ref{lemmadd} and \ref{llemma} and the definition of the breakpoint-drift equivalent transformation set, the set of all the optimal sequences can be obtained by spanning one optimal sequence and adjusting the orders of the breakpoint aircraft and their trailing aircraft for the optimization problem (\ref{optim1}). Lemma \ref{llemma} also shows that any independent operations between the optimal sequences should not change the value of the objective function $F(\cdot)$, which is a very important property to obtain all optimal sequences.}\end{remark}

\begin{theorem}\label{theorem1ds}{\rm Suppose that $(\phi, Sr(\phi))$ is an optimal solution of the optimization problem (\ref{optim1}) for aircraft
 $Tcf_{1},Tcf_{2},$ $ \cdots, Tcf_{i}$, where $f_1^{\min}\leq f_2^{\min}\leq \cdots\leq f_i^{\min}$. Then there are no resident-point aircraft that land/take off during the time interval $[f_{i}^{\min}, +\infty]$.}\end{theorem}

\noindent{Proof:} Note that $f_1^{\min}\leq f_2^{\min}\leq \cdots\leq f_i^{\min}$ and there is at least a sequence for all aircraft to land/take off without conflicts. Without loss of generality, denote $\phi=\langle \phi_0, \phi_1\rangle$, where $\phi_1=\langle Tcf_{k_0}, Tcf_{k_1}, \cdots, Tcf_{k_h}\rangle$ with $k_0=i$ for some integer $h>0$.
If aircraft $Tcf_{k_j}$ is an resident-point aircraft for some integer $0\leq j\leq h$, thus $S_{k_{j-1}k_j}(\phi)-Y_{k_{j-1}k_j}>0$. Let the landing/takeoff times of the aircraft $Tcf_{k_j}, Tcf_{k_{j+1}}, \cdots, Tcf_{k_h}$ be advanced by $S_{k_{j-1}k_j}(\phi)-Y_{k_{j-1}k_j}>0$. Let $\overline{Sr}(\phi)$ denote the new landing/takeoff times. It is clear that $F(\phi, Sr(\phi))>F(\phi, \overline{Sr}(\phi))$, implying that
$(\phi, Sr(\phi))$ is not the optimal solution of $Tcf_{1},Tcf_{2}, \cdots, Tcf_{i}$, which is a contradiction.

 \begin{remark}{\rm Using Theorem \ref{theorem1ds}, we can roughly determine which aircraft are resident-point aircraft and calculate out the total resident time. During a given time interval, when the aircraft operate at a high density, from the proof of Theorem \ref{theorem1ds}, there might be usually no resident-point aircraft. }\end{remark}

{\begin{remark}{\rm In practical systems,  the earliest landing/takeoff times are usually unable to be advanced while the choice of the latest landing/takeoff times involves more artificial factors and can be postponed without the occurrence of sudden change in most cases.
{In Algorithm 1, the optimal sequences for the optimization problem (\ref{optim1}) are searched in ascending order of their earliest landing/takeoff times, which is shown in Theorem \ref{theorem1ds} to be able to avoid the occurrence of the resident-point aircraft after the new inserted aircraft.
This property is very useful and important for the optimal convergence of Algorithm 1.}}\end{remark}}

In the following proposition, we make an estimation about the number of breakpoint aircraft for a landing/takeoff sequence.

\begin{proposition}{\rm Consider a landing/takeoff sequence $\phi=\langle Tcf_1, Tcf_2, \cdots, Tcf_n\rangle$. Suppose that $(\phi, Sr(\phi))$ is an optimal solution of the optimization problem (\ref{optim1}) for aircraft $Tcf_{1},Tcf_{2}, \cdots, Tcf_{n}$, and $[f^0, f^1]\subset \bigcap_{i=1}^n [f_i^{\min},$ $f_i^{\max}]$. In the time interval $[f^0, f^1]$, the number of breakpoint aircraft lies in $[0,2]$. }\end{proposition}

\noindent{Proof:} Since $[f^0, f^1]\subset \bigcap_{i=1}^n [f_i^{\min},f_i^{\max}]$, the orders of aircraft in the time interval $[f^0, f^1]$ can be adjusted arbitrarily.
Thus, there might be $0$ breakpoint aircraft. Note that the aircraft in the time interval $[t_0, f^0)$ and the aircraft in the time interval $[f^1, +\infty)$ might form  class-monotonically-decreasing sequences that finally generate $1$ or $2$ breakpoint aircraft in the time interval $[f^0, f^1]$. Summarizing all the analysis, in the time interval $[f^0, f^1]$, the number of breakpoint aircraft lies in $[0,2]$.

\begin{remark}{\rm When a sequence contains a resident-point aircraft, we can adjust aircraft orders to eliminate the resident-point aircraft or reduce their influences. When there are aircraft that cannot be moved forwards or backwards due to the constraints of time windows, there might be two methods to relax the time window constraints. One is
 to adjust the orders of aircraft by extraction and insertion, and the other is to select some proper aircraft according to the sequence features and possibly treat them as new aircraft for processing. We give a simple example to illustrate how to select aircraft according to the sequence features. Consider a landing sequence $\phi=\langle Tcf_1, Tcf_2, Tcf_3, Tcf_4\rangle$, where $cl_1=4$, $cl_2=3$, $cl_3=2$, $cl_4=4$, $S_{12}=Y_{12}$, $S_{23}=Y_{23}$, $S_1=t_0$, $S_{34}=Y_{34}$ and $S_3=f_{max}^3$. If we want to insert a new aircraft $Tcf_5$ with $f_{max}^5=f_{max}^3$ into the subsequence $\phi$, we can calculate the time increment of the sequence when $Tcf_1$ or $Tcf_2$ is extracted from the subsequence $\langle Tcf_1, Tcf_2, Tcf_3\rangle$ and inserted into the subsequence $\langle Tcf_3, Tcf_4\rangle$ to make a proper choice.
}\end{remark}

{\begin{theorem}{\rm Consider a mixed landing and takeoff sequence $\Phi_a=\langle Tcf_1, Tcf_2, Tcf_3, Tcf_4\rangle$ on a same runway, where $S_{12}(\Phi_a)\geq Y_{12}$, aircraft $Tcf_3$ is relevant to aircraft $Tcf_1$ or $Tcf_2$ and aircraft $Tcf_4$ is relevant to aircraft $Tcf_2$ or $Tcf_3$. Generate a new sequence $\Phi_b$ by inserting an aircraft $Tcf_5$ to be between aircraft $Tcf_2$ and $Tcf_3$, where $S_{1}(\Phi_a)=S_{1}(\Phi_b)$, $S_{2}(\Phi_a)=S_{2}(\Phi_b)$, aircraft $Tcf_5$ is relevant to aircraft $Tcf_1$ or $Tcf_2$, aircraft $Tcf_3$ is relevant to aircraft $Tcf_2$ or $Tcf_5$ and aircraft $Tcf_4$ is relevant to aircraft $Tcf_5$ or $Tcf_3$. Let $\omega_{min}(\Phi_a, Tcf_5)=\min\{S_3(\Phi_b)-S_3(\Phi_a),S_4(\Phi_b)-S_4(\Phi_a)\}$. The following statements hold.

(1) Suppose that aircraft $Tcf_2$, $Tcf_3$ and $Tcf_5$ are all landing (takeoff) aircraft. It follows that $\omega_{min}(\Phi_a, Tcf_5)=Y_{25}+Y_{53}-Y_{23}$.

(2) Suppose that aircraft $Tcf_2$ and $Tcf_5$ are both landing (takeoff) aircraft and $Tcf_3$ is a takeoff (landing) aircraft. 
It follows that $\omega_{min}(\Phi_a, Tcf_5)\geq\min\{T_D+D_T-Y_{13},T_D+D_T-Y_{24},0\}+Y_{25}$.

(3) Suppose that aircraft $Tcf_3$ and $Tcf_5$ are both landing (takeoff) aircraft and $Tcf_2$ is a takeoff (landing) aircraft. It follows that $\omega_{min}(\Phi_a, Tcf_5)\geq\min\{T_D+D_T-Y_{13},T_D+D_T-Y_{24},0\}+Y_{53}$.

(4)  Suppose that aircraft $Tcf_2$ and $Tcf_3$ are both landing (takeoff) aircraft and $Tcf_5$ is a takeoff (landing) aircraft. It follows that $\omega_{min}(\Phi_a, Tcf_5)\geq \max\{D_T+T_D-Y_{23},0\}$.
}\end{theorem}

\noindent{Proof:} In this proof, we only discuss the statement (2) and statements (1)(3)(4) can be proved in a similar way. When $Tcf_4$ is a takeoff (landing) aircraft and $Tcf_3$ is relevant to $Tcf_2$ in $\Phi_a$, $\omega_{min}(\Phi_a, Tcf_5)=Y_{25}$. When $Tcf_4$ is a takeoff (landing) aircraft and $Tcf_3$ is relevant to $Tcf_1$ in $\Phi_a$, $T_D+D_T-Y_{13}\leq0$ and $\omega_{min}(\Phi_a, Tcf_5)\geq Y_{25}+T_D+D_T-Y_{13}$. That is, $\omega_{min}(\Phi_a, Tcf_5)\geq Y_{25}+\min\{T_D+D_T-Y_{13},0\}$. Similarly, when $Tcf_4$ is a landing (takeoff) aircraft, it can be proved that $S_3(\Phi_b)-S_3(\Phi_a)\geq \min\{T_D+D_T-Y_{13},0\}+Y_{25}$ when $Tcf_4$ is relevant to $Tcf_3$ both in $\Phi_a$ and $\Phi_b$, and $S_4(\Phi_b)-S_4(\Phi_a)\geq \min\{T_D+D_T-Y_{13},T_D+D_T-Y_{24},0\}+Y_{25}$ when $Tcf_4$ is relevant to $Tcf_2$ in $\Phi_a$ and relevant to $Tcf_5$ or $Tcf_3$ in $\Phi_b$. That is, $\omega_{min}(\Phi_a, Tcf_5)\geq\min\{T_D+D_T-Y_{13},T_D+D_T-Y_{24},0\}+Y_{25}$.

\begin{remark}{\rm This theorem gives lower bounds for the time increments of the whole sequence when an aircraft is inserted, which can be used in Algorithm 1. Moreover, when defining $\omega_{min}$, since the aircraft $Tcf_3$ and $Tcf_4$ both have the possibilities of having no path from the last aircraft to themselves, we adopt the smaller value between $S_3(\Phi_b)-S_3(\Phi_a)$ and $S_4(\Phi_b)-S_4(\Phi_a)$ as the value of $\omega_{min}$ so as to implement a full space search for Algorithm 1.}\end{remark}}

\subsection{Algorithm 1}

From the definition of the objective function $F(\cdot)$, when the landing/takeoff scheduling problem is considered, the optimization problem (\ref{optim1}) can be rewritten as
\begin{eqnarray}\label{optim2}\begin{array}{lll}\min \sum_{i=1}^{n-1}(S_{i+1}-S_i-T_0)+(S_1-t_0)\\
\mbox{Subject to}~~~S_k \in Tf_k=[f_k^{\min}, f_k^{\max}]\\
\hspace {2cm}k =1, \cdots, n,\\
\hspace {2cm}S_{i+1}-S_i \geq Y_{{i}({i+1})},\\ \hspace {2cm}i=1, \cdots, n-1.\end{array}\end{eqnarray}

Note that $S_{i+1}-S_i-T_0\geq 0$ for any two aircraft $Tcf_i$ and $Tcf_{i+1}$. To minimize the objective function $F(\cdot)$, we can focus on the consecutive aircraft whose separation times are larger than $T_0$, including the breakpoint aircraft and the resident-point aircraft as well as the separation times between the consecutive aircraft of the same class that are larger than $T_0$. Motivated by this observation, we propose the following algorithm for landing scheduling problem, takeoff scheduling problem, scheduling problem of landing and takeoff aircraft on a same runway.

\textbf{Algorithm 1.} Suppose that there are totally $n$ aircraft, denoted by $Tcf_{01}, Tcf_{02}, \cdots, Tcf_{0n}$, and there is at least a feasible sequence for them to land without conflicts.

Step 1. Arrange the aircraft in ascending order of their earliest landing/takeoff times, denoted by $\langle Tcf_{11}, Tcf_{12}, \cdots,$ $ Tcf_{1n}\rangle$.

Step 2. Search for the optimal sequence of aircraft $Tcf_{11}$ and $Tcf_{12}$.

When no confusion arises, it is assumed by default that the mentioned sequence or subsequence contains the information of landing/takeoff times in Algorithms 1 and 2.

Step 3. Suppose that $\phi_2=\langle Tcf_1, Tcf_2\rangle$ is the optimal sequence of aircraft $Tcf_{11}$ and $Tcf_{12}$. Search for the optimal sequence of aircraft $Tcf_1$, $Tcf_2$ and $Tcf_{13}$.

Step $i$, $i=4,5, \cdots,n$. Suppose that $\phi_{i-1}=\langle Tcf_1,Tcf_2, \cdots, Tcf_{i-1}\rangle$ is an optimal sequence for aircraft $Tcf_{11},Tcf_{12}, \cdots, Tcf_{1(i-1)}$. Search for the optimal sequences of aircraft $Tcf_{11},Tcf_{12}, \cdots, Tcf_{1i}$ based on the obtained theoretical results. The main idea is to exclude most of the non-optimal solutions and narrow down the search for the optimal solutions within a small range. The main steps are as follows.

$(i-1)$. Based on the obtained optimal sequence $\phi_{i-1}$, by a series of equivalent transformations, record all possible optimal sequences based on the bases of all the breakpoint-drift equivalent transformation sets and the takeoff-landing/landing-takeoff transitions, and represent their formed set as $\Lambda_{i-1}$ for aircraft $Tcf_{11},Tcf_{12}, \cdots, Tcf_{1(i-1)}$. To find the optimal sequences of aircraft $Tcf_{11},Tcf_{12}, \cdots, Tcf_{1i}$, we can first insert aircraft $Tcf_{1i}$ between any two adjacent aircraft for each sequence in the set $\Lambda_{i-1}$, denoted by $\phi_{(i-1)k}$, under the time window constraints without changing the orders of the aircraft in $\phi_{(i-1)k}$. Let $f^i_{inc}$ be the smallest value of the objective functions $F(\cdot)$ among all the generated new sequences.

For the landing/takeoff scheduling problem, from Theorems \ref{lemma1.3s}, \ref{lemma2.3s} and \ref{lemma3.31s}, the value of the objective function $F(\cdot)$ is heavily related to the separation times between breakpoint aircraft and their trailing aircraft. To obtain or record all possible optimal sequences based on the bases of all the breakpoint-drift equivalent transformation sets, emphasis should be laid on the breakpoint aircraft. For the scheduling problem of landing and takeoff aircraft on a same runway, the number of breakpoint aircraft is usually small but the number of the takeoff-landing/landing-takeoff transitions might be large, and emphasis should be laid on the takeoff-landing/landing-takeoff transitions.

{The goal of this step is actually to find all possible combinations of four consecutive aircraft of different classes in all sequences of the set $\Lambda_{i-1}$ to obtain the smallest value of the objective functions $F(\cdot)$}, where each combination corresponds to a time increment when the aircraft $Tcf_{1i}$ is inserted. Starting from this goal, we can reduce some computation amount by giving up some unnecessary operations. For example,
since we are concerned about only the time increment when the aircraft $Tcf_{1i}$ is inserted into $\phi_{(i-1)k}$, we need not to repeatedly consider the insertion of aircraft $Tcf_{1i}$ into the same class combination of aircraft and the same separation times, and focus on only the different class combination of aircraft or different separation times, which might significantly reduce the computation amount.

Time increment would be used frequently in our algorithm. To illustrate its role specifically, we give an example. Generate two landing sequences $\phi_1=\langle Tcf_1, Tcf_4, Tcf_2, Tcf_3\rangle$ and $\phi_2=\langle Tcf_1, Tcf_2, Tcf_4, Tcf_3\rangle$ by inserting aircraft $Tcf_4$ into a sequence $\phi_0=\langle Tcf_1, Tcf_2, Tcf_3\rangle$. The time increment of the separation time between aircraft $Tcf_1$ and $Tcf_2$ in $\phi_1$ with respect to $\phi_0$ is $\Delta_1=S_{14}(\phi_1)+S_{42}(\phi_1)-S_{12}(\phi_0)$ and the time increment of the separation time between aircraft $Tcf_1$ and $Tcf_2$ in $\phi_1$ with respect to $\phi_0$ is $\Delta_2=S_{24}(\phi_2)+S_{43}(\phi_2)-S_{23}(\phi_0)$. If $S_{23}(\phi_1)=S_{23}(\phi_0)$ and $S_{12}(\phi_2)=S_{12}(\phi_0)$, then $F(\phi_2)-F(\phi_1)=\Delta_2-\Delta_1$.

{$(i-2)$. Construct a sequence set $F_{inc}=\{\langle h_1, h_2, h_3, h_4\rangle \mid h_1, h_2, h_3, h_4\in \{1,2,\cdots,i-1\}\}$ composed of all the sequences of length $4$ such that \begin{eqnarray}\label{eaf1}\begin{array}{lll}\omega_{min}(\langle Tcf_{h_1}, Tcf_{h_2}, Tcf_{h_3}, Tcf_{h_4}\rangle, Tcf_{1i})\\\hspace{-0.3cm}< f^i_{inc}-F(\phi_{i-1})\end{array}\end{eqnarray} for any $\langle h_1, h_2, h_3, h_4\rangle\in F_{inc}$. Note here that since the minimum separation times are related to only the classes and the operation tasks of the aircraft, we can classify the aircraft to form the set $F_{inc}$ by calculating the inequality (\ref{eaf1}) according to the classes and the operation tasks of the aircraft.

When the landing/takeoff scheduling problem is considered and each aircraft is relevant to its leading aircraft, the inequality $(\ref{eaf1})$ can be simplified into $Y_{h_2(1i)}+Y_{(1i)h_3}-Y_{h_2h_3}<f^i_{inc}-F(\phi_{i-1})$. When the scheduling problem of landing and takeoff aircraft on a same runway is considered,
the aircraft $Tcf_{h_3}$ might be relevant to $Tcf_{h_1}$ and not relevant to aircraft $Tcf_{h_2}$ in the original sequences, and the aircraft $Tcf_{h_4}$ might be relevant to $Tcf_{1i}$ and not relevant to aircraft $Tcf_{h_3}$ in the new sequences. This means that the time increment of an insertion of aircraft might be related to four consecutive aircraft except the aircraft $Tcf_{1i}$, which is the reason why four aircraft $Tcf_{h_1}, Tcf_{h_2}, Tcf_{h_3}, Tcf_{h_4}$ are involved. But it should be noted that the form of (\ref{eaf1}) can be simplified according to the aircraft relevance.}

When there exist resident-point aircraft, we need to additionally consider the influences of the resident times to define the set $F_{inc}$. For example, consider an insertion of an aircraft $Tcf_1$ between adjacent aircraft  $Tcf_2$ and $Tcf_3$ in a landing sequence $\phi$, where $S_{23}(\phi)-Y_{23}>0$ and $S_{23}(\phi)<Y_{21}+Y_{13}$.  In the new generated sequence after the insertion of $Tcf_1$, if $Tcf_1$ and $Tcf_3$ are relevant to $Tcf_2$ and $Tcf_1$ respectively, the time increment of the aircraft $Tcf_1$ should be $Y_{21}+Y_{13}-S_{23}(\phi)$.

$(i-3).$ According to the elements of the set $F_{inc}$, adjust the aircraft orders in the sequence $\phi_{(i-1)k}$ and insert the aircraft $Tcf_{1i}$ to generate a new sequence $\bar{\phi}_{ik}$ under the condition $F(\bar{\phi}_{ik})<f^i_{inc}$ such that $\bar{\phi}_{ik}$ contains the subsequence $\langle Tcf_{h_1}, Tcf_{h_2}, Tcf_{1i}, Tcf_{h_3}, Tcf_{h_4}\rangle$ and $\langle h_1, h_2, h_3, h_4\rangle\in F_{inc}$. Minimize $F(\bar{\phi}_{ik})$ based on the obtained theoretical results, and compare the values of all possible $F(\bar{\phi}_{ik})$ so as to find
the optimal sequence of aircraft $Tcf_{11},Tcf_{12}, \cdots, Tcf_{1i}$.

To illustrate this step in a more specific way, consider a group of landing aircraft with $T_0=60$ and $\delta=8$ under the minimum separation time standards at Heathrow Airport. From the previous analysis in (i-2), the inequality (\ref{eaf1}) is equivalent to the form \begin{eqnarray}\label{ear2}Y_{p(1i)}+Y_{(1i)j}-Y_{pj}<f^i_{inc}-F(\phi_{i-1}).\end{eqnarray} If $f^i_{inc}-F(\phi_{i-1})=T_0$, by simple calculations, the inequality (\ref{ear2}) holds when $(cl_p, cl_{1i}, cl_j)=(2,5,4), (3,6,5), (4,6,5)$. Clearly, when $(cl_p, cl_{1i}, cl_j)=(2,5,4), (3,6,5), (4,6,5)$, each aircraft $Tcf_p$ is a breakpoint aircraft and the aircraft $Tcf_{1i}$ is needed to be inserted to be between aircraft $Tcf_{p}$ and $Tcf_{j}$. To generate a subsequence $\langle Tcf_{p}, Tcf_{1i}, Tcf_{j}\rangle$, we can focus on the breakpoint aircraft and their trailing aircraft, and
based on the theoretical results in previous sections, by adjusting the aircraft orders properly, the optimal sequence can be found under the condition that $F(\bar{\phi}_{ik})<f^i_{inc}$ in an easy way. 
If $f^i_{inc}-F(\phi_{i-1})=T_0+\delta$, the generation of a subsequence $\langle Tcf_{p_1}, Tcf_{1i}, Tcf_{j_1}\rangle$ with $\langle {p_1}, {j_1}\rangle\in F_{inc}$ would result in a new breakpoint aircraft or some changes of the separation times between existing breakpoint aircraft and their trailing aircraft. Considering the values of the separation times between all possible breakpoint aircraft and their trailing aircraft, note that $T_{45}=T_0+\delta$ and $T_{ij}\geq 1.5T_0$ for all $1\leq i<j\leq 6$ with $(i,j)\neq (4,5)$. We can generate the whole sequence based on the subsequence $\langle Tcf_{p_1}, Tcf_{1i}, Tcf_{j_1}\rangle$ by considering the total time increments caused by the insertion of the aircraft $Tcf_{1i}$ under the condition that $F(\bar{\phi}_{ik})<f^i_{inc}$.

Since each $\phi_{(i-1)k}$ is an optimal sequence for aircraft $Tcf_{11},Tcf_{12}, \cdots, Tcf_{1(i-1)}$, 
if there is no utilizable subsequence $\langle Tcf_{h_1}, Tcf_{h_2}, Tcf_{h_3}, Tcf_{h_4}\rangle$ formed by four consecutive aircraft of $Tcf_{h_1}, Tcf_{h_2}, Tcf_{h_3}, Tcf_{h_4}$ in all $\phi_{(i-1)k}$, the generation of $\langle Tcf_{h_1}, Tcf_{h_2}, Tcf_{h_3}, Tcf_{h_4}\rangle$ in $\phi_{(i-1)k}$ would increase the value of the objective function $F(\cdot)$ and $F(\bar{\phi}_{ik})\geq f^i_{inc}$ when $\omega_{min}(\langle Tcf_{h_1}, Tcf_{h_2}, Tcf_{h_3}, Tcf_{h_4}\rangle, Tcf_{1i})\geq f^i_{inc}-F(\phi_{i-1})$. 
This is the reason why only the aircraft combinations satisfying (\ref{eaf1}) are taken into account. Moreover, the optimal insertion of $Tcf_{1(i-1)}$ might usually generate a new breakpoint/landing-takeoff/takeoff-landing transition or change an existing separation time between a breakpoint aircraft and its trailing aircraft with some necessary aircraft order adjustments. These properties and related discussions in previous sections can be used to reduce the computation amount of Algorithm 1.

Based on the above analysis, we can combine the following cases to discuss the optimality of the sequence $\bar{\phi}_{ik}$.

Scenario 1. Landing/takeoff sequences. Let $Nb_{\bar{\phi}_{ik}}$ and $Nb_{\phi_{(i-1)k}}$ denote the numbers of breakpoint aircraft in the sequences $\bar{\phi}_{ik}$ and $\phi_{(i-1)k}$.

Case 1.1. $Nb_{\bar{\phi}_{ik}}\leq Nb_{\phi_{(i-1)k}}$.

Case 1.2. $Nb_{\bar{\phi}_{ik}}>Nb_{\phi_{(i-1)k}}$.

Scenario 2. Mixed landing and takeoff sequences on a same runway. Let $Nt_{\bar{\phi}_{ik}}$ and $Nt_{\phi_{(i-1)k}}$ denote the numbers of landing-takeoff and takeoff-landing transitions in the sequences $\bar{\phi}_{ik}$. and $\phi_{(i-1)k}$.

Case 2.1. $Nt_{\bar{\phi}_{ik}}\leq Nt_{\phi_{(i-1)k}}$. 

~~~Subcase 2.1.1. $Nb_{\bar{\phi}_{ik}}\leq Nb_{\phi_{(i-1)k}}$.

~~~Subcase 2.1.2. $Nb_{\bar{\phi}_{ik}}>Nb_{\phi_{(i-1)k}}$.

Case 2.2. $Nt_{\bar{\phi}_{ik}}>Nt_{\phi_{(i-1)k}}$.

~~~Subcase 2.2.1. $Nb_{\bar{\phi}_{ik}}\leq Nb_{\phi_{(i-1)k}}$.

~~~Subcase 2.2.2. $Nb_{\bar{\phi}_{ik}}>Nb_{\phi_{(i-1)k}}$.

When the number of breakpoint aircraft in $\bar{\phi}_{ik}$ is 2 larger than that in $\phi_{(i-1)k}$, $\bar{\phi}_{ik}$ is usually not the optimal sequence when the aircraft subsequence $\langle Tcf_i, Tcf_j\rangle$ with $cl_i=\rho_2-1$ and $cl_j=\rho_2$ is not involved.

{Note that $\phi_{(i-1)k}$ might be composed of multiple class-monotonically-decreasing subsequences, denoted by $\phi^1, \phi^2, \cdots, \phi^{s}$. We can consider the insertion of $Tcf_{1i}$ on each subsequence $\phi^k$ with some corresponding aircraft order adjustments so as to reduce the complexity of the analysis of the optimality of the sequence. }

{\begin{remark}{\rm Since different class combinations might correspond to different time increments, we can search for the desired aircraft of proper classes in $F_{inc}$ in ascending order of time increments to minimize the value of the objective function $F(\cdot)$. Moreover, note that the separation times between aircraft of different classes are time-invariant and
the number of the total aircraft classes is not very large. The computation cost of time increments is also not very large, which can be calculated out offline. }\end{remark}}

\begin{remark}{\rm It should be noted that in Algorithm 1, when there is no resident-point aircraft, the optimal sequence $\phi_i$ satisfies that $F(\phi_{i-1}, Sr(\phi_{i-1}))\leq F(\phi_i, Sr(\phi_i))\leq F(\phi_{i-1}, Sr(\phi_{i-1}))+Y_{(i-1)i}$ for the landing/takeoff scheduling problem, and the landing and takeoff scheduling problem on a same runway, where $T_0\leq Y_{(i-1)i}\leq 3T_0$.
For example, consider a sequence $\phi_0=\langle Tcf_1, Tcf_2\rangle$, where $Tcf_1, Tcf_2$ are both landing aircraft and $Y_{12}>T_D+D_T$. Generate a new sequence $\phi_1$ by inserting a takeoff aircraft $Tcf_3$ to be between $Tcf_1$ and $Tcf_2$, where $S_{13}(\phi_1)=T_D$ and $S_{12}(\phi_1)=S_{12}(\phi_0)$. It is clear that $F(\phi_1)-F(\phi_0)=0$.
}\end{remark}

\begin{remark}{\rm When the orders of aircraft need to be adjusted, it is better to fully consider the classes of the aircraft rather than the aircraft themselves. }\end{remark}

\subsection{Algorithm 2}

In the following, we propose the following algorithm to deal with the scheduling problem of landing and takeoff aircraft on dual runways with spacing no larger than $760$ $m$.

\textbf{Algorithm 2.} Suppose that there are totally $n+m$ aircraft composed of $n$ landing aircraft and $m$ takeoff aircraft, denoted by $Tcf_{1}, Tcf_{2}, \cdots, Tcf_{n+m}$, and there is at least a feasible sequence for them to land/takeoff without conflicts. Let $Tcf^0_1, Tcf^0_2, \cdots, Tcf^0_n$ denote all the landing aircraft, and $Tcf^1_{1}, Tcf^1_{2}, \cdots, Tcf^1_{m}$ denote all the takeoff aircraft.  The main idea of this algorithm is also to exclude most of the non-optimal solutions and narrow down the search for the optimal solutions within a small range. The main steps are as follows.

Step 1. Use Algorithm 1 to find the optimal sequence for the landing aircraft, denoted by $\Phi_a$, and the optimal sequence for the takeoff aircraft, denoted by $\Phi_b$. Without loss of generality, we suppose that $F(\Phi_a)\geq F(\Phi_b)$ and the case of $F(\Phi_a)< F(\Phi_b)$ can be similarly discussed.

Based on the obtained optimal sequence $\Phi_a$, by a series of aircraft order adjustments, record all possible optimal sequences based on the bases of all the breakpoint-drift equivalent transformation sets and represent their formed set as $\Lambda_{a}$ for aircraft $Tcf^0_{1},Tcf^0_{2}, \cdots, Tcf^0_{n}$.

The goal of recording all possible optimal sequences is to make preparations to generate all possible block subsequences in Steps 3 and 4. Therefore, we can reduce some computation amount according to the needs of Steps 3 and 4.

Step 2. (2.1) According to the landing/takeoff time increments of block subsequences and the number of breakpoint aircraft, classify all possible block subsequences into several sets, which can be processed offline.

When the number of breakpoint aircraft in a block subsequence is fixed, all possible block subsequences can be obtained by an enumeration method. It should be noted that the block subsequences with more breakpoint aircraft can be generated based on the block subsequences with less breakpoint aircraft. It should be also noted that
the number of all possible block subsequences of different class combinations and the range of time lengths and the landing/takeoff increments of block subsequences are both not large under the minimum separation time standards at Heathrow Airport and for the RECAT-EU system are considered without considering resident-point aircraft, when the number of breakpoint aircraft in a block subsequence is fixed.

(2.2) Focus on combinations of different block subsequences, and study the landing/takeoff time increments of any two given consecutive block subsequences, where the last two aircraft of the first block subsequence is the first two aircraft of the second block subsequence. According to the landing/takeoff time increments and the number of breakpoint aircraft, classify all possible two consecutive subsequences into several sets.

In our simulations for the case of $F(\Phi_a)\geq F(\Phi_b)$, the combination of a D-block subsequence and a T-block subsequence is frequently used, which has the form of $\langle DB, TB\rangle$ where $DB$ denotes a D-block subsequence and $TB$ denotes a T-block subsequence, the first two aircraft of $DB$ have the same landing/takeoff time, and the last two aircraft of $TB$ have the same landing/takeoff time (see Tables \ref{tab3} and \ref{tab1} for related discussions). 
The time length and the landing time increment of D-block subsequence in such a combination are key quantities that needs to be studied for the optimal sequence.

Step 3. (3.1) Match the takeoff aircraft with a landing sequence in $\Lambda_{a}$ to generate a new sequence, denoted by $\Phi_c=\langle Tcf_{c1}, Tcf_{c2}, \cdots, Tcf_{c(n+m)}\rangle$, and minimize $\Phi_c$ under the condition that $F(\Phi^c_{LA}, Sr_{LA}(\Phi_c))=F(\Phi_a)$, where $\Phi^c_{LA}$ and $Sr_{LA}(\Phi_c)$ denote the landing sequence and the corresponding landing time vector in $\Phi_c$.

Under Assumption \ref{ass4.21}, the sequence $\Phi_c$ can be divided into a series of block subsequences. Since $\!F(\Phi^c_{LA},\! Sr_{LA}(\Phi_c))$ $=F(\Phi_a)$, we can use a dynamic programming approach to find the optimal sequence $\Phi_c$ based on the sets defined in Step 2. It should be noted that at each step of the dynamic programming approach, we need consider all the possible block subsequences that can decrease the objective function $F(\cdot)$ of the whole sequence, and if the generation of one block subsequence might disrupt the existing block subsequences, we can first generate the block subsequence and then use it as a basis to generate other block subsequences in an optimal way.

(3.2) {Let $n_z$ denote the number of the aircraft of $\Phi_c$ which take off in $(t_0+F(\Phi_a),+\infty)$.} Insert the last $n_z$ aircraft of $\Phi_c$ into the subsequence of $\Phi_c$, $\Phi^{sub}_{c}=\langle Tcf_{1}, Tcf_{2}, \cdots, Tcf_{n+m-n_z}\rangle$, to generate a new sequence $\Phi_{c3}$ and minimize $F(\Phi_{c3})$ without changing the orders of the landing aircraft based on the sets defined Step 2 by combing the capacity, the time length, the initial redundant time and the landing/takeoff time increment of each block subsequence. Note here that the sequence $\Phi_{c3}$ can contain D-block subsequences.

It should be noted that the minimum value of the objective function $F(\cdot)$ for the optimization problem (\ref{optim1}) lies in $[F(\Phi_a), F(\Phi_{c3})]$ and the final optimal sequence usually contains at least $n_z-1$ D-block subsequences, which can be used as an initial condition to find a feasible solution near the optimal solutions to narrow the search range in step 4.

{Step 4. The main idea of this step is to implement a full space search for an optimal solution $(\Phi_{c4}, Sr(\Phi_{c4}))$ for the optimization problem (\ref{optim1}) under the condition that $F(\Phi_{a})\leq F(\Phi_{c4})\leq F(\Phi_{c3})$ starting from the initial condition that $\Phi^{c4}_{LA}\in \Lambda_{a}$ base on T-block and D-block subsequences.

Since different D-block subsequences might have different landing time increments, we can search for the optimal solutions for the optimization problem (\ref{optim1}) by increasing or decreasing the number and the landing time increments of D-block subsequences of the sequence $\Phi_{c4}$.
Based on this idea, since $F(\Phi_{a})\leq F(\Phi^{c4}_{LA})\leq F(\Phi_{c4})$, the condition $F(\Phi_{a})\leq F(\Phi_{c4})\leq F(\Phi_{c3})$ can be relaxed to be $F(\Phi_{a})\leq F(\Phi^{c4}_{LA})\leq F(\Phi_{c3})$.}

When optimizing the sequence, we can first consider the block subsequences containing the breakpoint aircraft since the number of all possible block subsequences containing breakpoint aircraft is usually larger than the number of all possible block subsequences containing no breakpoint aircraft without considering the aircraft classes.

{When D-block subsequences of $\Phi_{c4}$ are considered, some necessary aircraft order adjustments might be needed to be imposed on the subsequence $\Phi^{c4}_{LA}$} to fully consider all the possibilities of the block subsequences to reduce the value of the objective function $F(\cdot)$. For example, suppose that $\phi_{LA}^{0}=\langle Tcf_1, Tcf_2, Tcf_3, Tcf_4, Tcf_5\rangle$ is a landing subsequence of a sequence $\phi_0$, and $\phi_0$ is composed of a D-Block subsequence $DB_1$ and a T-Block subsequence $TB_1$ with $Tcf_1, Tcf_2\in DB_1$ and $Tcf_2, Tcf_3, Tcf_4, Tcf_5\in TB_1$. According to the definition of D-block subsequence, $S_{12}>Y_{12}$. If $Y_{15}>Y_{12}$ and the objective function value of the landing sequence can be decreased by exchanging the orders of $Tcf_2$ and $Tcf_5$, the D-block subsequences of the sequence $\phi=\langle Tcf_1, Tcf_5, Tcf_3, Tcf_4, Tcf_2\rangle$ are needed to be considered.

There are several method, e.g., binary search method, to realize the search in this step, and the most basic method is to generate the optimal sequence completely according to block subsequences and the combinations of block subsequences discussed in Step 2, which might have more computation amount. To reduce the computation amount, we can use the number of breakpoint aircraft to partition the search space as stated in Algorithm 1. 

\begin{remark}{\rm {The main idea of Algorithm 2 is to decompose the aircraft sequence into several block subsequences and fully explore combinations of the block subsequences starting from an optimal landing/takeoff subsequence to consider the optimization problem (\ref{optim1}).} Note that the number of all possible block subsequences of different class combinations is usually not large under the minimum separation time standards at Heathrow Airport and for the RECAT-EU system without considering resident-point aircraft. The use of block subsequences might significantly reduce the computation amount of the algorithm.}\end{remark}

\begin{remark}{\rm The complexities of Algorithms 1 and 2 are heavily related to the aircraft number, the aircraft classes, and the constraints of the aircraft time windows.
But the algorithms essentially exhibit the features of polynomial algorithms and can be applied in actual systems in real time, and based on these two algorithms, the optimal solutions for the optimization problem (\ref{optim1}) can be obtained.}\end{remark}

\section{Simulation}\label{simulations}
In this section, we evaluate the efficiency of the proposed algorithms by comparing its computation time and objective function performance against a standard MIP solver for two scheduling problems: takeoff and landing operations on a single runway and on dual runways. For the single-runway case, we consider four scenarios with $|A|= 30,40,50,60$ aircraft. For the dual-runway case, we consider scenarios with $|A| = 70, 80, 90, 100$ aircraft. The aircraft are categorized into $6$ classes based on the RECAT-EU framework $(A, B, C, D, E, F)$, with class proportions of $10\%$, $20\%$, $25\%$, $15\%$, $20\%$, and $10\%$, respectively.
The minimum separation times for aircraft pairs with the same operation task are detailed in Tables \ref{tab:Heathrow_separation} and \ref{tab:Taking-off_separation}. For both scheduling problems, the minimum separation times between a takeoff aircraft and a trailing landing aircraft are set to $D_T = D_P = 60$ seconds, while the minimum separation time between a landing aircraft and a trailing takeoff aircraft is set to $T_D = 75$ seconds for the single-runway case and $P_D = 0$ seconds for the dual-runway case.
All simulations are conducted on a computer with an AMD Ryzen 7 7840H processor (3.8 GHz, 16 GB RAM). The MIP formulations are solved using CPLEX 12.10, with a computation time limit of $600$ seconds per instance. The proposed algorithms are implemented in MATLAB R2021b.

First, we consider the aircraft scheduling problem on a single runway. The earliest landing or takeoff times for each aircraft are generated randomly. These times follow a uniform distribution within the interval $[0, T_E]$ minutes. The time window lengths are set to $T_W = 30$, $45$ and $60$ minutes to evaluate the algorithms's performance under different levels of scheduling flexibility. We evaluate three types of operations: takeoff only, landing only, and mixed operations (both takeoff and landing). The values of $T_E$ are set to $30$ minutes for the takeoff-only and landing-only cases and $20$ minutes for the mixed operation case.
Table~\ref{tab:cost_function CPU Time2} presents a comparison of the objective function values and computation times between Algorithm 1 and the MIP solver, along with the percentage gap in objective values between the two methods. The values of the objective function are denoted as NFS when the MIP solver fails to find a solution within the $10$-minute time limit.

Next, we consider the aircraft scheduling problem on dual runways. The earliest landing or takeoff times for each aircraft are also generated randomly, following a uniform distribution within the interval $[0,T_E]$ minutes. The time windows are again set to $T_W = 30$, $45$ and $60$ minutes.
Table~\ref{tab:cost_function CPU Time3} presents a comparison of the objective function values and computation times between Algorithm 2 and the MIP solver, along with the percentage gap in objective values between the two approaches.

\begin{table}[!htb]
		\caption{Minimum landing separation times in Heathrow Airport (Sec)}
		\label{tab:Heathrow_separation}
		\centering
		\begin{tabular}{ccccccccccccc}
			\toprule
			\midrule
            \multicolumn{2}{c}{} & \multicolumn{6}{c}{Trailing Aircraft} \\
            \midrule 
            \multicolumn{2}{c}{} & A & B & C & D & E & F \\
            \midrule 
            \multirow{6}*{Leading Aircraft}  & A & 90 & 135& 158& 158 & 158 & 180 \\
                                 & B & 90 & 90 & 113& 113 & 135 & 158 \\
                                 & C & 60 & 60 & 68 & 90  & 90  & 135 \\
                                 & D & 60 & 60 & 60 & 60  & 68  & 113 \\
                                 & E & 60 & 60 & 60 & 60  & 68  & 90 \\
                                 & F & 60 & 60 & 60 & 60  & 60  & 60 \\

			\midrule
			\bottomrule
		\end{tabular}
\end{table}
\begin{table}[!htb]
		\caption{Minimum takeoff separation times based on RECAT-EU (Sec)}
		\label{tab:Taking-off_separation}
		\centering
		\begin{tabular}{ccccccccccccc}
			\toprule
			\midrule
            \multicolumn{2}{c}{} & \multicolumn{6}{c}{Trailing Aircraft} \\
            \midrule 
            \multicolumn{2}{c}{} & A & B & C & D & E & F \\
            \midrule 
            \multirow{6}*{Leading Aircraft}  & A & 80 & 100& 120& 140 & 160 & 180 \\
                                 & B & 80 & 80 & 100& 100 & 120 & 140 \\
                                 & C & 60 & 60 & 80 & 80  & 100 & 120 \\
                                 & D & 60 & 60 & 60 & 60  & 60  & 120 \\
                                 & E & 60 & 60 & 60 & 60  & 60  & 100 \\
                                 & F & 60 & 60 & 60 & 60  & 60  & 80 \\

			\midrule
			\bottomrule
		\end{tabular}
\end{table}

\begin{table*}[!htb]
        \scriptsize
		\caption{{{Comparison of performance and computation times for single-runway aircraft scheduling problem}}}
		\label{tab:cost_function CPU Time2}
		\centering
		\begin{tabular}{ccccccccccc}
			\toprule
			\midrule
             \multirow{2}*{$T_W$(min)} & Aircraft Number  & \multirow{2}*{$T_E$(min)} & Operation & \multicolumn{2}{c}{Objective Function(s)} & \multicolumn{2}{c}{Computation Times(s)} & Gap\\
             & $|A|$  &  &  Task & MIP & Our Algorithm & MIP & Our Algorithm & \%\\
            \midrule 
            \multirow{12}*{30} & \multirow{3}*{30} & 30 & Takeoff & 2162 & 2162 & 600 & 0.25 & 0 \\
                              & & 30 & Landing & 2235 & 2221 & 600 & 0.15 & 0.63 \\
                             &  & 20 & Mixed & 1856 & 1856 & 600 & 0.96 & 0\\
            \cmidrule{2-9}
             & \multirow{3}*{40} & 30 & Takeoff & 2971 & 2951 & 600 & 0.33 & 0.68  \\
                              & & 30 & Landing & 3061 & 3015 & 600 & 0.21 & 1.53  \\
                            &   & 20 & Mixed & 2601 & 2540 & 600 & 1.79 & 2.40\\
            \cmidrule{2-9}
             &\multirow{3}*{50} & 40 & Takeoff & NFS & 3709 & 600 & 0.52 & / \\
                            &   & 40 & Landing & NFS  & 3833 & 600 & 0.43 & / \\
                            &   & 30 & Mixed & 3183 & 3095 & 600 & 2.72 & 2.84\\
            \cmidrule{2-9}
             &\multirow{3}*{60} & 50 & Takeoff & NFS & 4461 & 600 & 0.83 & / \\
                           &    & 50 & Landing & NFS & 4586 & 600 & 1.11 & / \\
                           &    & 40 & Mixed   & NFS & 3737 & 600 & 4.94 & / \\
			\midrule
            \multirow{12}*{45} & \multirow{3}*{30} & 20 & Takeoff & 2213 & 2193 & 600 & 0.03 & 0.91 \\
                              & & 20 & Landing & 2217 & 2188 & 600 & 0.09 & 1.32 \\
                              & & 20 & Mixed & 1975 & 1975 & 600 & 0.95 & 0\\
            \cmidrule{2-9}
             & \multirow{3}*{40}& 20 & Takeoff & 2997 & 2886 & 600 & 0.09 & 3.84  \\
                              & & 20 & Landing & 3084 & 2875 & 600 & 0.08 & 7.27  \\
                            &   & 20 & Mixed & 2502 & 2487 & 600 & 1.4 & 0.60\\
            \cmidrule{2-9}
             &\multirow{3}*{50} & 30 & Takeoff & 3696 & 3656 & 600 & 0.34 & 1.09 \\
                            &   & 30 & Landing & 3892 & 3736 & 600 & 0.51 & 4.17 \\
                            &   & 30 & Mixed   & 3217 & 3164 & 600 & 2.25 & 1.68\\
            \cmidrule{2-9}
             &\multirow{3}*{60} & 40 & Takeoff & 4548 & 4532 & 600 & 0.61 & 0.35 \\
                           &    & 40 & Landing & 4686 & 4476 & 600 & 0.29 & 4.69 \\
                           &    & 30 & Mixed   & 4129 & 3780 & 600 & 3.65 & 9.17 \\
            \midrule
            \multirow{12}*{60} & \multirow{3}*{30} & 30 & Takeoff & 2220 & 2220 & 600 & 0.18 & 0 \\
                            &   & 30 & Landing & 2253 & 2245 & 600 & 0.15 & 0.36\\
                            &   & 20 & Mixed & 1931 & 1931 & 600 & 1.21 & 0\\
            \cmidrule{2-9}
             &\multirow{3}*{40} & 30 & Takeoff & 3134 & 3074 & 600 & 0.13 & 1.95  \\
                            &   & 30 & Landing & 3330 & 3085 & 600 & 0.23 & 7.94  \\
                            &   & 20 & Mixed & 2564 & 2534 & 600 & 2.06 & 1.18\\
            \cmidrule{2-9}
             &\multirow{3}*{50} & 30 & Takeoff & 3684 & 3624 & 600 & 0.09 & 1.66  \\
                            &   & 30 & Landing & 3831  & 3642 & 600 & 0.13 & 5.19 \\
                           &    & 20 & Mixed & 3220 & 3141 & 600 & 2.49 & 2.52\\
            \cmidrule{2-9}
             &\multirow{3}*{60} & 30 & Takeoff & 4573 & 4373 & 600 & 0.55 & 4.57 \\
                         &      & 30 & Landing & 4620 & 4367 & 600 & 0.42 & 5.79 \\
                          &     & 20 & Mixed & 3890 & 3709 & 600 & 3.04 & 4.88\\
			\bottomrule
		\end{tabular}
\end{table*}

\begin{table*}[!htb]
		\caption{{Comparison of performance and computation times for dual-runway aircraft scheduling problem}}
		\label{tab:cost_function CPU Time3}
		\centering
		\begin{tabular}{ccccccccc}
			\toprule
			\midrule
            $T_W$ & $T_E$ & Aircraft Number  & \multicolumn{2}{c}{Objective Function(s)} & \multicolumn{2}{c}{Computation Times(s)} & Gap\\
            (min)& (min) & $|A|$   & MIP & Our Algorithm & MIP & Our Algorithm & \%\\
            \midrule 
            \multirow{4}*{30} & 30& 70  & 2986 & 2775 & 600 & 1.28 & 7.60 \\
             &40& 80  & 3329 & 3196 & 600 & 1.75 & 4.16  \\
             &45& 90  & 4193 & 3670 & 600 & 2.14 & 14.25  \\
             &60& 100 & 4533 & 4241 & 600 & 3.03 & 6.89  \\
			\midrule
            \multirow{4}*{45} & 20& 70  & 2952 & 2738 & 600 & 1.23 & 7.82 \\
             &30& 80  & 3224 & 3020 & 600 & 1.67 & 6.75  \\
             &40& 90  & 3701 & 3589 & 600 & 3.05 & 3.12  \\
             &50& 100 & 4259 & 3947 & 600 & 4.36 & 7.90  \\
            \midrule
            \multirow{4}*{60} & 20& 70  & 2814 & 2668 & 600 & 1.90 & 5.47 \\
             &20& 80  & 3562 & 3151 & 600 & 1.38 & 13.04  \\
             &20& 90  & 4111 & 3528 & 600 & 3.33 & 16.52  \\
             &20& 100 & 4505 & 3873 & 600 & 3.54 & 16.32  \\
			\bottomrule
		\end{tabular}
\end{table*}

From Tables~\ref{tab:cost_function CPU Time2} and \ref{tab:cost_function CPU Time3}, it is evident that the proposed algorithms consistently obtains solutions within $5$ seconds, whereas the MIP solver requires the full $10$-minute time limit. Moreover, our algorithmx consistently yield better objective function values than those produced by the MIP solver, with the performance gap widening as the number of aircraft increases. In addition, the runtime of the MIP solver has been extended to over one hour but no significant improvements were found in the objective function values compared to the results obtained within the $10$-minute time limit.

In the single-runway scenario, our algorithm shows particularly strong performance in the takeoff-only and landing-only cases. The performance advantage is less pronounced in the mixed takeoff-and-landing case. This is primarily due to the relatively short separation times required between different operation tasks: $60$ seconds between a takeoff aircraft and a trailing landing aircraft, and $75$ seconds between  a landing aircraft and a trailing takeoff aircraft. These moderate separation values enable the MIP solver to generate relative high-quality solutions by frequently using takeoff-landing and landing-takeoff transitions. Nonetheless, as the number of aircraft increases, the superiority of our algorithm remains increasingly evident.

For the dual-runway scenario, we evaluate $12$ instances with progressively larger numbers of aircraft, reflecting the higher operational capacity of dual-runway airports. In particular, when the problem size reaches $90$ and $100$ aircraft, the performance gap in objective function values between our algorithm and the MIP solver exceeds $16\%$, exhibiting the scalability and efficiency of our proposed algorithms.

\section{Conclusions}
In this paper, scheduling problems of landing and takeoff aircraft on a same runway and on dual runways were addressed. A new theoretical framework for scheduling problem of aircraft was established, which is completely different from the framework of mixed-integer optimization problem.
Two real-time optimal algorithms are proposed for the four scheduling problems by fully exploiting the combinations of different classes of aircraft, which can even be applied to the RECAT-EU systems. Numerical examples are presented to show the effectiveness of the algorithms. In particular, when $100$ aircraft are considered, by using the algorithm in this paper, the optimal solution can be obtained in less than $5$ seconds, while by using the CPLEX software to solve the mix-integer optimization model, the optimal solution cannot be obtained within $1$ hour.


\end{document}